\documentclass[aos, preprint]{imsart}

\RequirePackage[OT1]{fontenc}
\RequirePackage{amsthm}
\RequirePackage[colorlinks,citecolor=blue,urlcolor=blue]{hyperref}
\arxiv{arXiv:1708.03964}

\usepackage{xr} % for cross-referencing (Nestor)
\usepackage{lscape,exscale,german,amsthm,enumerate,rawfonts,epsfig,multirow,latexsym,dsfont,color}
\usepackage[intlimits]{amsmath}
\usepackage[cp850]{inputenc}
\usepackage{bbm}
\usepackage{footnote}
\usepackage{float}

\usepackage{latexsym}
\usepackage{mathtools}
\usepackage{amsmath}
\usepackage{amssymb}

\usepackage{url}

\usepackage[authoryear]{natbib}
\startlocaldefs

\makeatletter
\providecommand*{\diff}%
{\@ifnextchar^{\DIfF}{\DIfF^{}}}
\def\DIfF^#1{%
\mathop{\mathrm{\mathstrut d}}%
\nolimits^{#1}\gobblespace}
\def\gobblespace{%
\futurelet\diffarg\opspace}
\def\opspace{%
\let\DiffSpace\!%
\ifx\diffarg(%
\let\DiffSpace\relax
\else
\ifx\diffarg[%
\let\DiffSpace\relax
\else
\ifx\diffarg\{%
\let\DiffSpace\relax
\fi\fi\fi\DiffSpace}

\numberwithin{equation}{section}
\theoremstyle{plain}

\newcommand{\er}{\mathbb R}
\newcommand{\C}{\mathbb C}
\newcommand{\bSigma}{\boldsymbol{\Sigma}}

\newcommand{\bOmega}{\boldsymbol{\Omega}}

\newcommand{\bY}{\mathbf{Y}}
\newcommand{\bS}{\mathbf{S}}

\newcommand{\bX}{\mathbf{X}}
\newcommand{\bx}{\mathbf{x}}

\newcommand{\bW}{\mathbf{W}}
\newcommand{\bZ}{\mathbf{Z}}
\newcommand{\bU}{\mathbf{U}}
\newcommand{\bT}{\mathbf{T}}
\newcommand{\bP}{\mathbf{P}}

\newcommand{\bR}{\mathbf{R}}

\newcommand{\bzero}{\mathbf{0}}
\newcommand{\bI}{\mathbf{I}}

\newcommand{\Cov}{\mbox{Cov}}
\newcommand{\E}{\mbox{E}}
\newcommand{\Var}{\mbox{Var}}

\newcommand{\bF}{\mathbf{F}}

\newcommand{\bD}{\mathbf{D}}

\newcommand{\bM}{\mathbf{M}}

\usepackage{ulem}

\newtheorem{theorem}{Theorem}

\newtheorem{proposition}{Proposition}
\newtheorem{lemma}{Lemma}

\newtheorem{remark}{Remark}

\endlocaldefs

\begin{document}

\begin{frontmatter}
\title{Testing for Independence of Large Dimensional Vectors}
\runtitle{Testing for Independence of Large Vectors}
%\thankstext{T1}{Footnote to the title with the ``thankstext'' command.}

\begin{aug}
  \author{\fnms{Taras} \snm{Bodnar}\thanksref{m1}\ead[label=e1]{taras.bodnar@math.su.se}},
  \author{\fnms{Holger} \snm{Dette}\thanksref{t1,m2}\ead[label=e2]{holger.dette@ruhr-uni-bochum.de}}
\and
\author{\fnms{Nestor} \snm{Parolya}\thanksref{m3}\ead[label=e3]{parolya@statistik.uni-hannover.de}}

%\thankstext{t1}{Some comment}
\thankstext{t1}{DFG Research Unit 1735, DE 502/26-2}
%\thankstext{t3}{Second supporter of the project}
\runauthor{T. Bodnar, H. Dette and N. Parolya}

\affiliation{Stockholm University\thanksmark{m1}, Ruhr University Bochum\thanksmark{m2} and Leibniz University Hannover\thanksmark{m3}}

\address{Address of the First author\\
 % \thanksmark{m1}
\footnotesize  Department of Mathematics, Stockholm University, SE-10691 Stockhom, Sweden\\
\printead{e1}\vspace{0.2cm}}

\address{Address of the Second author\\
  %\thanksmark{m2}
\scriptsize Department of Mathematics,  Ruhr University Bochum, D-44870 Bochum, Germany\\
\printead{e2}\vspace{0.2cm}}

\address{Address of the Third author\\
%\thanksmark{m3}
  \scriptsize Institute of Statistics, Leibniz University Hannover, D-30167 Hannover, Germany\\
\printead{e3}}

\end{aug}

\begin{abstract}
\quad  In this  paper new tests for the  independence of two high-dimensional vectors are investigated.
We consider the case where the dimension of the vectors increases with the sample size
and propose  multivariate analysis of variance-type statistics for  the hypothesis of a block diagonal covariance matrix.
The asymptotic properties of the new test statistics are investigated under the null hypothesis and the alternative hypothesis using random matrix theory. For this purpose we study  the weak convergence of linear spectral statistics of central and (conditionally) non-central Fisher matrices.
In particular, a central limit theorem for linear spectral statistics of large dimensional (conditionally) non-central Fisher matrices is derived which is then used to analyse
the power of the tests under the alternative.

The theoretical results are illustrated  by means of a simulation study where we also compare the new tests with several alternative, in particular
with the commonly used corrected likelihood ratio
 test. It is demonstrated that the latter test does not keep its nominal level, if the dimension of one sub-vector is relatively small compared to the dimension of the other sub-vector. On the other hand the tests proposed in this paper provide a reasonable approximation of the nominal level in such situations. Moreover, we observe that one of the proposed tests is most powerful under  a variety of correlation scenarios.
\end{abstract}



\end{frontmatter}
\section{Introduction} \label{sec1}
\def\theequation{1.\arabic{equation}}
\setcounter{equation}{0}

Estimation and testing the structure of the covariance matrix are important problems that have a number of applications in practice. For instance, the covariance matrix plays an important role in the determination of the optimal portfolio structure following the well-known mean-variance analysis of \cite{markowitz1952}. It is also used in  prediction theory where the problem of forecasting future values of the process based on its previous observations arises.
In such applications misspecification of  the covariance matrix might lead to significant errors in the optimal portfolio structure and predictions. The problem becomes even more difficult if  the  dimension is of similar order or even larger as  the sample size.
A number of such situations are present in biostatistics, wireless communications and finance (see, e.g., \cite{Fan2006}, \cite{Johnstone} and references therein).

The sample covariance matrix is the commonly used  estimator in practice.
However, in the case of large dimension (compared to  the sample size), a number of studies demonstrate that the sample covariance does not perform well
as an estimator of the population covariance matrix and
numerous authors have recently addressed this problem. One approach is based on the construction of improved estimators
in particular shrinkage type estimators  which
 reduce the variability of  the sample covariance matrix at the cost of an additional bias (see, \cite{ledoitwolf2012}, \cite{wangpantongzhu2015}  or \cite{bodnaretal2014,bodnaretal2016} among others). Alternatively several authors impose structural assumptions
 on the population covariance matrix such  as a block diagonal structure (e.g., \cite{devigall2016}),  Toeplitz  matrix  (see, \cite{cairenzhou2013}),
 band matrix (see, \cite{bicklevi2008}) or general sparsity assumptions  (see \cite{cailuiluo2011}, \cite{caishen2011},  \cite{caizhou2012} among others) and show, that the population covariance matrix can be estimated consistently in these cases, even for large dimensions.
However, these techniques may fail if the structural assumptions are not satisfied and
consequently it is desirable to validate  the corresponding assumptions regarding  the postulated
  structure of the covariance matrix.

In the present paper we consider the problem of testing for a  block diagonal structure of the covariance, which  has found considerable interest in the literature.  Early work in this direction  has been done by \cite{mauchly1940}, who proposed a likelihood ratio test for the hypothesis of sphericity of a normal distribution, that is the independence of all components. This method has been extended by \cite{guptaxu2006} to the non-normal case  and by
 \cite{baietal2009} and \cite{jiangyang2013} to the high-dimensional case. An alternative approach is based on the empirical distance between the sample covariance matrix and the target (e.g., a multiplicity of the identity matrix)
  and was  initially suggested by \cite{john1971}  and \cite{nagao1973}. These tests can also be
  extended for testing the  corresponding hypotheses in the  high-dimensional
  setup (see, \cite{ledoitwolf2002}, \cite{birkedette2005}, \cite{fishsungall2010}, \cite{chenzhangzhong2010}). Other authors use the distributional properties of the largest eigenvalue of the sample   covariance matrix to construct tests
  (see \cite{johnstone2001,johnstone2008} for example).

In the problem of testing  the independence between two (or more)  groups of random variables  under the assumption of normality
the  likelihood ratio approach  has also found considerable interest in the literature.  The main results for a fixed dimension
can be found in the text books of \cite{muirhead1982} and \cite{anderson2003}. Recently,  \cite{jiangyang2013} have extended the
 likelihood ratio approach to the case of high-dimensional data, while \cite{hyodoetal2015}  and \cite{yamadaetal2017}
 used an empirical distance approach  to test for a block diagonal covariance matrix.

 In Section \ref{sec2} we introduce the testing problem (in the case of two blocks)  and
demonstrate by means of a small simulation study   that the likelihood ratio test does not yield a reliable approximation of the nominal level,
if the size of one block is small compared to the other  one. In Section \ref{sec3} we introduce three alternative test statistics which are motivated from classical
   multivariate analysis of variance (MANOVA) and are defined as linear spectral statistics of a
  Fisher matrix. We derive their asymptotic distributions under  the null hypotheses
  and illustrate the approximation of the nominal level by means of a simulation study. A comparison  with the commonly used likelihood ratio test
  shows  that the new tests  provide a reasonable approximation of the nominal level in  situations where the likelihood ratio test
 fails. Section \ref{sec4} is devoted to the analysis of statistical properties of the new tests under the alternative hypothesis.
 For this purpose, we present a new central limit theorem for a (conditionally) non-central Fisher random matrix which is of own interest and
 can be  used  to study some properties of the power of the new tests.  Finally, most technical details and proofs are given in the appendix (see, Section \ref{sec5}) and in  the supplementary material (see, \cite{supplement}).

\section{Testing for independence} \label{sec2}
\def\theequation{2.\arabic{equation}}
\setcounter{equation}{0}

Let $\bx_1,...,\bx_n$ be a sample of i.i.d. observations from  a  $p$-dimensional normal distribution with zero mean vector and covariance
matrix $\bSigma$, i.e.
$\bx_1\sim\mathcal{N}_p(\bzero,\bSigma)$. We define the $p \times n$ dimensional observation matrix $\bX=(\bx_1,...,\bx_n)$ and denote by
\begin{equation} \nonumber 
\bS=\frac{1}{n}\bX\bX^\top
\end{equation}
 the sample covariance matrix which is used as an estimate of $\bSigma$.
 It is well known that  $n\bS$ has a $p$-dimensional Wishart distribution with $n$ degrees of freedom and covariance matrix $\bSigma$, i.e., $n\bS \sim W_p(n,\bSigma)$.
In the following we consider partitions of the population and the sample covariance matrix given by
\begin{eqnarray}\label{SigmaS}
\bSigma=\left(
          \begin{array}{cc}
            \bSigma_{11} & \bSigma_{12} \\
            \bSigma_{21} & \bSigma_{22}
         \end{array}
        \right)
\quad\text{and} \quad
n\bS=\left(
          \begin{array}{cc}
            \bS_{11} & \bS_{12}\\
            \bS_{21} & \bS_{22}
          \end{array}
        \right) \,,
\end{eqnarray}
respectively,
where $\bSigma_{ij} \in \er^{p_i\times p_j}$ and $\bS_{ij} \in \er^{ p_i\times p_j}$ with $i,j=\{1,2\}$ and $p_1+p_2=p$.
We are interested  in the hypothesis that the sub-vectors $\bx_{1,1}$
and $\bx_{1,2}$ of size $p_1$ and $p_2$ in the vector $\bx_1 =(\bx_{1,1}^\top,\bx_{1,2}^\top)^\top$
 are independent,
or equivalently that the covariance matrix is block diagonal, i.e.
\begin{equation}\label{hyp_H0}
H_0:~\bSigma_{12}=\mathbf{O} \quad \mbox{ versus } \quad
H_1:~\bSigma_{12}\neq \mathbf{O}\,.
\end{equation}
Here the symbol $\mathbf{O}$ denotes a matrix of an appropriate order with all entries equal to $0$.
It is worthwhile to mention that the case of non-zero mean vector can be treated exactly in the same way
observing  that the centred  sample covariance matrix, has a  $ \tfrac{1}{n-1 }W_p(n-1,\bSigma)$ distribution.
Thus, one needs to normalize the sample covariance matrix by $1/(n-1)$ instead of $1/n$ due to the {\it substitution principle} of \citet{zhengbaiyao2015b} and the results presented in our paper will still remain valid.

Throughout this paper we consider the case where the dimension of the blocks is increasing with the sample size, that is
$p=p(n)$, $p_i=p_i(n)$,  such that
\begin{equation}\label{ass1}  \nonumber
\lim_{n \to\infty}
\frac{p_i}{n} = c_i < 1~,~i=1,2
\end{equation}
and define $c= c_1+c_2$.  For further reference we also introduce the quantities
\begin{eqnarray}
\label{g1}
 \gamma_{1,n}& =& \frac{p-p_1}{p_1},\\
\label{g2}   \gamma_{2,n} &=& \frac{p-p_1}{n-p_1},\\[2mm]
\label{g3}   h_n &= & \sqrt{\gamma_{1,n}+\gamma_{2,n}-\gamma_{1,n}\gamma_{2,n}}\,.
\end{eqnarray}
A common approach in testing for independence is the likelihood ratio test based on the statistic
$$
  V_n = \dfrac{|\bS|}{|\bS_{11}||\bS_{22}|} = \frac{|\bS_{11}|\left|\bS_{22} - \bS_{21}\bS_{11}^{-1}\bS_{12}  \right|}{|\bS_{11}||\bS_{22}|} = \left|\bI_{p-p_1} - \bS_{21}\bS_{11}^{-1}\bS_{12}\bS^{-1}_{22}  \right|\,.
$$
 The null hypothesis is rejected for small values of $V_n$. \cite{jiangetal2013}  showed  that under the assumptions made in this section
 $V_n$ can be written in terms of a determinant of  a central Fisher matrix, that is
\begin{eqnarray}\label{Vn}
 V_n = \Big|\bI_{p-p_1}-\bF(\bF+\frac{\gamma_{1,n}}{\gamma_{2,n}}\bI_{p-p_1})^{-1} \Big|= 
 \Big|\frac{\gamma_{2,n}}{\gamma_{1,n}} \bF+ \bI_{p-p_1}\Big|^{-1} ~,
\end{eqnarray}
where $\bF=\frac{1}{p_1}\bS_{21}\bS^{-1}_{11}\bS_{12}\big ( \frac{1}{n-p_1}(\bS_{22}-  \bS_{21}\bS^{-1}_{11}\bS_{12})\big )^{-1}$.
Under the null hypothesis of independent blocks, the matrix $\bF$
 is a ''ratio'' of two central Wishart matrices with $p_1$ and $n-p_1$ degrees of freedom. Naturally, it is called a central Fisher matrix with $p_1$ and $n-p_1$ degrees of freedom, an analogue to its one dimensional counterpart  (see, \cite{Fisher1939}).
 In particular, we have the following result (see, Theorem 8.2 in \cite{yaobaizheng2015})
\begin{proposition} \label{prop1}
Under the  null hypothesis we have for $T_{LR}=\log(V_n)$
  \begin{eqnarray*}
\frac{T_{LR} -(p-p_1)s_{LR} -\mu_{LR}}{\sigma_{LR}} \overset{\mathcal{D}}{\longrightarrow} \mathcal{N}(0, 1 )\,,
  \end{eqnarray*}
where the quantities $\mu_{LR}, \sigma_{LR}^2$ and $s_{LR}$ are defined by
{\small
\begin{eqnarray*}
     \mu_{LR} &=& 1/2\log\left[ \frac{(w_n^{*\;2}-d_n^{*\;2})h_n^2}{(w_n^{*}h_n-\gamma_{2,n}d_n^{*\;2})^2}\right],\qquad\qquad \sigma^2_{LR} =2\log\left[\frac{w_n^{*\;2}}{w_n^{*\;2}-d_n^{*\;2}}\right],\\
  s_{LR} &=& \log\Big(\frac{\gamma_{1,n}}{\gamma_{2,n}}(1-\gamma_{2,n})^2 \Big)+ \frac{1-\gamma_{2,n}}{\gamma_{2,n}}\log(w_n^*) - \frac{\gamma_{1,n}+\gamma_{2,n}}{\gamma_{1,n}\gamma_{2,n}}\log(w_n^*-d_n^*\gamma_2/h_n)\\
  &+&
\left\{
\begin{array}{ll}
\frac{1-\gamma_{1,n}}{\gamma_{1,n}}\log(w_n^*-d_n^*h_n), & \gamma_{1,n} \in (0, 1)\\
0, & \gamma_{1,n}=1\\
-\frac{1-\gamma_{1,n}}{\gamma_{1,n}}\log(w_n^*-d_n^*/h_n), & \gamma_{1,n} >1\,
\end{array}                                       \right.
\end{eqnarray*}
}
with  $w_n^*=\frac{h_n}{\sqrt{\gamma_{2,n}}}$ and $d_n^*=\sqrt{\gamma_{2,n}}$.
\end{proposition}

\begin{figure}[H]
  \centering
  \includegraphics[scale=0.22]{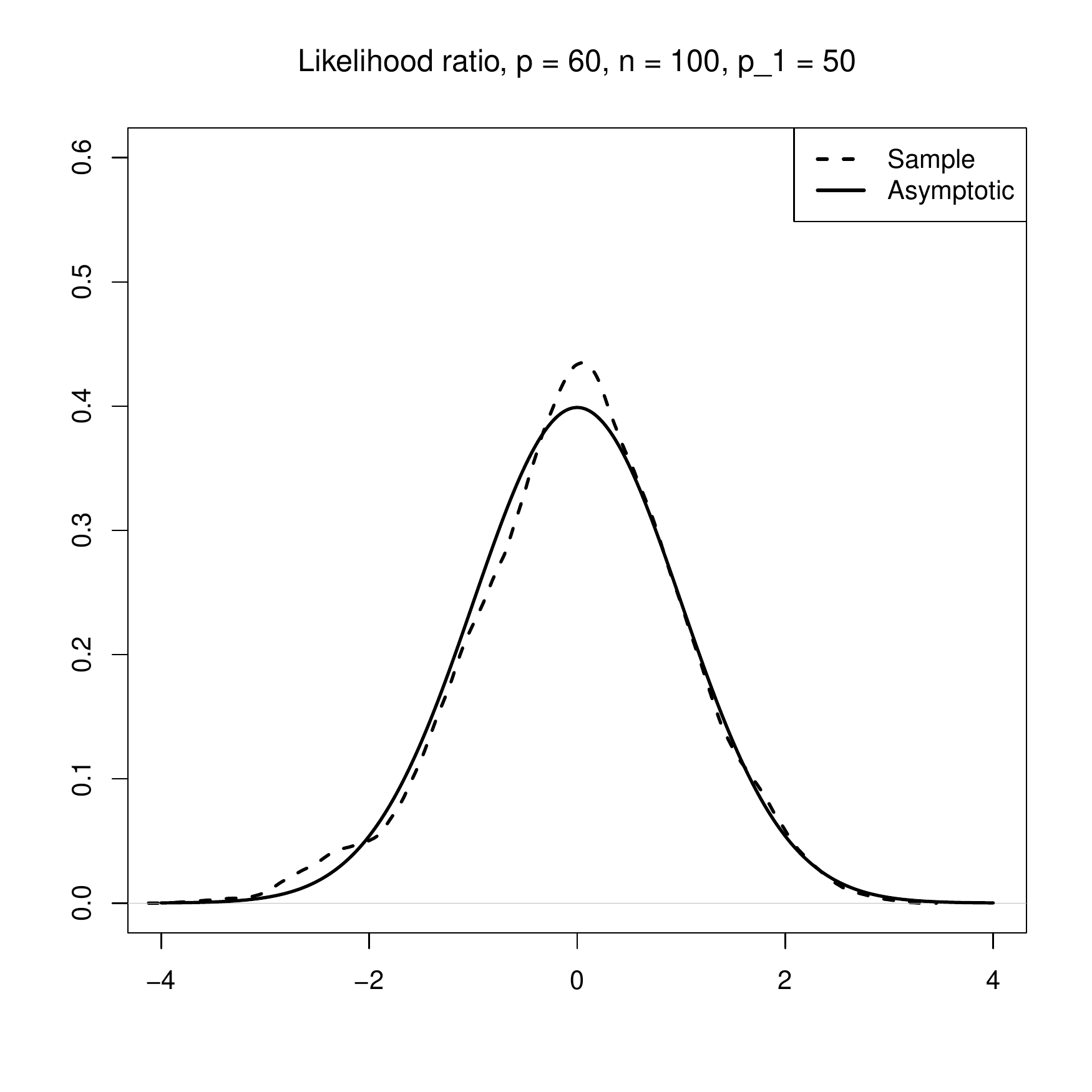} ~~
  \includegraphics[scale=0.22]{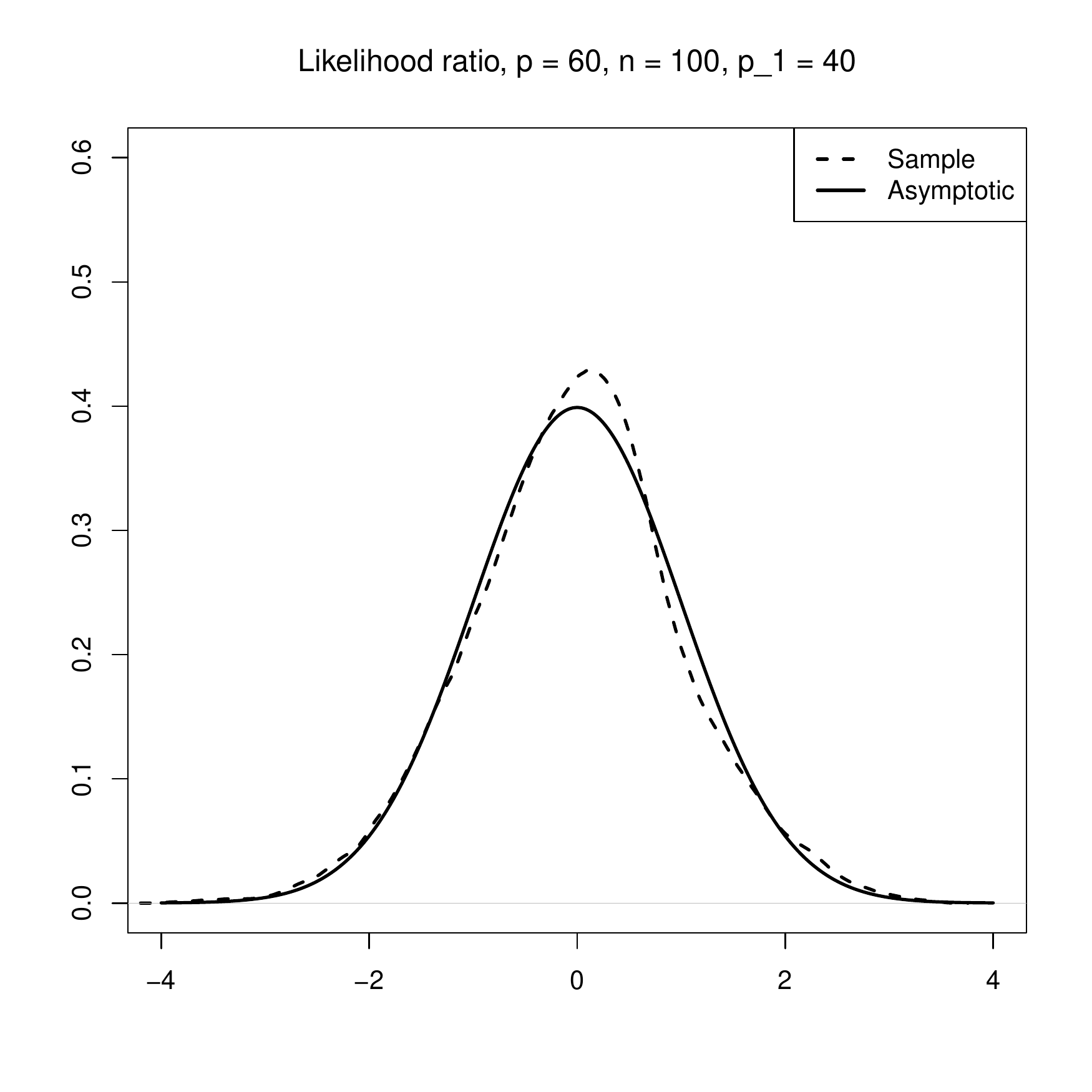} ~~
  \includegraphics[scale=0.22]{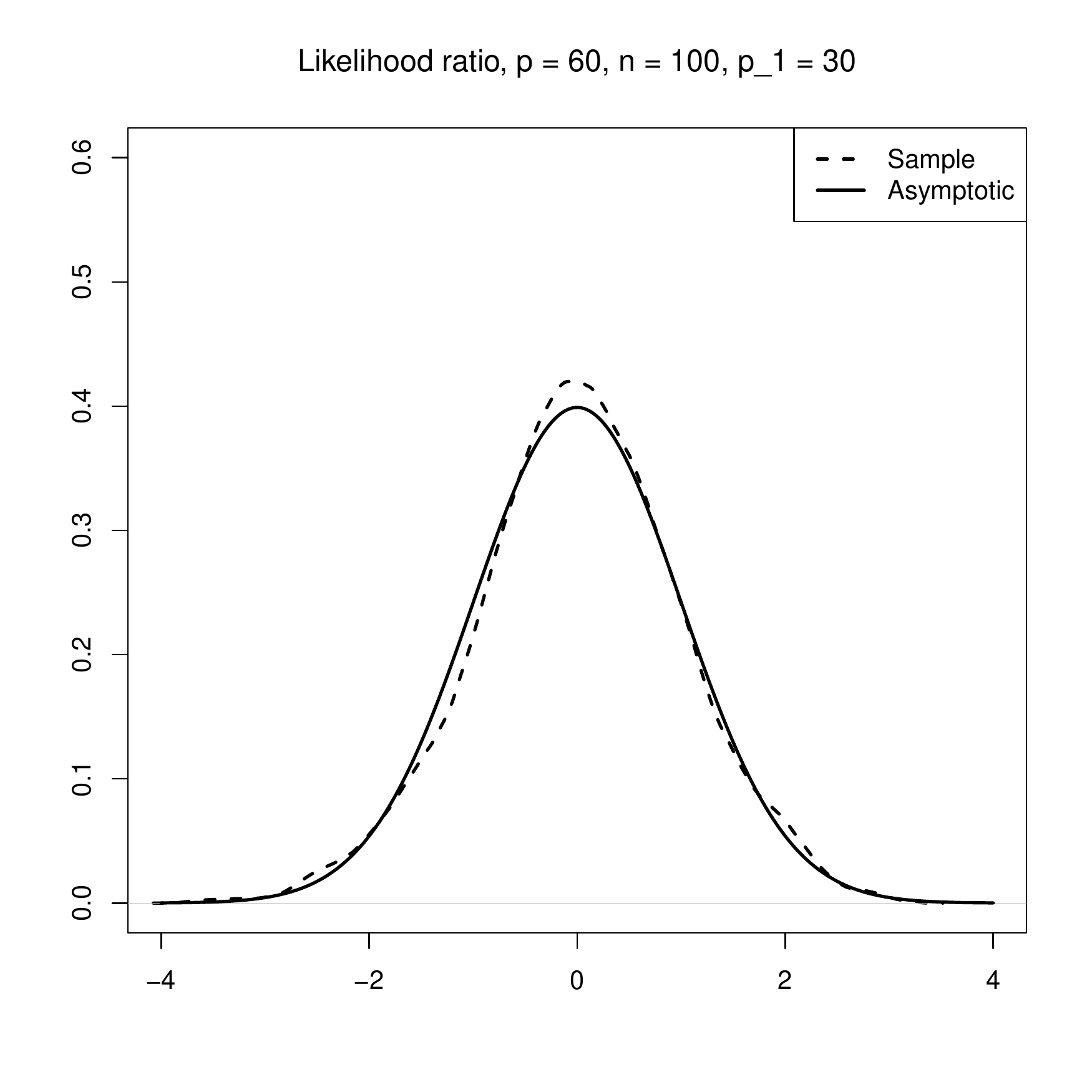}\\
    \includegraphics[scale=0.22]{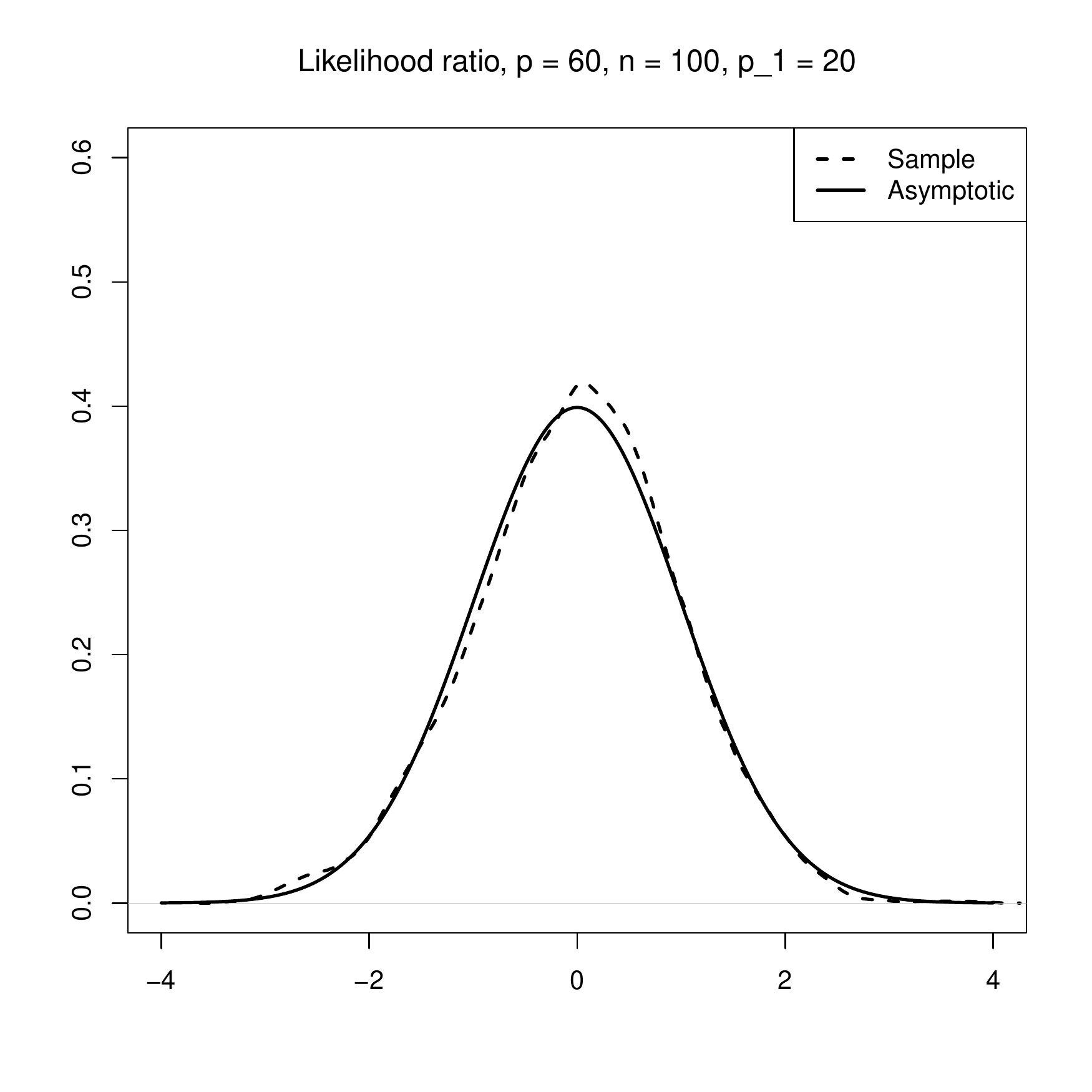} ~~
    \includegraphics[scale=0.22]{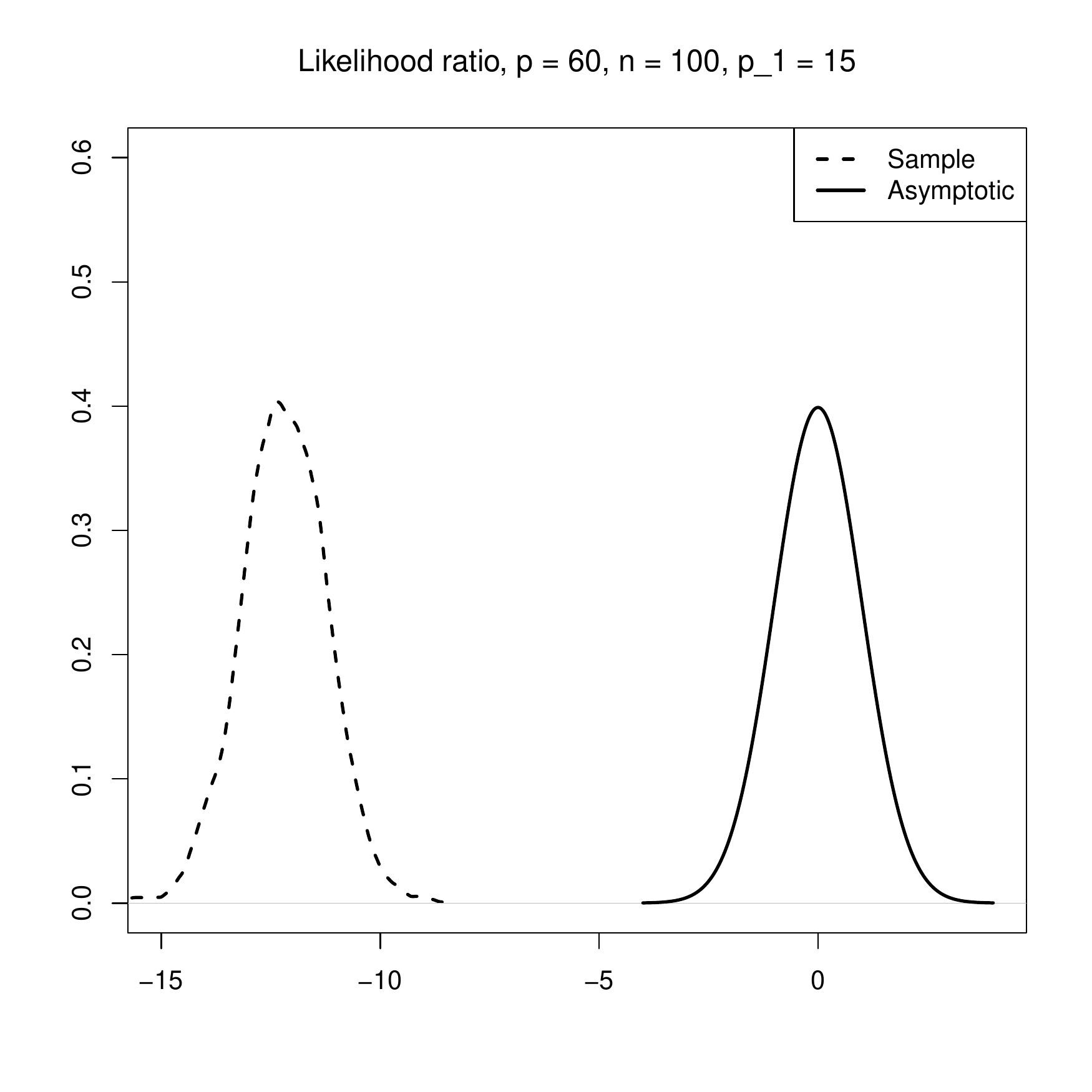} ~~
  \includegraphics[scale=0.22]{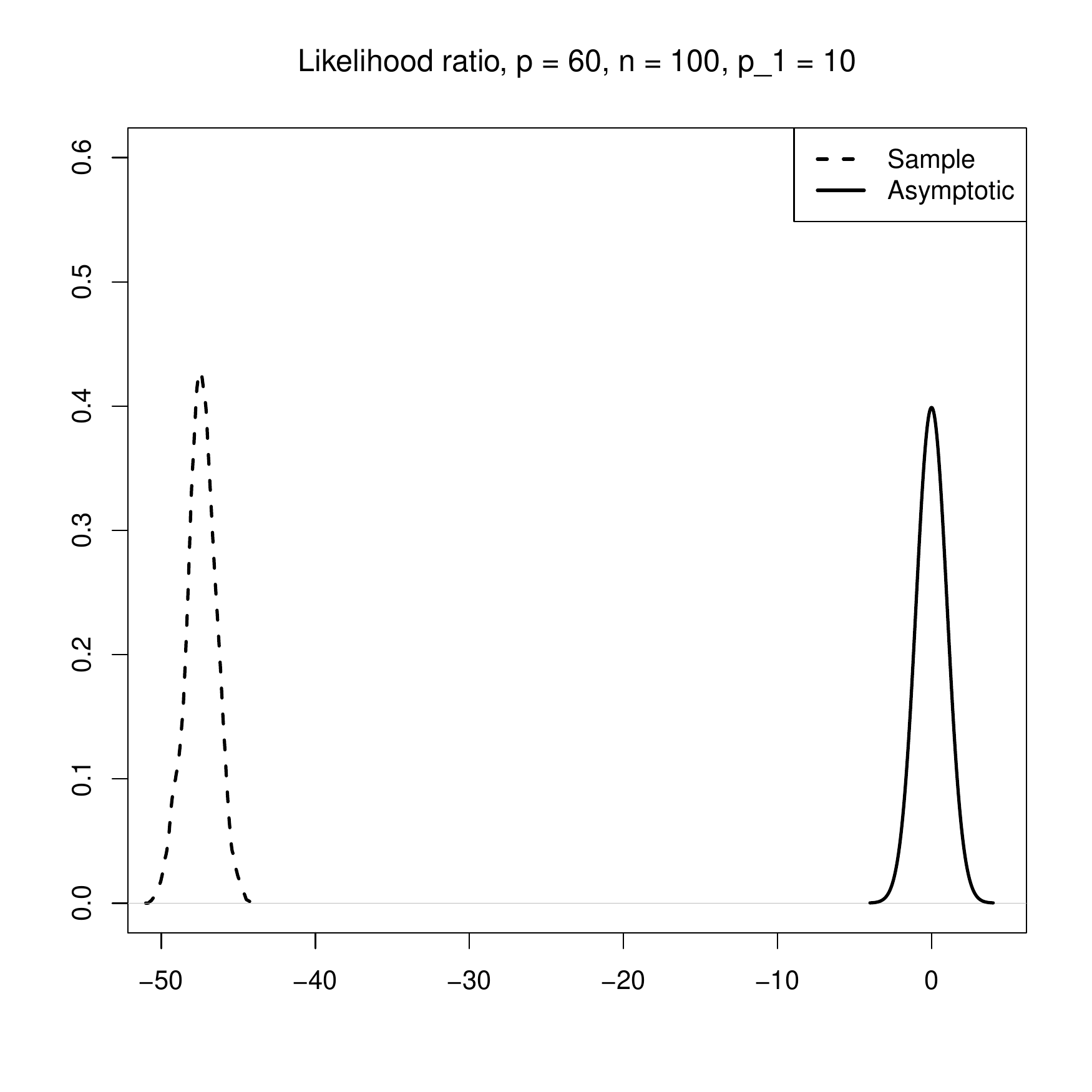}
\vspace{-.5cm}
  \caption{\it Simulated distribution of the statistic $(T_{LR} -(p-p_1)s_{LR}  - \mu_{LR}) /  \sigma_{LR} $ und the null hypothesis
  for  sample size $n=100$,  dimension $p=60$ and  various values of $p_1=50, 45, 40, 30,15, 10$. The solid curve shows the standard normal distribution.
  \label{fig1}}
\end{figure}

\noindent
Proposition \ref{prop1}  shows that the likelihood ratio test,
 which rejects the null hypothesis, whenever
\begin{equation} \label{textLR}
\dfrac{T_{LR} - (p-p_1)s_{LR}-\mu_{LR}}{ \sigma_{LR}}  < -u_{1-\alpha},
\end{equation}
is an asymptotic level $\alpha$ test (here and throughout this paper $u_{1-\alpha}$
denotes the $(1-\alpha)$-quantile of the standard normal distribution).
In Figure \ref{fig1} we illustrate the  approximation of the nominal level of the test \eqref{textLR}  by means of a small simulation study for the sample size $n=100$,
dimension $ p=60$ and different values  of $p_1$ and $p_2$. We considered a centered $p$-dimensional normal distribution
where the blocks $\bSigma_{11}$ and  $\bSigma_{22}$  in the block diagonal
matrix $\bSigma$  are constructed as follows. For the first block $\bSigma_{11}$ we took $p_1$ uniformly distributed  eigenvalues on the interval $(0, 1]$
 while the corresponding eigenvectors are simulated from the Haar distribution on the unit sphere.
 The $p_2$ eigenvalues of the second block  $\bSigma_{22}$  are drawn from a uniform distribution on the interval
 $[1, 10]$ while the corresponding eigenvectors are  again Haar distributed. The matrices $\bSigma_{11}$ and $\bSigma_{22}$ are then fixed
 for the generation of multivariate normal distributed random variables ($\bSigma_{12} = \mathbf{O}$).
   The plots show the empirical  distribution of the statistic $(T_{LR} -(p-p_1)s_{LR}  - \mu_{LR}) /  \sigma_{LR} $
using $1000$ simulation runs and the density of a standard normal distribution.
We observe a reasonable approximation  if the dimension $p_1$ of the sub-vector  $\bx_{1,1}$
is large compared to the  dimension $p$ of the vector  $\bx_1$, that is $\gamma_{1,n} \leq 1$ (see, the upper part of Figure \ref{fig1}). However, if $\gamma_{1,n} >>1$,
 there arises a strong bias (see, the lower part of Figure \ref{fig1}) and the asymptotic statement in Proposition \ref{prop1} cannot be used to obtain  critical value for the test \eqref{textLR}.

\noindent
Motivated by the poor quality of the  approximation of the finite sample distribution of the likelihood ratio test by a normal distribution if the
 dimension $p_1$ is small compared to the dimension $p_2$  we now construct alternative tests for the  hypothesis (\ref{hyp_H0}), which
 will yield a more stable approximation of the nominal
level. For this purpose, we first note that a non-singular partitioned matrix $\bSigma$ in \eqref{SigmaS} is block diagonal (i.e. $\bSigma_{21}=\mathbf{O}$) if and only if
$ \bSigma_{21}\bSigma_{11}^{-1}\bSigma_{12}=\mathbf{O} $. Therefore, a test for independence can also be obtained by testing the hypotheses
\begin{equation}\label{hyp_H0i_2}
H_{0}:~\bSigma_{21}\bSigma_{11}^{-1}\bSigma_{12}=\mathbf{O}   ~ \mbox{ versus } ~
H_{1}:~ \bSigma_{21}\bSigma_{11}^{-1}\bSigma_{12}\neq \mathbf{O}.
\end{equation}
In the following section we will propose three tests for the hypothesis \eqref{hyp_H0i_2} as an alternative to the likelihood ratio test.

\section{Alternative tests for independence and their null distribution} \label{sec3}
\def\theequation{3.\arabic{equation}}
\setcounter{equation}{0}

Recall the definition of the matrices $\bSigma$ and $\bS$ in \eqref{SigmaS} and
denote by
$\bSigma_{22\cdot 1}=\bSigma_{22}-\bSigma_{21}\bSigma_{11}^{-1}\bSigma_{12}$ and $\bS_{22\cdot 1}=\bS_{22}-\bS_{21}\bS_{11}^{-1}\bS_{12}$ the corresponding Schur complements.
From Theorem 3.2.10 of \cite{muirhead1982}, it follows that
\begin{eqnarray*}
\bS_{21}\bS_{11}^{-1/2}|\bS_{11}&\sim& \mathcal{N}_{p-p_1,p_1}(\bSigma_{21}\bSigma_{11}^{-1}\bS_{11}^{1/2},\bSigma_{22\cdot 1} \otimes \bI_{p_1}),\\
\bS_{22\cdot 1}&\sim& W_{p-p_1}(n-p_1, \bSigma_{22\cdot 1}),
\end{eqnarray*}
and the Schur complement $\bS_{22\cdot 1}$ is independent of $\bS_{21}\bS_{11}^{-1/2}$ and $\bS_{11}$. Hence, under the null hypothesis,
\begin{eqnarray}
\nonumber
\widehat{\bW}&=&\bS_{21}\bS_{11}^{-1}\bS_{12}\sim W_{p-p_1}(p_1,\bSigma_{22\cdot 1}),\\
\widehat{\bT}&=&\bS_{22\cdot 1}\sim W_{p-p_1}(n-p_1,\bSigma_{22\cdot 1}), \label{centW}  \nonumber
\end{eqnarray}
and $\widehat{\bW}$ and $\widehat{\bT}$ are independent. Under the alternative hypothesis $H_{1}$, $\widehat{\bW}$ and $\widehat{\bT}$ are still independent as well as $\widehat{\bT}\sim W_{p-p_1}(n-p_1,\bSigma_{22\cdot 1})$, but $\widehat{\bW}$ has a non-central Wishart distribution conditionally on $\bS_{11}$, i.e.,
\begin{equation} \label{noncW}  \nonumber
\widehat{\bW}|\bS_{11}\sim W_{p-p_1}(p_1,\bSigma_{22\cdot 1},\bOmega_1(\bS_{11}))
\end{equation}
where the non-centrality parameter is given by
\begin{equation} \label{noncpar}  \nonumber
\bOmega_1= \bOmega_1(\bS_{11})=\bSigma_{22\cdot 1}^{-1}\bSigma_{21}\bSigma_{11}^{-1}\bS_{11}\bSigma_{11}^{-1}\bSigma_{12}  .
\end{equation}
For technical reasons we will use the normalized versions of $\widehat{\bW}$ and $\widehat{\bT}$ throughout this paper. Thus, the distributional properties of $\bW=\frac{1}{p_1}\widehat{\bW}$ and $\bT=\frac{1}{n-p_1}\widehat{\bT}$ are very similar to the ones observed for the within and between covariance matrices in the multivariate
analysis of variance  (MANOVA) model (see \cite{fujhimwak2004}, \cite{schott2007}, \cite{kakiiwas2008}). More precisely, $p_1\bW$ and $(n-p_1)\bT$ are independent (under both hypotheses) and they possess Wishart distributions under the null hypothesis. However under the alternative hypothesis the matrix $p_1\bW$ has only conditionally on $\bS_{11}$ a non-central Wishart distribution, while the unconditional distribution appears to be a more complicated matrix-variate distribution.
The similarity to MANOVA motivates the application of three tests which are usually used in this context and are given by
\begin{enumerate}[(i)]
\item Wilks' $\Lambda$ statistic:
  \begin{equation}\label{W}
T_W = -\log(|\bT|/|\bT +\bW |)=\log(|\bI +\bW\bT^{-1}|)=\sum_{i=1}^{p-p_1} \log(1+v_i)
\end{equation}
\item Lawley-Hotelling's trace criterion:
\begin{equation}\label{LH}
T_{LH} = tr(\bW\bT^{-1})=\sum_{i=1}^{p-p_1} v_i
\end{equation}
\item Bartlett-Nanda-Pillai's trace criterion:
\begin{equation}\label{BNP}
T_{BNP} = tr(\bW\bT^{-1}(\bI +\bW\bT^{-1})^{-1})=\sum_{i=1}^{p-p_1} \frac{v_i}{1+v_i}
\end{equation}
\end{enumerate}
where $v_1 \ge v_2 \ge ... \ge v_{p-p_1}$ denote the ordered eigenvalues of the matrix $\bW\bT^{-1}$.
A statistic very similar  to \eqref{BNP} was proposed by \cite{jiangetal2013}, who used
 $$
 tr(\bW\bT^{-1}(\frac{\gamma_1}{\gamma_2}\bI +\bW\bT^{-1})^{-1})=\sum_{i=1}^{p-p_1} \frac{v_i}{\frac{\gamma_1}{\gamma_2}+v_i}
 $$
  instead of $tr(\bW\bT^{-1}(\bI +\bW\bT^{-1})^{-1})$.
It is remarkable that all proposed test statistics are functions of the eigenvalues of $\bW\bT^{-1}$ and
can be presented as linear spectral statistics calculated for the random matrix $\bW\bT^{-1}$, which is the so-called Fisher matrix under the null hypothesis $H_{0}$ (see \cite{zheng2012}).

A linear spectral statistics  for the matrix $\bW\bT^{-1}$ is generally defined by
\begin{equation}\label{def_LSS}
LSS_n=(p-p_1)\int_0^\infty f(x) \diff F_n(x) = \sum_{i=1}^{p-p_1} f(v_i) \,,
\end{equation}
where  $v_1 \ge v_2 \ge ... \ge v_{p-p_1}$ are the ordered eigenvalues of the matrix $\bW\bT^{-1}$. The symbol
\begin{equation*}
F_n(x)= \frac{1}{p-p_1} \sum_{i=1}^{p-p_1} \mathds{1}_{(-\infty, x]}(v_i)
\end{equation*}
denotes the corresponding empirical spectral distribution and the symbol $\mathds{1}_{\mathcal{A}}$ is the indicator function of the set $\mathcal{A}$.  Define
{\small
\begin{eqnarray}\label{centralFemp} \nonumber
  F_n^*(dx) &=& q_n(x)\mathbbm{1}_{[a_n, b_n]}(x)dx + (1-1/\gamma_{1,n} )\mathbbm{1}_{\gamma_{1,n}>1}\delta_0(dx)~~\text{with}\\[5mm]
  q_n(x)   &=& \frac{1-\gamma_{2,n}}{2\pi x(\gamma_{1,n}+\gamma_{2,n}x)}\sqrt{(b_n-x)(x-a_n)},\nonumber \label{WTlsdemp}\\
  a_n&=& \frac{(1-h_n)^2}{(1-\gamma_{2,n})^2},~~b_n=\frac{(1+h_n)^2}{(1-\gamma_{2,n})^2}\,, \nonumber
\end{eqnarray}
}
where $\gamma_{1,n}$,  $\gamma_{2,n}$ and $h_n$ are defined by \eqref{g1}, \eqref{g2} and \eqref{g3}, respectively.
Note that     $F_n^*$ is  a finite sample proxy of limiting spectral distribution  $F$ of $F_n$, which is obtained
by   replacing $\gamma_{1,n}$ and $\gamma_{2,n}$ by their corresponding limits (see \cite{baisilverstein2010}), that is
{\small
\begin{eqnarray}\label{centralF}
  F(dx) &=& q(x)\mathbbm{1}_{[a, b]}(x)dx + (1-1/\gamma_1)\mathbbm{1}_{\gamma_1>1}\delta_0(dx)~~\text{with}\\[5mm]
  q(x)   &=& \frac{1-\gamma_2}{2\pi x(\gamma_1+\gamma_2x)}\sqrt{(b-x)(x-a)},\label{WTlsd}\\
  a&=& \frac{(1-h)^2}{(1-\gamma_2)^2},~~b=\frac{(1+h)^2}{(1-\gamma_2)^2}\nonumber\,.
\end{eqnarray}
}{\small
where
\begin{eqnarray*}
 \gamma_1& =&  \lim_{n\to\infty} \gamma_{1,n}
 = \lim_{n\to\infty}\frac{p-p_1}{p_1} ~,~~
  \gamma_2 = \lim_{n\to\infty} \gamma_{2,n} =\lim_{n\to\infty}\frac{p-p_1}{n-p_1},\\[2mm]
 h &= & \lim_{n\to\infty} h_n = \sqrt{\gamma_{1}+\gamma_{2}-\gamma_{1}\gamma_{2}}\,.
\end{eqnarray*}
}
The representations of $T_W$, $T_{LH}$, and $T_{BNP}$ in terms of the eigenvalues of the random matrix $\bW\bT^{-1}$ are used intensively in the proof of our first main result, which provides their asymptotic distribution under the null hypothesis in \eqref{hyp_H0i_2}. The details of the proof are deferred to Appendix B of the 
 supplementary material (see, \cite{supplement}).

\begin{theorem}\label{th2} Under the assumptions stated in Section \ref{sec2}  we have
      \begin{eqnarray*}
    &&\frac{T_{a} - (p-p_1)s_{\alpha} -\mu_{a} }{\sigma_{a}} \overset{\mathcal{D}}{\longrightarrow} \mathcal{N}(0, 1)    \,,
    \end{eqnarray*}
    where  the index $a\in\{W, LH, BNP \}$ represents the statistic under  consideration defined in \eqref{W}, \eqref{LH} and \eqref{BNP}, respectively. The asymptotic means and variances are given by
    \begin{eqnarray*}
      \mu_{W} &=& 1/2\log\left[ \frac{(w_n^2-d_n^2)h_n^2}{(w_nh_n-\gamma_{2,n}d_n)^2}\right],\qquad\qquad \sigma^2_{W} =2\log\left[\frac{w_n^2}{w_n^2-d_n^2}\right],\\
      \mu_{LH} &=&  \frac{\gamma_{2,n}}{(1-\gamma_{2,n})^2},\qquad\qquad\qquad\qquad\quad~ \sigma^2_{LH} =\frac{2h_n^2}{(1-\gamma_{2,n})^4},\\
      \mu_{BNP} &=&-\frac{(1-\gamma_{2,n})^2w_n^2(d_n^2-\gamma_{2,n})}{(w_n^2-d_n^2)^2},\\
      \sigma^2_{BNP}& =& 2\frac{d^2(1-\gamma_{2,n})^4(w_n^2(w_n^2+d_n)+d_n^3(w_n^2-1))}{w_n^2(1+d_n)(w_n^2-d_n^2)^4} \, ,
\end{eqnarray*}
where $w_n>d_n>0$ satisfy $w_n^2+d_n^2=(1-\gamma_{2,n})^2+1+h_n^2$, $w_nd_n=h_n$, and the quantities  $\gamma_{1,n}$, $\gamma_{2,n}$ and $h_n$
are defined by \eqref{g1}, \eqref{g2} and  \eqref{g3}, respectively. The centering parameters are given by
{\small
\begin{eqnarray*}
  s_W&=&-\log\left((1-\gamma_{2,n})^2 \right)- \frac{1-\gamma_{2,n}}{\gamma_{2,n}}\log(w_n) + \frac{\gamma_{1,n}+\gamma_{2,n}}{\gamma_{1,n}\gamma_{2,n}}\log(w_n-d_n\gamma_{2,n}/h_n)\\
  &-&
\left\{
\begin{array}{ll}
\frac{1-\gamma_{1,n}}{\gamma_{1,n}}\log(w_n-d_nh_n), & \gamma_{1,n} \in (0, 1)\\
0, & \gamma_{1,n}=1\\
-\frac{1-\gamma_{1,n}}{\gamma_{1,n}}\log(w_n-d_n/h_n), & \gamma_{1,n} >1\,
\end{array}
\right.,\\
      s_{LH}&=& \frac{1}{1-\gamma_{2,n}},\\[5mm]
      s_{BNP}&=&  \frac{1-\gamma_{2,n}}{w_n^2-\gamma_{2,n}}\,.
    \end{eqnarray*}
 }
\end{theorem}
\medskip

\begin{figure}[H]
  \centering
  \includegraphics[scale=0.222]{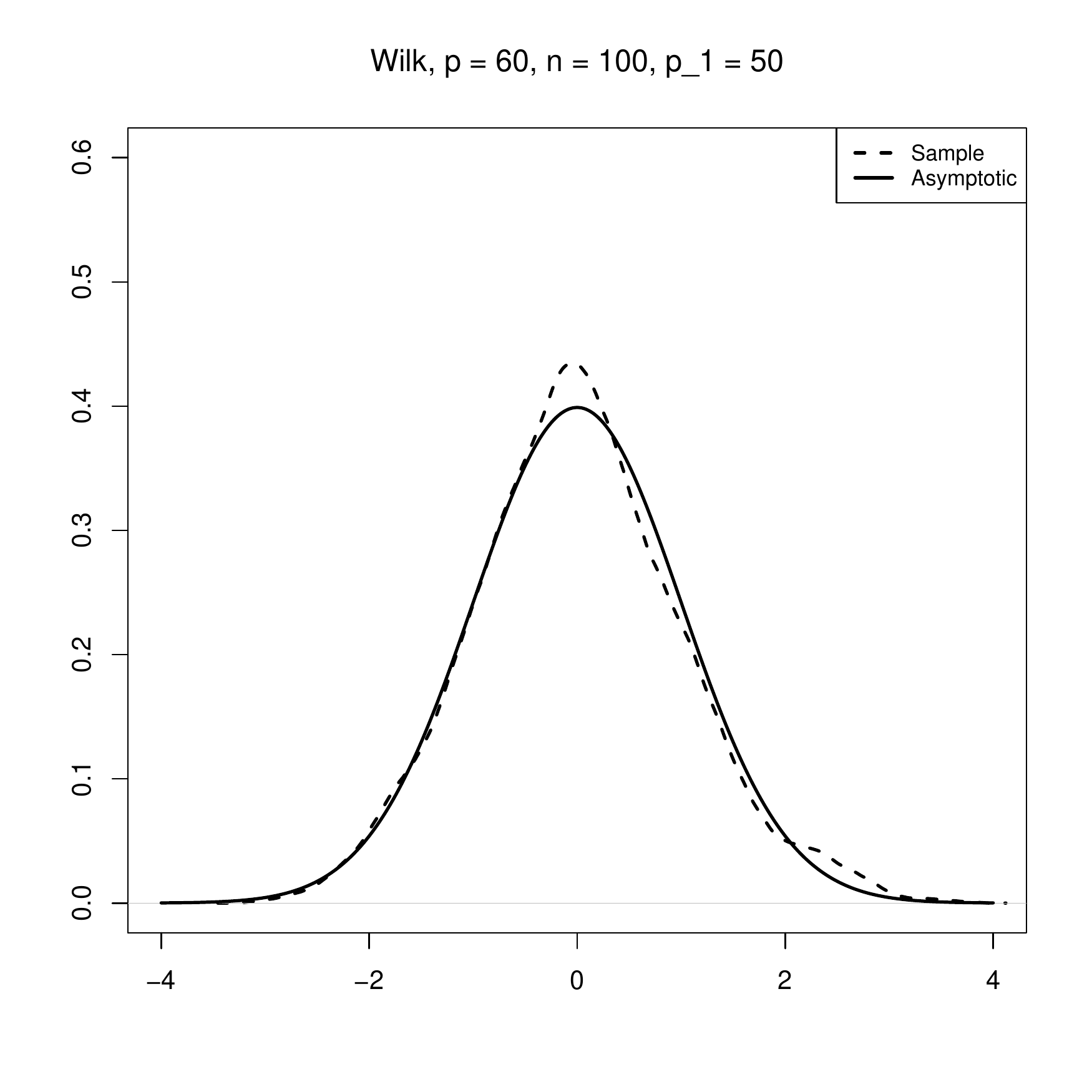} ~~
  \includegraphics[scale=0.222]{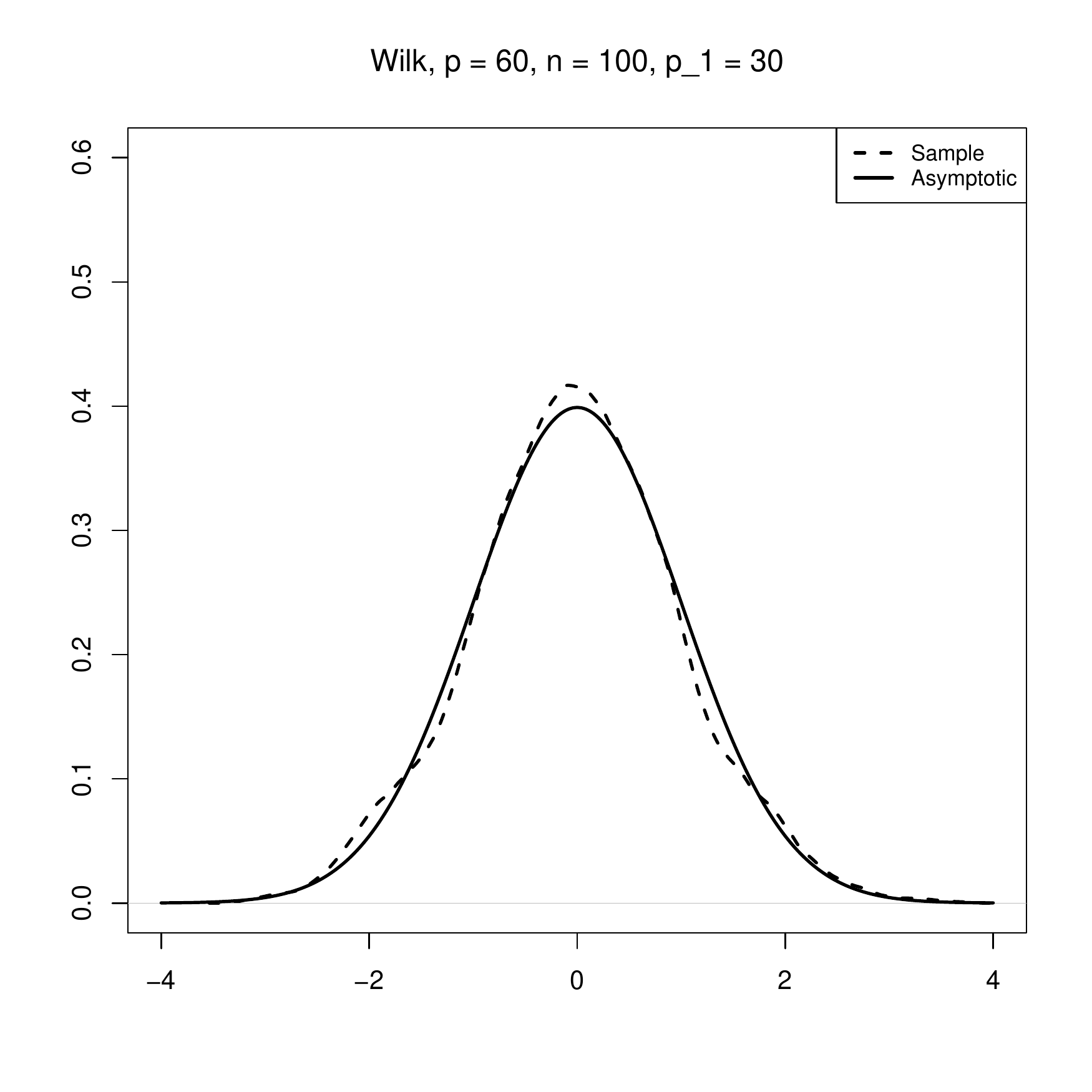} ~~
  \includegraphics[scale=0.222]{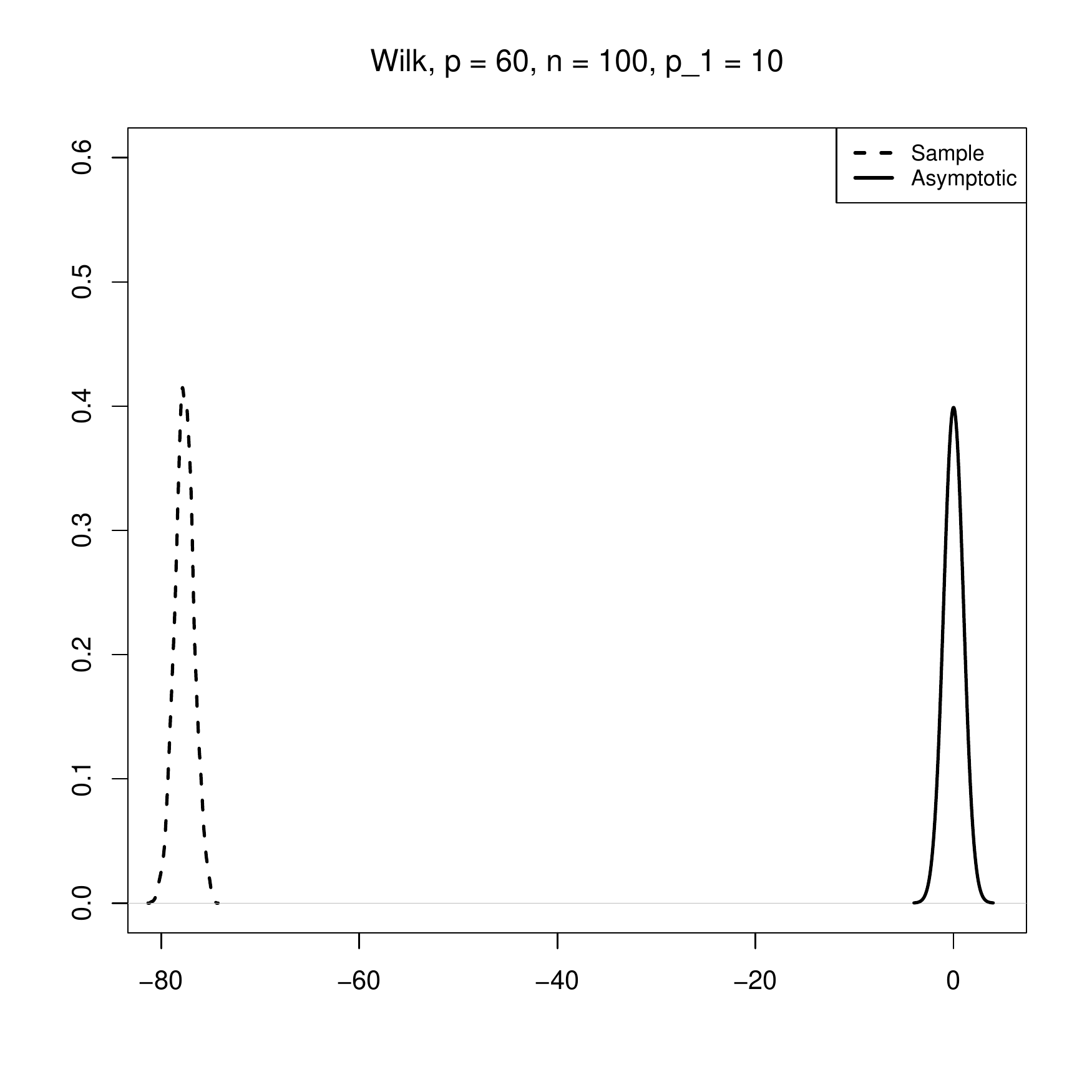}\\

\vspace{-.5cm}
  \caption{\it Simulated distribution of the statistic $(T_{W} -(p-p_1)s_{W}  - \mu_{W}) /  \sigma_{W} $ und the null hypothesis
  for  sample size $n=100$,  dimension $p=60$ and  various values of $p_1=50,  30, 10$.
  The solid curve shows the standard normal distribution.
  \label{fig2}}
\end{figure}

\begin{figure}[H]
  \centering
  \includegraphics[scale=0.222]{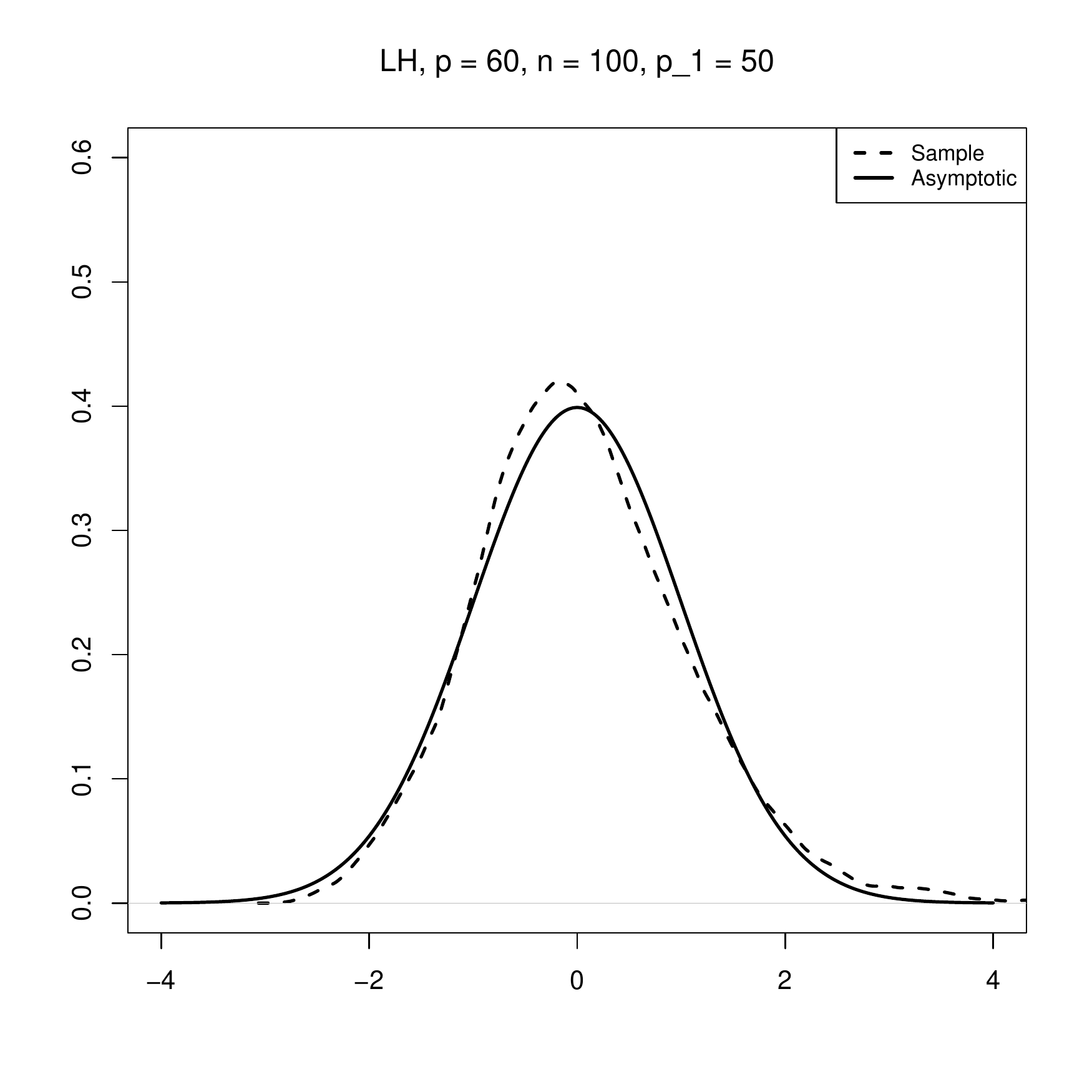} ~~
  \includegraphics[scale=0.222]{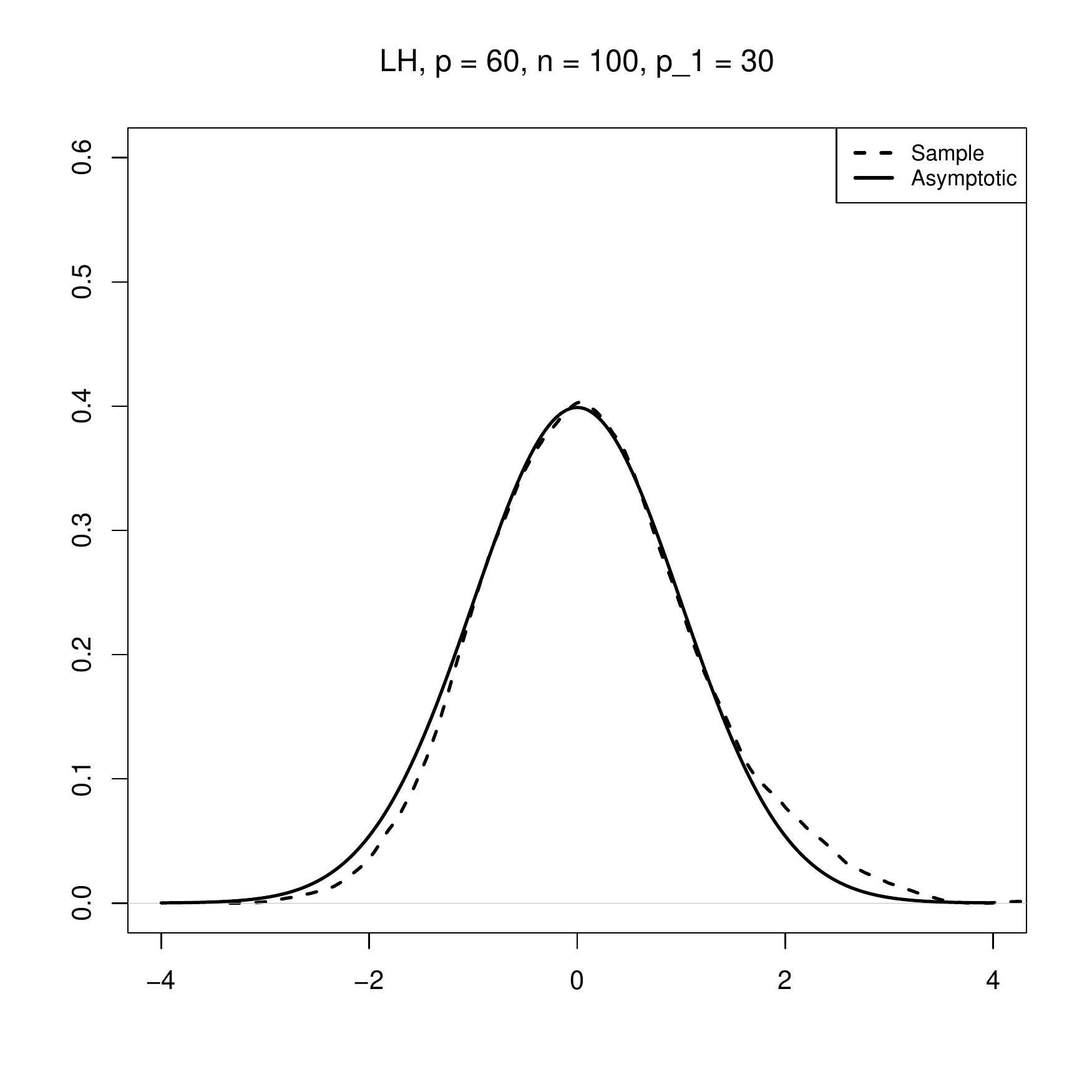} ~~
  \includegraphics[scale=0.222]{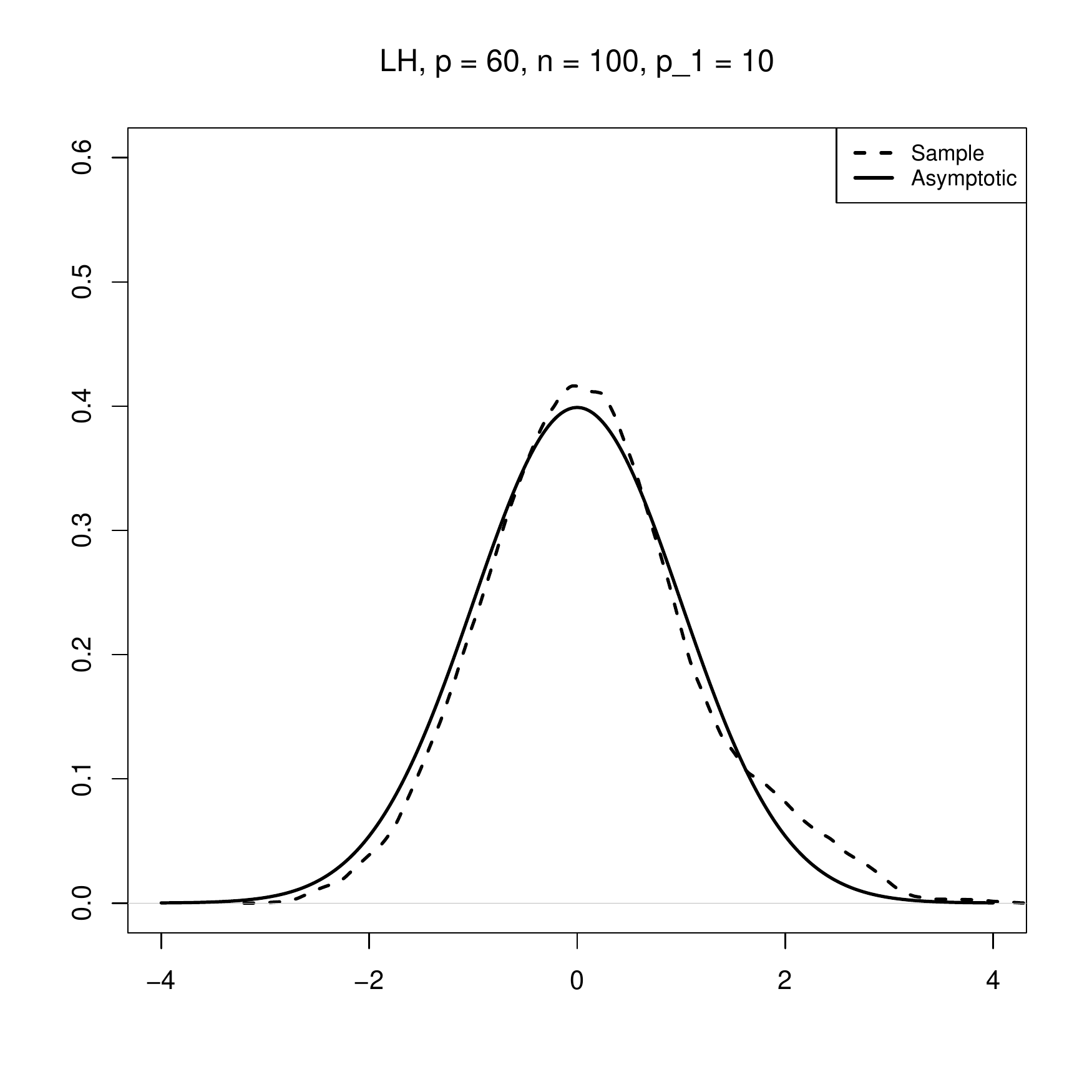}\\

\vspace{-.5cm}
  \caption{\it Simulated distribution of the statistic $(T_{LH} -(p-p_1)s_{LH}  - \mu_{LH}) /  \sigma_{LH} $ und the null hypothesis
 for  sample size $n=100$,  dimension $p=60$ and  various values of $p_1=50,  30, 10$. The solid curve shows the
  standard normal distribution.
  \label{fig3}}
\end{figure}

\begin{figure}[H]
  \centering
  \includegraphics[scale=0.222]{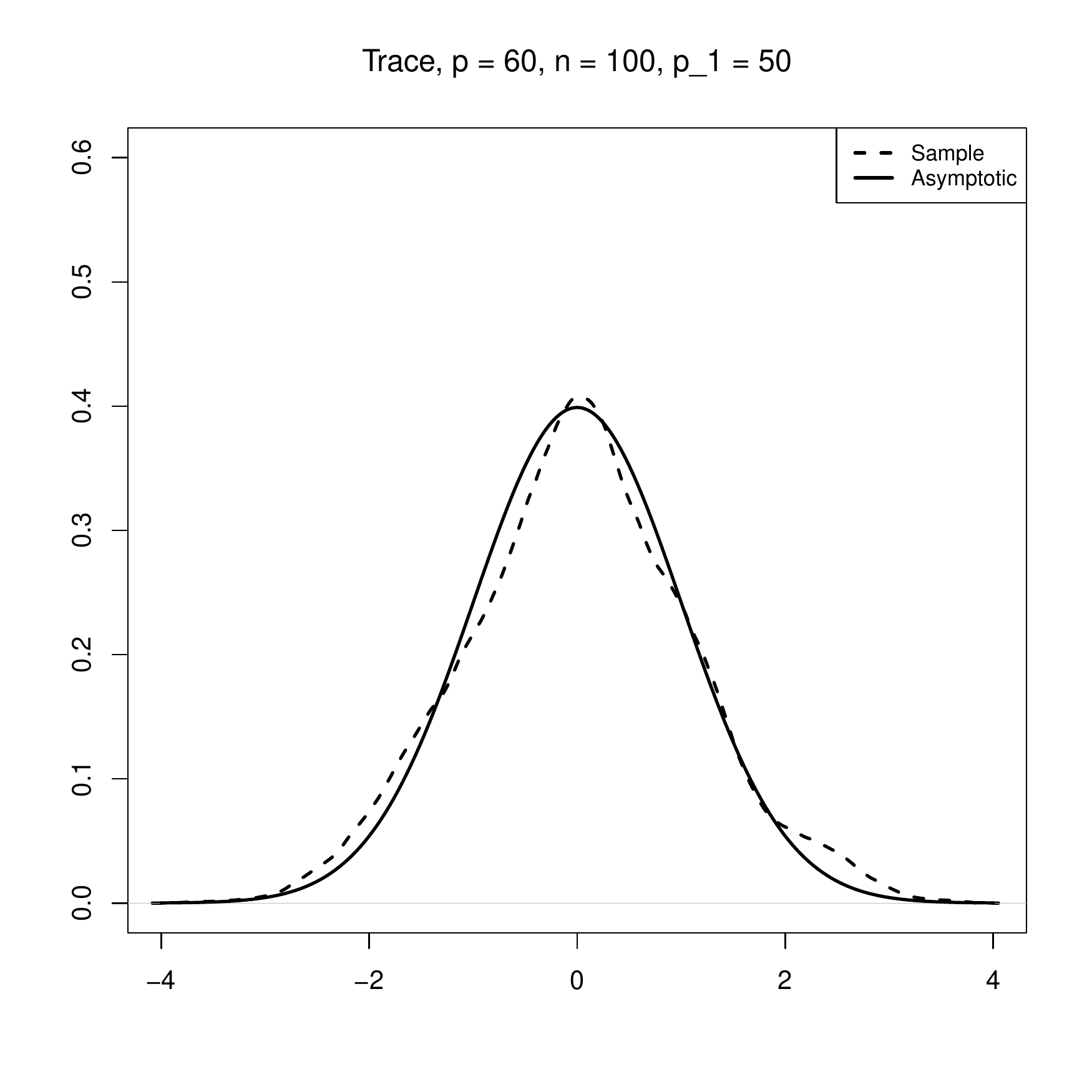} ~~
  \includegraphics[scale=0.222]{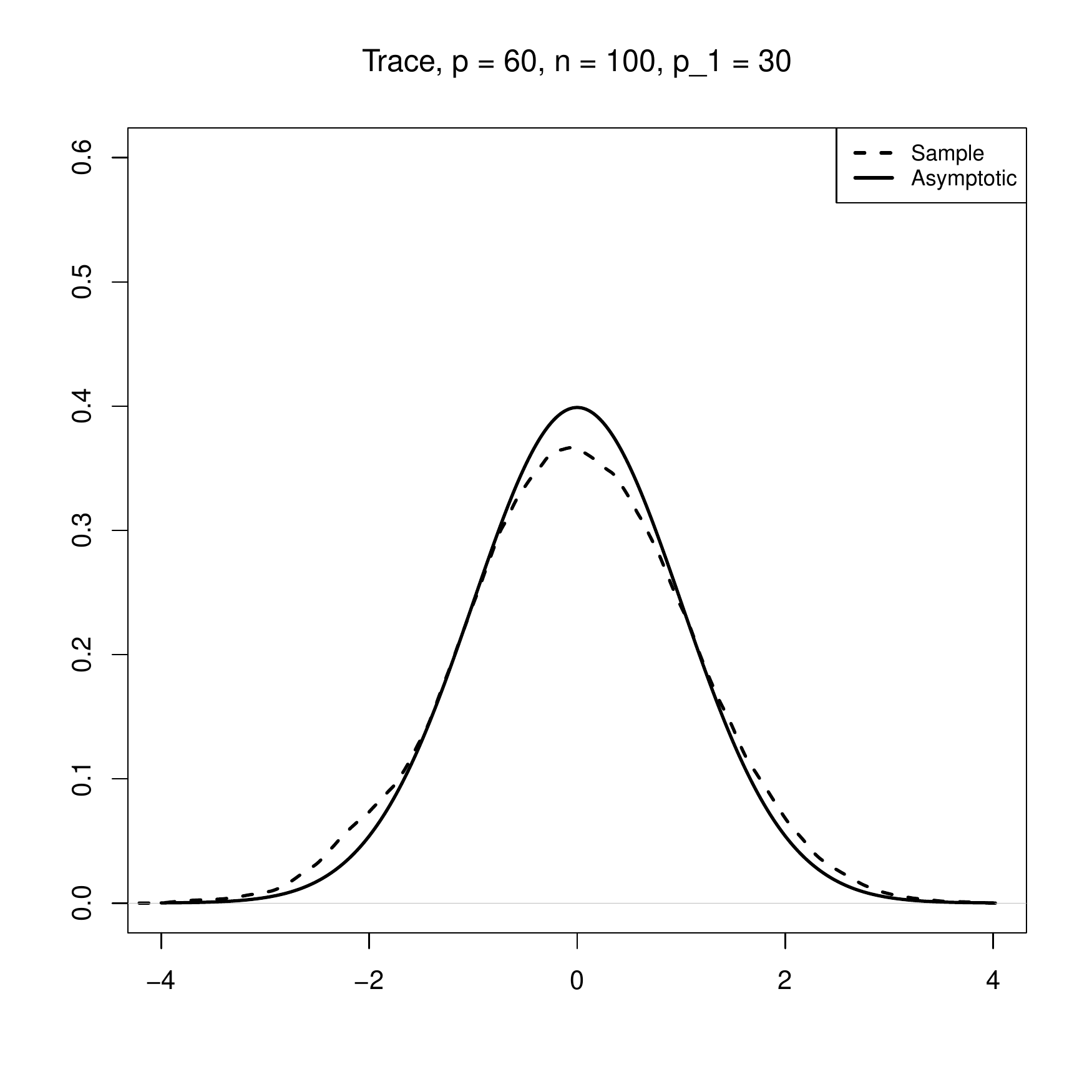} ~~
  \includegraphics[scale=0.222]{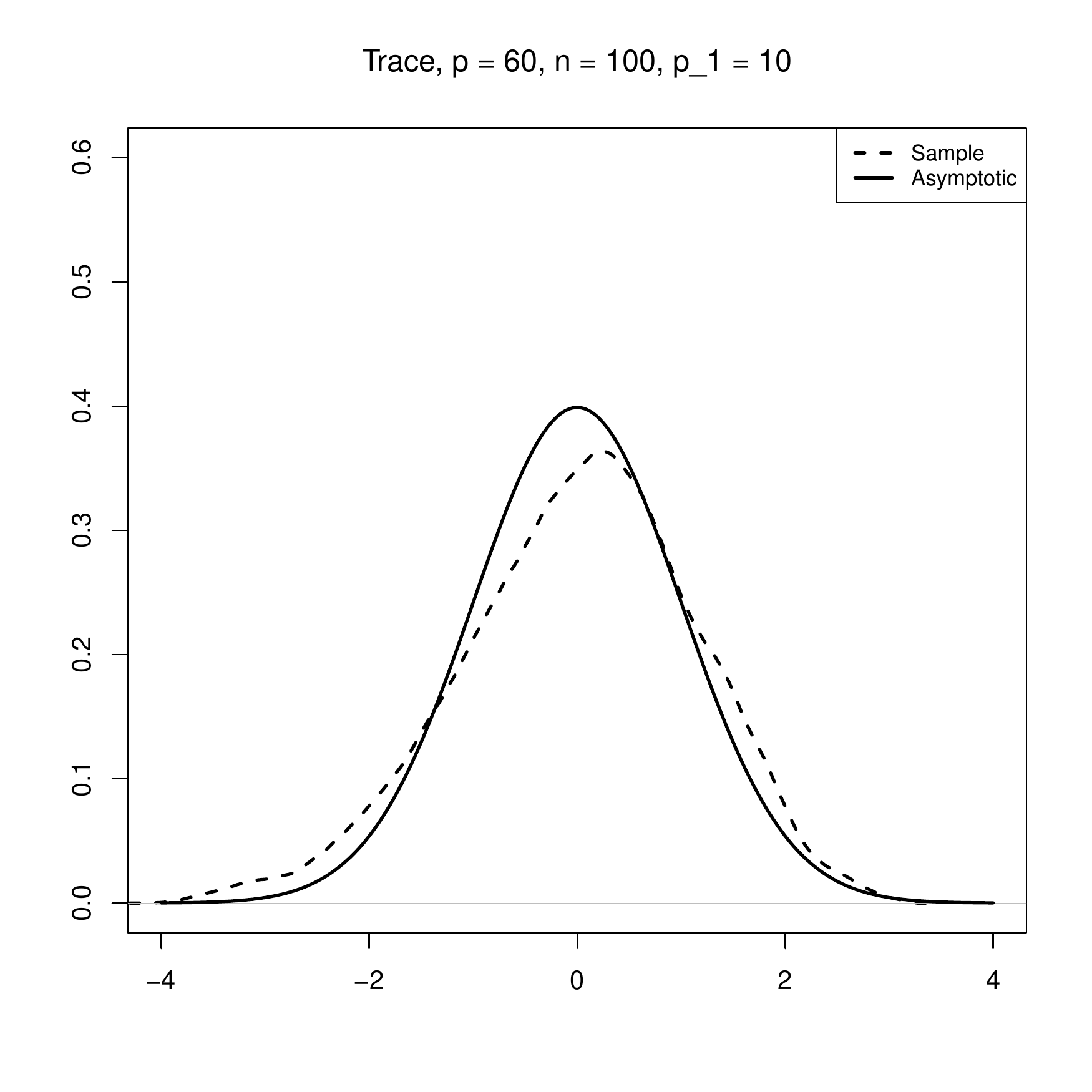}\\

\vspace{-.5cm}
  \caption{\it Simulated distribution of the statistic $(T_{BNP} -(p-p_1)s_{BNP}  - \mu_{BNP}) /  \sigma_{BNP} $ und the null hypothesis
 for  sample size $n=100$,  dimension $p=60$ and  various values of $p_1=50, 30, 10$. The solid curve shows the
  standard normal distribution.
  \label{fig4}}
\end{figure}

\noindent
Theorem \ref{th2} provides a simple asymptotic level $\alpha$ test by rejecting the null hypothesis $H_{0}$ if
\begin{equation} \label{testall}
\dfrac{T_{a} - (p-p_1)s_{a}-\mu_{a}}{ \sigma_{a}} > u_{1-\alpha}
\end{equation}
We illustrate the quality of the approximation in Theorem \ref{th2} by means of a small simulation study. For the sake of comparison with the
likelihood ratio test, we use the same scenario as in Section \ref{sec2}, that is $n=100$, $p=60$
and different values for $p_1$. In Figure \ref{fig2} - \ref{fig4} we display the rejection probabilities of the test
\eqref{testall} under the null hypothesis in the case of  the Wilk test, the Lawley-Hotelling's, and the Bartlett-Nanda-Pillai's trace criterion.
From the results depicted in Figure \ref{fig2} we observe
that the statistic $T_W$ exhibits similar problems as the statistic of the likelihood ratio test. If the dimension $p_1$ is too small the approximation
provided by Theorem \ref{th2} is not reliable. This fact seems to be related to the use of the log determinant criterion.
 On the other hand, the Lawley-Hotelling's and the Bartlett-Nanda-Pillai's trace criterion yield test statistics which do not possess these drawbacks. The results in  Figures \ref{fig3}  and  \ref{fig4}  show a reasonable approximation of the nominal level in all considered scenarios.

\begin{figure}[H]
  \centering
  \includegraphics[scale=0.222]{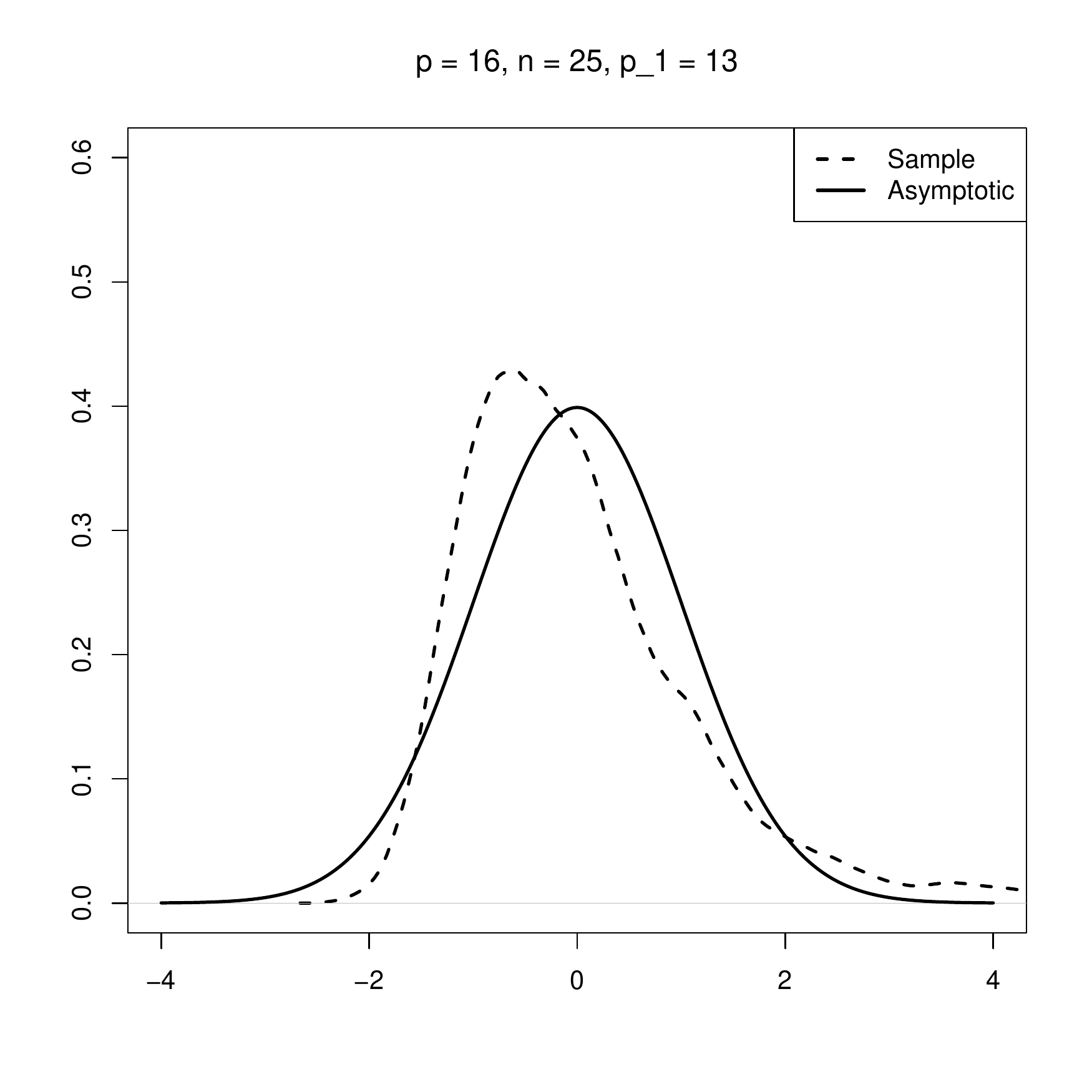} ~~
  \includegraphics[scale=0.222]{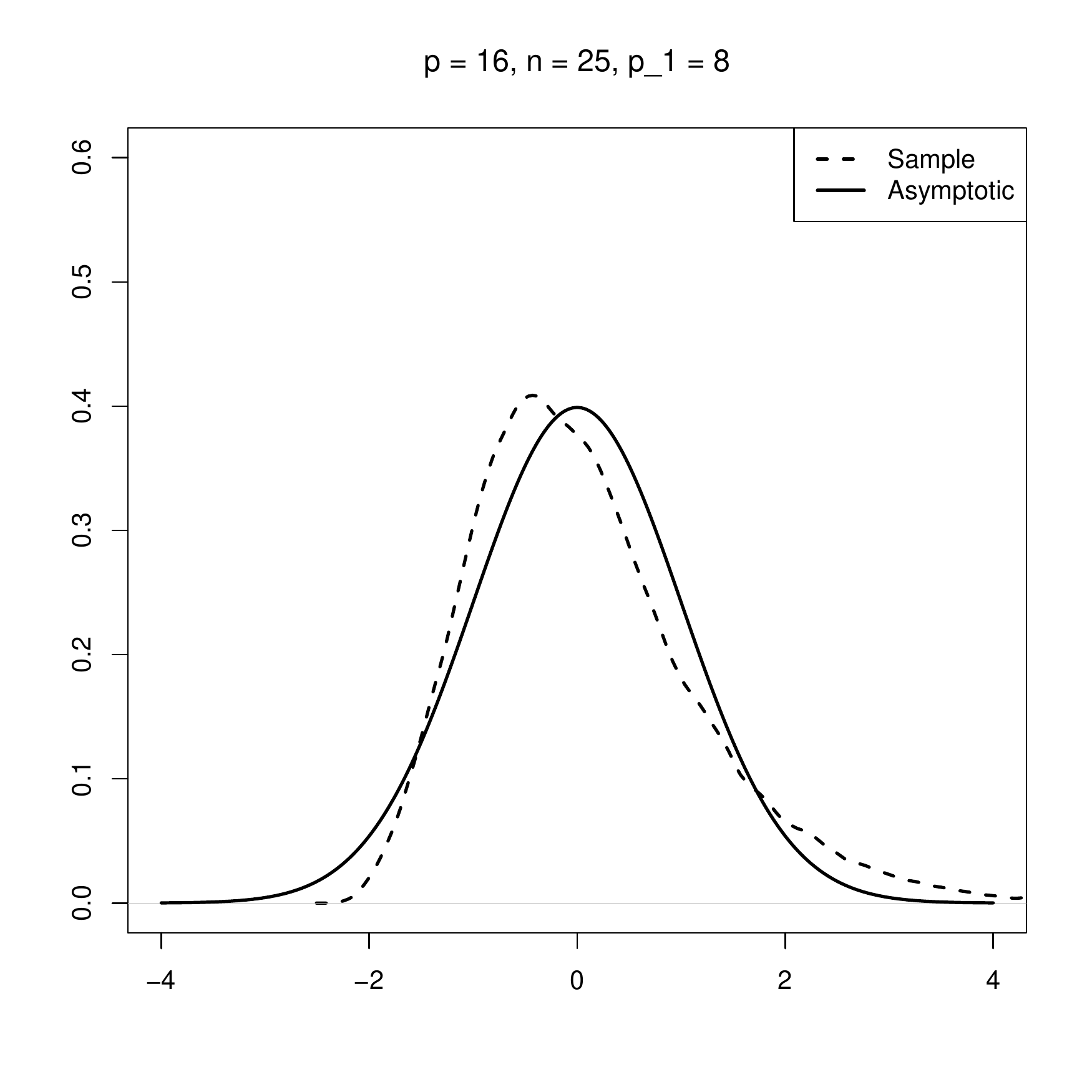} ~~
  \includegraphics[scale=0.222]{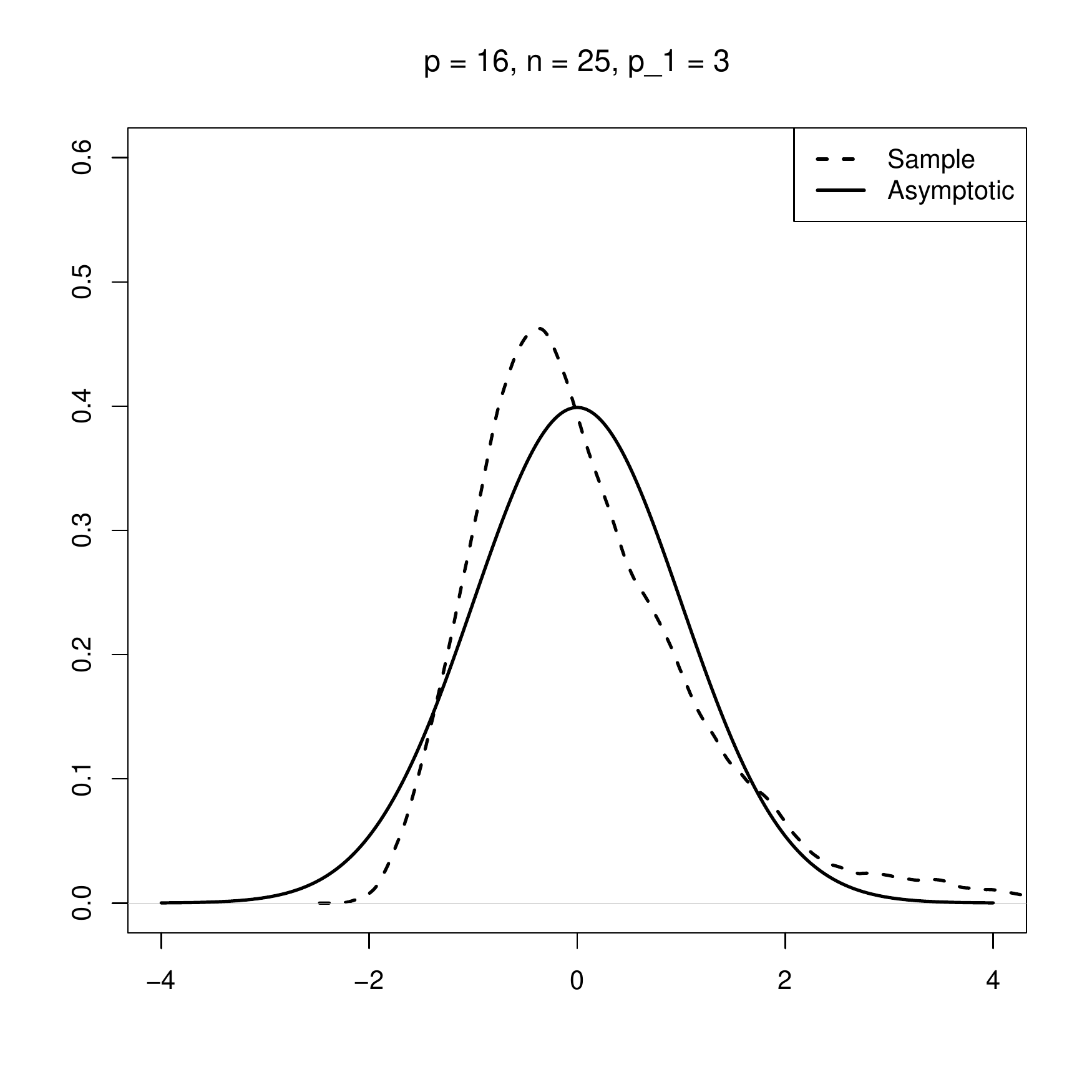}\\

\vspace{-.5cm}
  \caption{\it Simulated distribution of the statistic $(T_{LH} -(p-p_1)s_{LH}  - \mu_{LH}) /  \sigma_{LH} $ und the null hypothesis
 for  sample size $n=25$,  dimension $p=16$ and  various values of $p_1=13,  8, 3$. The solid curve shows the
  standard normal distribution.
  \label{fig51}}
\end{figure}

\begin{figure}[H]
  \centering
  \includegraphics[scale=0.222]{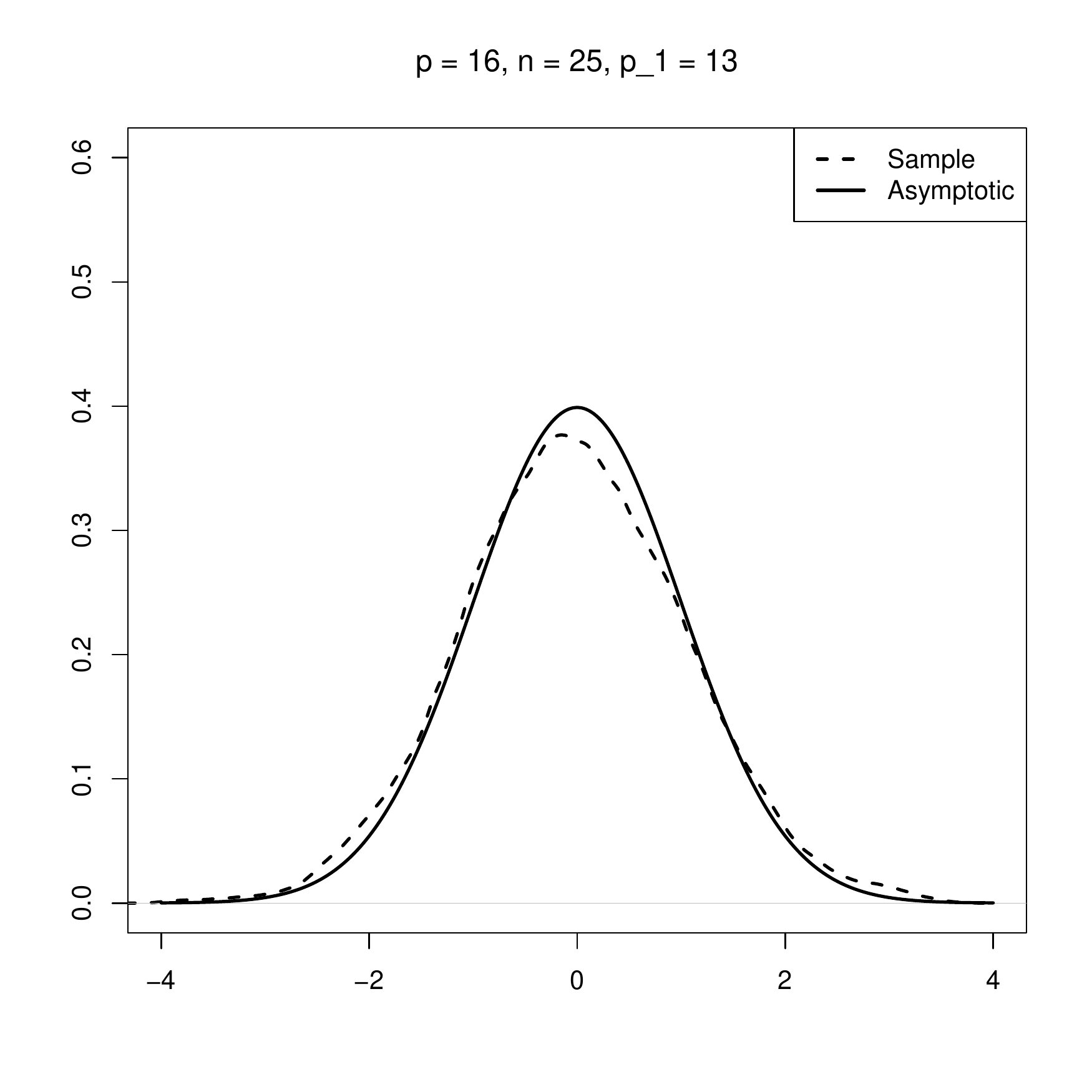} ~~
  \includegraphics[scale=0.222]{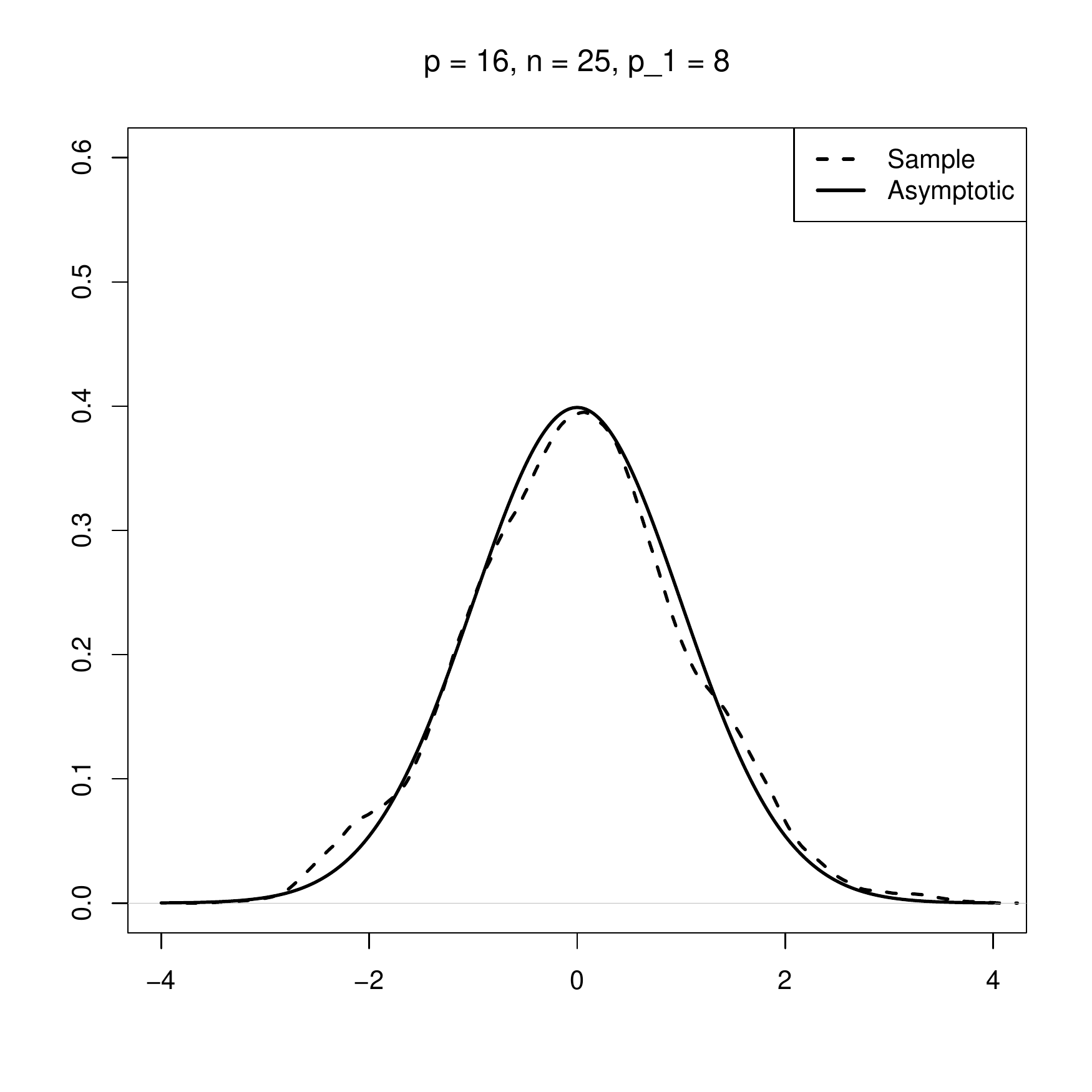} ~~
  \includegraphics[scale=0.222]{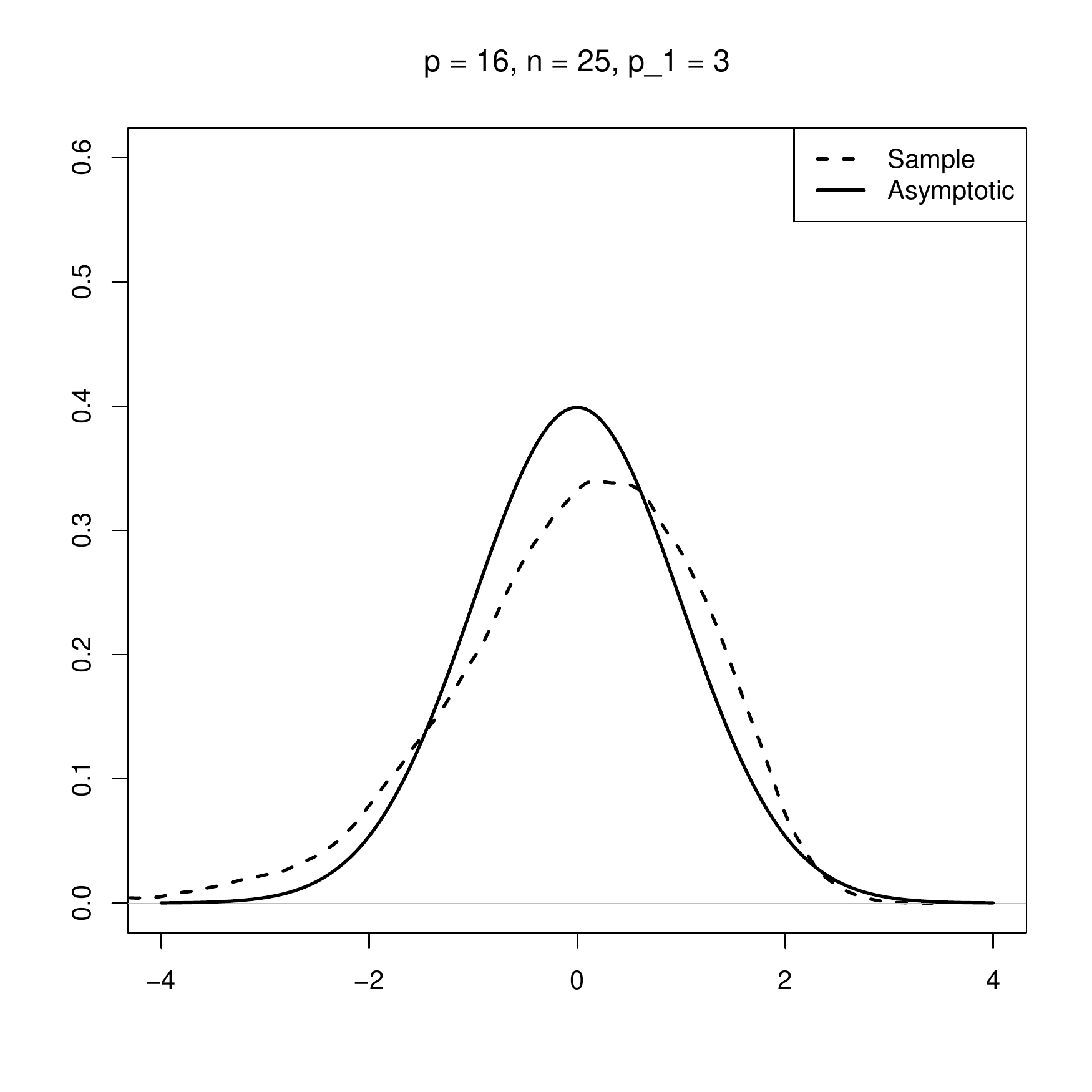}\\

\vspace{-.5cm}
  \caption{\it Simulated distribution of the statistic $(T_{BNP} -(p-p_1)s_{BNP}  - \mu_{BNP}) /  \sigma_{BNP} $ und the null hypothesis
 for  sample size $n=25$,  dimension $p=16$ and  various values of $p_1=13, 8, 3$. The solid curve shows the
  standard normal distribution.
  \label{fig61}}
\end{figure}

In order to investigate the properties of two adjusted tests $T_{BNP}$ and $T_{LH}$ for small dimensions and small sample sizes we provide additional results for $p=16$, $n=25$ and different values of $p_1=13, 8, 3$. The results are depicted in Figures \ref{fig51}  and \ref{fig61} and indicate a good approximation of the
nominal level although a small-sample effect is present. Note that this effect is more pronounced for the LH test as for the  BNP. 
Thus the results are still reliable and there is again no large bias as in case of LR and Wilk's statistics when the dimension $p_1$ is much smaller than $p-p_1$.

\section{Distributional properties under alternative hypothesis}\label{sec4}
\def\theequation{4.\arabic{equation}}
\setcounter{equation}{0}

In this section we derive the distribution of the considered linear spectral statistics under the alternative hypothesis.
The main difficulty consists in the fact that under the alternative the random matrix $\bW\bT^{-1}$
has a (conditionally) non-central Fisher distribution in this case.

The following two  results, which are  proved in the Appendix and of  independent  interest, specify  the asymptotic distribution of the empirical spectral distribution of the matrix $\bW\bT^{-1}$ under $H_{1}$. Throughout the paper
\begin{eqnarray*}
  m_Q(z) = \int_{-\infty}^{+\infty}\frac{dQ(t)}{t-z}\,
\end{eqnarray*}
denotes the Stieltjes transform of a distribution function $Q$.

\begin{theorem}\label{lsdH1}
Consider the  alternative hypothesis $H_{1}$ in \eqref{hyp_H0} and
assume that the assumptions of Section \ref{sec2} are satisfied.
If  the the matrix $\bR= \bSigma_{22\cdot 1}^{-1/2}\bSigma_{21}\bSigma_{11}^{-1}\bSigma_{12}\bSigma_{22\cdot 1}^{-1/2}$ is bounded in spectral norm and its spectral distribution  converges weakly to some function $G$, then for any $z \in  \C \setminus \er$
   the Stieltjes transform of the empirical spectral distribution of the matrix
   $\bW\bT^{-1}$ converges almost surely to some deterministic function  $s$, which is the unique solution of the following system of equations
   {\footnotesize
  \begin{eqnarray}
 && \frac{s(z)}{1+\gamma_2zs(z)}=  m_{H}\left(z(1+\gamma_2zs(z))\right),\label{s(z)}  \nonumber \\[5mm]
  && \frac{ m_H(z)}{1+\gamma_1m_H(z)}= m_{\tilde{H}}((1+\gamma_1m_H(z))[(1+\gamma_1m_H(z))z-(1-\gamma_1)]),\label{H}\\[5mm]
 &&m_{\tilde{H}}(z)(1-(c-c_1)-(c-c_1)zm_{\tilde{H}}(z))c_1^{-1}= m_G\left(\frac{c_1z}{1-(c-c_1)-(c-c_1)zm_{\tilde{H}}(z)} \right)\label{tilde_H}\,,
  \end{eqnarray}
  }
subject to the condition that $\Im \{s(z)\}$ is of the same sign as $\Im\{z\}$. The functions $H$ and $\tilde{H}$ denote the limiting spectral distributions of the matrices $\bW$ and $\tilde{\bR}=1/p_1\bSigma_{22\cdot 1}^{-1/2}\bSigma_{21}\bSigma_{11}^{-1}\bS_{11}\bSigma_{11}^{-1}\bSigma_{12} \bSigma_{22\cdot 1}^{-1/2}$, respectively.
  \end{theorem}

\medskip

\noindent
Note that the matrix $\tilde{\bR}$ from Theorem \ref{lsdH1} can be interpreted as the sample covariance matrix generated from a population with the covariance matrix equal to $\frac{p_1}{n}\bR$.

We will use this result to derive a CLT for the linear spectral statistics of the matrix $\bW\bT^{-1}$, which can be used
for the analysis of the test proposed in Section \ref{sec3} under the alternative hypothesis.
For this purpose we introduce   some  useful notations as follows
\begin{eqnarray}
 \delta(z) &=& \gamma_1m_{H}(z)\label{dz}\\
\tilde{\delta}(z) &=& \delta(z)-\frac{1-\gamma_1}{z}\label{dtildez}  \nonumber\\
  \eta(z) &=& (1+\delta(z))(1+\tilde{\delta}(z))\label{etaz}\\
 \xi(z) &=& \frac{\delta'(z)}{(z\eta(z))'}\label{xiz}\\
%\tilde{\xi}(z) &=& \xi(z) - \frac{1-\gamma_1}{z^2(z\eta(z))^\top}\\
%\Delta(z) &=&  \left(\frac{1}{1+\delta(z)} - \frac{z\eta(z)}{(1+\delta(z))^2}\xi(z)\right)^2-z^2\xi(z)\tilde{\xi}(z)\\
\Psi(z) &=& \left(\frac{1}{1+\delta(z)}-2\xi(z)z+\frac{1-\gamma_1}{1+\delta(z)}\xi(z)\right)^{-1}\,, \label{Psiz}\\
r&=& 2\frac{(1+\sqrt{\gamma_1})^2+\lambda_{max}(\bR){(1+\sqrt{c_1})^2}}{(1-\sqrt{\gamma_2})^2}  \label{rtilde}
\end{eqnarray}

\noindent
\begin{theorem}\label{cltH1}
  If the  assumptions of Theorem \ref{lsdH1} are satisfied, then for any pair $ f, g$ of  analytic functions  in an open region of the complex plane containing the interval $[0, r]$
 the random vector
\begin{equation}  \nonumber
\Big ((p-p_1)\int_0^\infty f(x) \diff (F_n(x)-F^*_n(x)),  ~ (p-p_1)\int_0^\infty g(x) \diff (F_n(x)-F^*_n(x))   \Big )^\top
\end{equation}
converges weakly to a Gaussian vector $(X_{f}, X_{g})^\top$ with mean and covariances given by
{\footnotesize
\begin{eqnarray}
  \E[X_{f}] &=&\frac{1}{4\pi i}\oint f(z) \diff\log(q(z)) + \frac{1}{2\pi i}\oint f(z)B(zb(z)) \diff(zb(z))\nonumber\\
            &+& \frac{1}{2\pi i}\oint f(z)\theta_{b,H}(z)\nonumber\\
  &\times&\left(\theta_{\tilde{b},\tilde{H}}(zb(z))\frac{c^2_1\int \underline{m}^3_{\tilde{H}}(zb(z))t^2(c_1+t\underline{m}_{\tilde{H}}(zb(z)))^{-3}dG(t)}
{(1-c_1\int\underline{m}^2_{\tilde{H}}(zb(z))t^2(c_1+t\underline{m}_{\tilde{H}}(zb(z)))^{-2}dG(t))^2}\right) \diff z
   \nonumber \\[0.2cm]
&&  \label{asym_mean}\\
  \Cov[X_{f}, X_{g}] &=&-\frac{1}{2 \pi^2} \oint \oint f(z_1)g(z_2)\frac{\partial^2 \log(z_1b(z_1)-z_2b(z_2))}{\partial z_1 \partial z_2}\diff z_1 \diff z_2\nonumber\\
             &-&\frac{1}{2 \pi^2} \oint \oint  f(z_1)g(z_2)\frac{\partial^2\log(z_1b(z_1)\eta(z_1b(z_1))-z_2b(z_2)\eta(z_2b(z_2)))}{\partial z_1 \partial z_1}\diff z_1 \diff z_2\nonumber\\
            &-&\frac{1}{2 \pi^2} \oint \oint  f(z_1)g(z_2)\nonumber\\
  &\times&\left[ \theta_{\tilde{b},\tilde{H}}(z_1b(z_1))\theta_{\tilde{b},\tilde{H}}(z_2b(z_2))\left(\frac{\partial^2 \log\left[\frac{\underline{m}_{\tilde{H}}(z_2b(z_2))-\underline{m}_{\tilde{H}}(z_1b(z_1))}{(z_2b(z_2)-z_1b(z_1))}\right]}{\partial z_1 \partial z_2} \right)\right]\diff z_1 \diff z_2\nonumber\\ \label{asym_cov}\,
\end{eqnarray}
}
respectively, where
{\small
  \begin{eqnarray}  
    b(z)&=& 1+\gamma_2zs(z)\nonumber\\
     \tilde{b}(z) &=& 1+\gamma_1m_{H}(z)\label{btildez}\\
    q(z)& =&  1- \gamma_2\int \frac{b^2(z)dH(t)}{(t/z-b(z))^2}\nonumber
  \end{eqnarray}
  {\footnotesize
  \begin{eqnarray*}
 % \nonumber   \theta_{b, H}(z)&=& \frac{b(z)}{1-\gamma_{2}z m_{H}(zb(z))-b(z)\gamma_{2}z^2\int\frac{dH(t)}{(t-zb(z))^2}}\\
    \nonumber  &&\theta_{\tilde{b},\tilde{H}}(z)\\
                               &=&\frac{\tilde{b}(z)}{1-\gamma_{1} m_{\tilde{H}}\left(\tilde{b}(z)(\tilde{b}(z)z-(1-\gamma_{1}))\right)-\tilde{b}(z)\gamma_{1}(2z\tilde{b}(z)-(1-\gamma_{1}))\int\frac{d\tilde{H}(t)}{\left[t-\left(\tilde{b}(z)(\tilde{b}(z)z-(1-\gamma_{1}))\right)\right]^2}}
  \end{eqnarray*}}
  \begin{eqnarray}
 \nonumber    \underline{m}_{\tilde{H}}(z)&=&-\frac{1-c_1}{z}+c_1m_{\tilde{H}}(z)\\
  B(z) &=&\Psi^2(z) \left[ -\widetilde{\omega}(z)N(z)(1-\delta(z)) + \frac{1}{1+\delta(z)}N(z)+\xi(z)\Psi^{-1}(z)+z\xi^2(z)\right.\nonumber\\
        &+&  \left.   z^2\tilde{\delta}^2(z)\left( \xi^2(z) - \delta(z)N(z)\left(z-\frac{1-\gamma_1}{1+\delta(z)}+1\right)\right)     \right]\label{B}\,
\end{eqnarray}
}with
{\small
\begin{eqnarray} \nonumber
N(z) = \frac{\xi'(z)\Psi^{-1}(z) }{2} - \xi^2(z)~~~\text{and}~~~\widetilde{\omega}(z)=z^2\xi(z) + \frac{1-\gamma_1}{1+\delta(z)}\Psi^{-1}(z)\,.
\end{eqnarray}
}
Here the integrals are taken over an arbitrary positively oriented contour which contains the interval $[0 ,r]$, moreover the contours in \eqref{asym_cov} are non-overlapping.
\end{theorem}

There are substantial  differences between the  CLT  derived here  and  the recent results  in \citet{zheng2017}. In particular the matrix $\bW$ does not possess the usual properties of the covariance matrix under normality anymore. Indeed, the conditional distribution of $\bW$ given  $\bS_{11}$ is 
a non-central Wishart distribution, while the unconditional distribution is defined by a very complicated integral expression. 
As a consequence  $\bW\bT^{-1}$ can be interpreted as a conditionally non-central Fisher matrix, while \citet{zheng2017} considered a rescaled Fisher matrix. In general, the CLT presented in \citet{zheng2017} is constructed for studying the asymptotic power of the test for the equality of two population covariance matrices. In contrast, the CLT derived in Theorem \ref{cltH1} is used to investigate the power of the test for block-diagonality, i.e., $H_0: \bSigma_{12}=\mathbf{O}$.

\medskip

\noindent
It follows from the proof of  Theorem \ref{lsdH1} that
\begin{eqnarray}\label{tail_ineq}
   \bW\stackrel{d}{\leq } 2\big(\frac{1}{p_1}\bX\bX^\top + \bM\bM^\top \big ),
  \end{eqnarray} 
where $n\bM\bM^\top\sim\mathcal{W}_{p-p_1}(n, \bR)$ and all entries of $\bX$ are independent and standard normally distributed. Consequently the  largest eigenvalue of the matrix $\bW$ will  asymptotically
be smaller than $$2\left((1+\sqrt{\gamma_1})^2+\lambda_{max}(\bR){(1+\sqrt{c_1})^2}\right)$$ and the quantity  $r$ defined in \eqref{rtilde} is an upper bound for the limiting spectrum of the matrix $\bW \bT^{-1}$.

This observation is quite important for controlling the tail estimates of the extreme eigenvalues of the matrix $\bW\bT^{-1}$, which play a vital role for the application of the Cauchy's integral formula \eqref{cauchy} at the end of the proof of Theorem \ref{cltH1}.  The proof of the following result 
is given in the appendix.

  \begin{proposition}\label{upper_bound}  Let $l_r>r$, where $r$ is given in \eqref{rtilde}, then  
\begin{eqnarray*}
    \forall k\in \mathbbm{N}:\quad \mathbbm{P}(\lambda_{max}(\bW\bT^{-1})>l_r)=o\bigl(n^{-k}\bigr)\,.
  \end{eqnarray*}
  \end{proposition}

 Although, the limiting mean and variance presented in Theorem \ref{cltH1} are very difficult to calculate in a closed form even for simple  cases, there are several important implications of Theorem \ref{cltH1}.
 \begin{remark}[Eigenvectors]
\rm Going through the proof of Theorem \ref{cltH1} one can see that Lemma \ref{Taras} in Section \ref{sec5} reveals an interesting though quite expected fact that the resulting asymptotic distributions depend neither on the eigenvectors of the non-centrality matrix $\bOmega_1$ nor on the eigenvectors of the  matrix $\bR=\bSigma_{22\cdot 1}^{-1/2}\bSigma_{21}\bSigma_{11}^{-1}\bSigma_{12}\bSigma_{22\cdot 1}^{-1/2}$  for the normally distributed data. Loosely speaking, without loss of generality (w.l.g.), we can restrict ourselves  to the case when $\bOmega_1$ and $\bR$ are diagonal matrices, which simplifies the simulations in a remarkable way. 
 \end{remark}
 \begin{remark}[Generalizations and simplifications]
\rm   The non-central Fisher matrix in our case arises only conditionally on $\bS_{11}$ where the non-centrality matrix $\bOmega_1$ is random in our framework. As a consequence  Theorem \ref{cltH1} generalizes the result of \cite{yao2013},  where a deterministic non-centrality matrix was considered. Moreover, all the asymptotic quantities including $\delta(z)$ are expressed  in a more convenient form, like, $\delta (z)=\gamma_1m_H(z)$.
Finally, the expression of the bias term $B(z)$ is significantly simplified which makes it possible to do numerical computations more efficiently and to investigate the results of Theorem \ref{cltH1} deeper in the future.
\end{remark}
\begin{remark}[Finite rank alternatives] \label{rem3}
\rm   Combining Theorem \ref{lsdH1} and Theorem \ref{cltH1} one observes
that finite rank alternatives with a bounded spectrum have no influence on the asymptotic power of the tests,
because the asymptotic means and variances under the null hypothesis and  alternative hypothesis coincide. Indeed, assuming that the matrix $\bR$ has a finite rank, say  $k$, and  a bounded spectrum we get
  \begin{eqnarray*}
    m_{F^{\bR}}(z)&=&\int \frac{dF^{\bR}(t)}{t-z} = \frac{1}{p-p_1}\sum\limits_{i=1}^{p-p_1}\frac{1}{\lambda_i(\bR)-z}\\
    &=&\frac{1}{p-p_1}\sum\limits_{i=1}^k\frac{1}{\lambda_i(\bR)-z}-\frac{p-p_1-k}{p-p_1}\frac{1}{z}\to -\frac{1}{z}\,.
  \end{eqnarray*}
  Thus, it follows  that $m_G(z)=-\frac{1}{z}$, and therefore $G$ is the distribution function of the Dirac measure concentrated at the point $0$.
Consequently we obtain  $\underline{m}_{\tilde{H}}(z)=-1/z$ and  
   the third summands in \eqref{asym_mean} and in \eqref{asym_cov}  vanish, that  is 
  \begin{eqnarray*}
   && \int \frac{t^2}{(c_1+t\underline{m}_{\tilde{H}}(z))^{3}}dG(t)=\int \frac{t^2}{(c_1+t\underline{m}_{\tilde{H}}(z))^{3}}\delta_0(t)d(t)=0,\\
    && \frac{\partial^2 \log\left[\frac{\underline{m}_{\tilde{H}}(z_2)-\underline{m}_{\tilde{H}}(z_1)}{z_2-z_1}\right]}{\partial z_1\partial z_2}=\frac{\underline{m}'_{\tilde{H}}(z_1)\underline{m}'_{\tilde{H}}(z_2)}{(\underline{m}_{\tilde{H}}(z_1)-\underline{m}_{\tilde{H}}(z_2))^2}-\frac{1}{(z_1-z_2)^2} 
    =0\,, 
  \end{eqnarray*}
  for any $z, z_1, z_2 \in\mathbbm{C}^+$. The other summands in \eqref{asym_mean} and in \eqref{asym_cov} do not depend  
  on the eigenvalues of matrix $\bR$, which reflects the alternative hypothesis $H_1$ via $\bSigma_{12}$, thus, they are expected to be equal to the corresponding quantities under the null hypothesis $H_0$ given in Theorem \ref{th2}. 
  %  Showing this claim precisely is not an easy task, so we check it via simulations.
Consequently, all tests based on a linear spectral statistic cannot detect the alternative hypothesis $H_1$ if the matrix $\bR$ has no large eigenvalues. 

On the other hand, if $\lambda_{max}(\bR)$  is an increasing function of the dimension $p-p_1$ the 
spectrum of $\lambda_{max}(\bR)$ is not bounded and Theorem \ref{cltH1} is not applicable. Although we have no theoretical result in this case
we expect that the power of the tests will be an increasing function of  $\lambda_{max}(\bR)$. These properties have been verified numerically by means of a simulation study.
\end{remark}

 \begin{remark}[Full rank alternatives]
\rm   As we have already mentioned, the formulas in Theorem \ref{lsdH1} and Theorem \ref{cltH1} are very complex, which makes
it difficult  to calculate the power functions of the considered tests in an analytic form. For instance, we need to solve the system of three equations in Theorem \ref{lsdH1} which leads to the cubic equation already for $m_H(z)$  even in the simple case $\bR=\rho^2\bI$. On the other hand, the whole system in Theorem \ref{lsdH1} simplifies to a quadratic equation under the null hypothesis $H_0$. Nevertheless, we believe that these results may be useful for future investigations of the power of the considered tests on the block diagonality of the covariance matrix. For example, one may consider the numerical approximations discussed in \cite{zheng2017}.
 \end{remark}

\begin{figure}[H]
  \centering
  \includegraphics[scale=0.222]{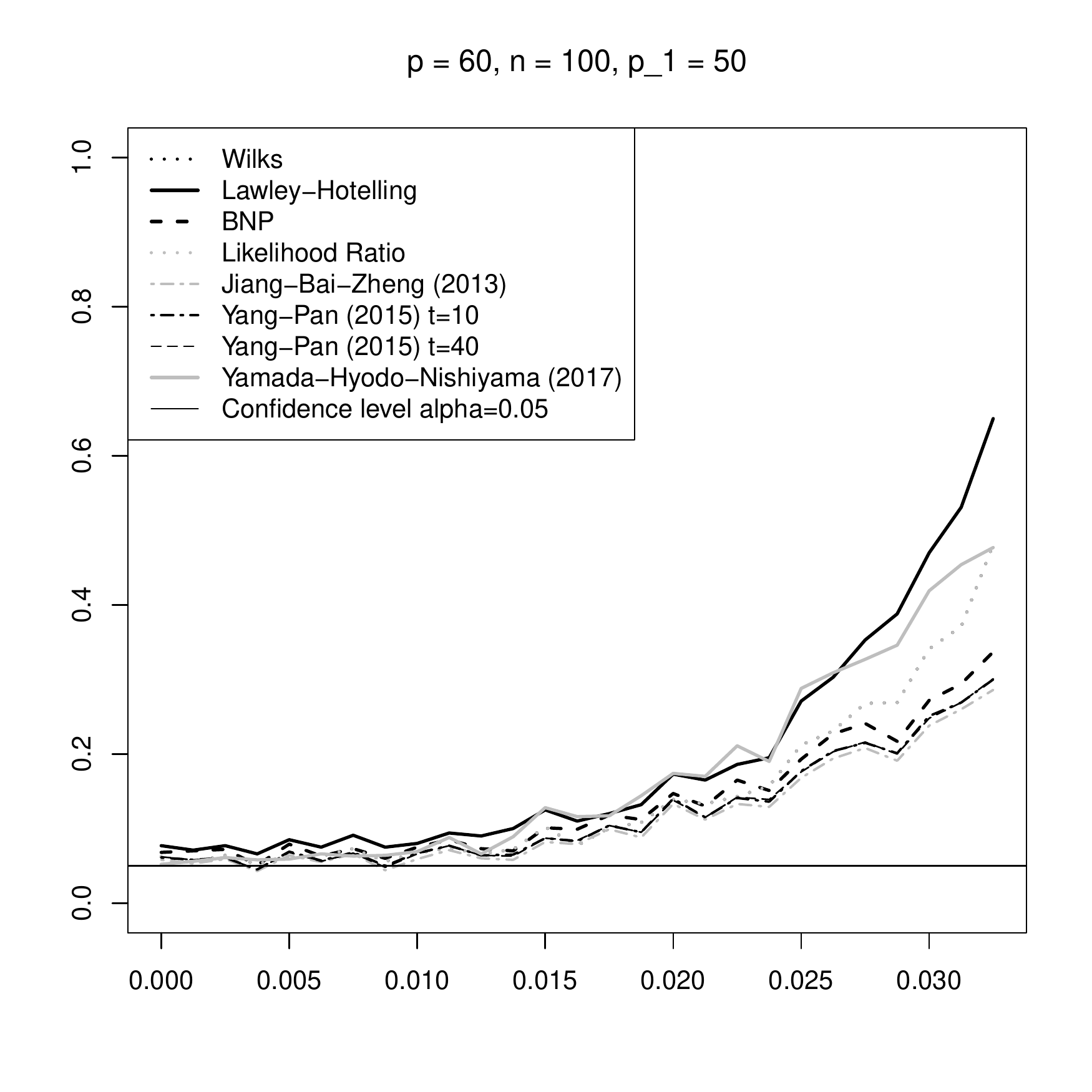} ~~
  \includegraphics[scale=0.222]{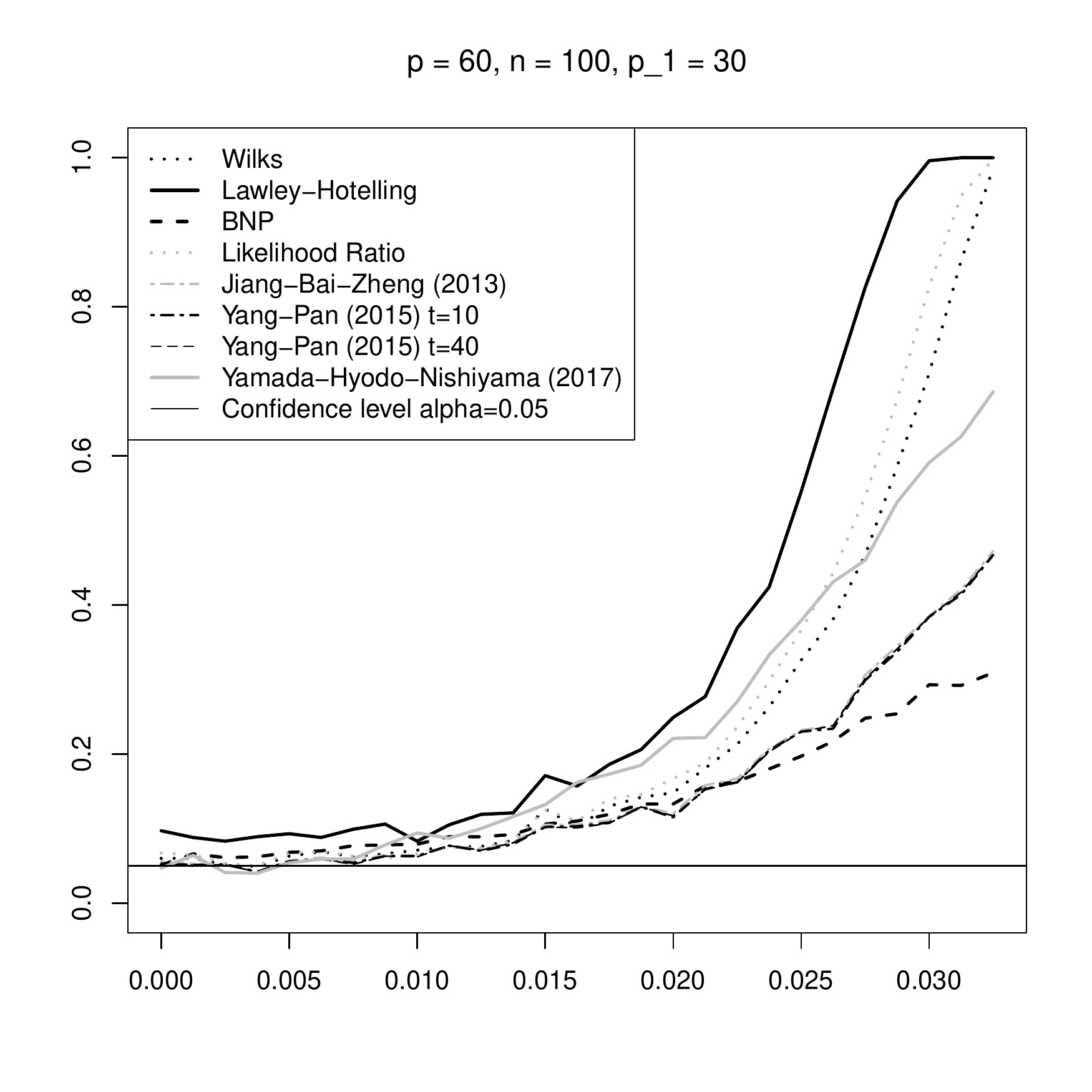} ~~
  \includegraphics[scale=0.222]{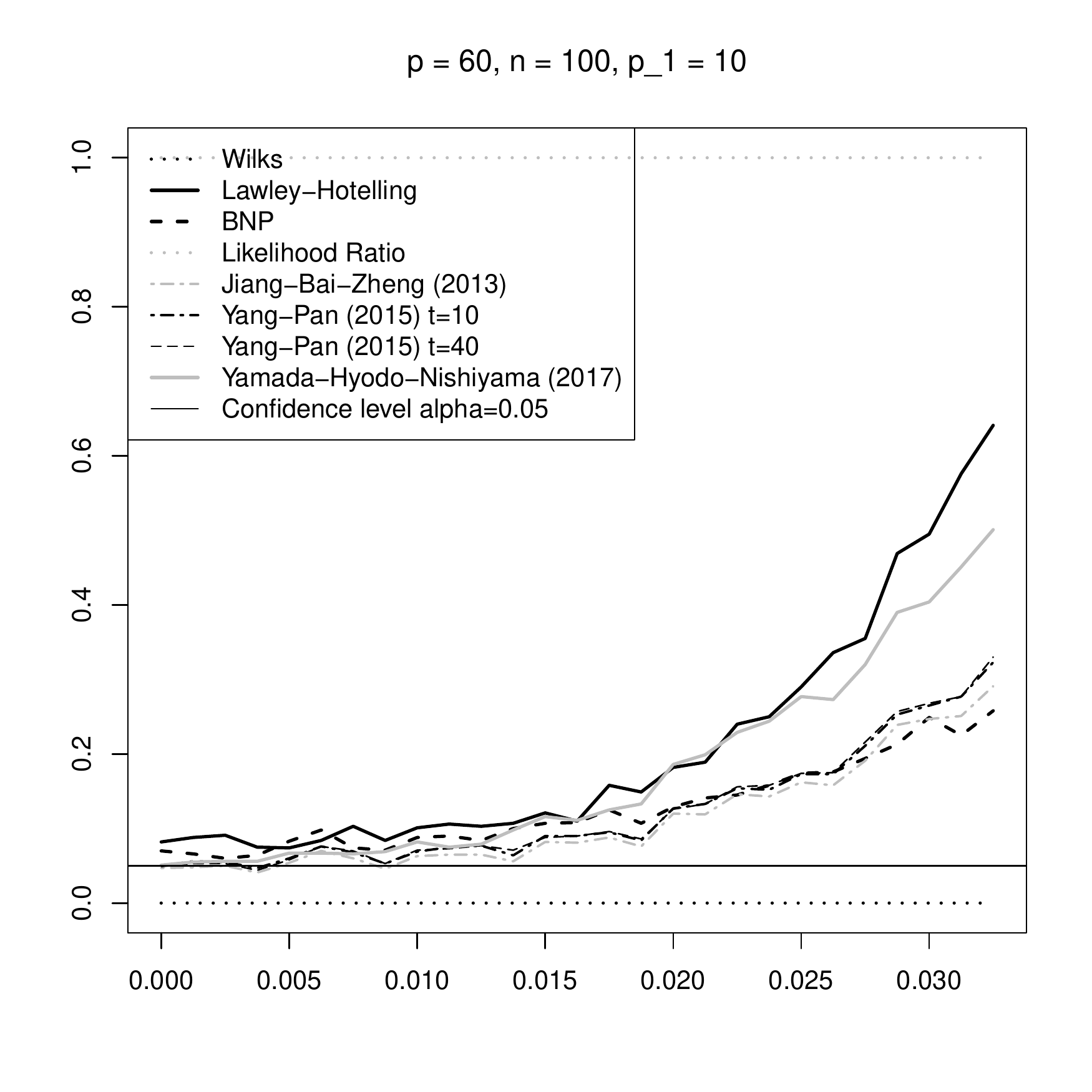}\\

\vspace{-.5cm}
  \caption{\it Empirical power  of  different tests for block diagonality   for  sample size $n=100$,  dimension $p=60$ and  various values of $p_1=50, 30, 10$ as a function of the correlation coefficient $\rho=\frac{\sigma_{12}}{\sigma}$ in $[0, 0.0325]$.
  \label{fig5}}
\end{figure}
To illustrate these remarks and comments, we present a comparison of  the power of the  different tests under consideration by means of
 a small simulation study. In order to demonstrate the results in a clear way we assume for simplicity that $\bSigma_{11}=\bSigma_{22}=\sigma\bI$
 which yields
  $$
  \bR=\frac{1}{\sigma}(\sigma\bI-\frac{1}{\sigma}\bSigma_{21}\bSigma_{12})^{-1/2}\bSigma_{21}\bSigma_{12}(\sigma\bI-\frac{1}{\sigma}\bSigma_{21}\bSigma_{12})^{-1/2}.
  $$
  Note that the spectrum of matrix $\bR$ is the same as that of  the matrix  $\bSigma_{21}\bSigma_{12}(\sigma^2\bI-\bSigma_{21}\bSigma_{12})^{-1}$. First, we take $\bSigma_{12}$ as a rank  $1$ matrix with all components equal to $\sigma_{12}\in[0, 1.3]$ (equicorrelation) and in order to assure positive definiteness of $\bSigma$ in that range we choose $\sigma=40$. Note that if $\sigma_{12}$ varies in the interval  $[0, 1.3]$ the correlation coefficient $\rho=\sigma_{12}/\sigma$ will change in the interval  $[0, 0.0325]$. 
  Further, we increase the rank of $\bSigma_{12}$ by setting some of its elements to zero (sparsifying). The empirical rejection probailities of  the  proposed  tests in the case of rank  $1$ alternatives are given in Figure \ref{fig5}.
  
   For the sake of comparison  we also included the trace criterion  recently proposed  by \cite{jiangetal2013}, the test  introduced by \cite{yamadaetal2017}, which is based on  an empirical distance between the full and a block diagonal covariance matrix; and  the test suggested 
  by  \cite{yang2015} built on the sum of the canonical correlations coefficients. Note that there exists a regularized and a non-regularized  version of 
  the latter test. In general, the statistic of \cite{yang2015} is defined by the sum of eigenvalues of the matrix $(\bS_{22}+t\bI_{p-p_1})^{-1}\bS_{21}\bS^{-1}_{11}\bS_{12}$ for some $t\geq 0$. Thus, in case $t=0$ this test is equivalent to the sum of canonical correlation coefficients. Moreover, 
   the test of \cite{jiangetal2013} and \cite{yang2015} coincide for $t=0$ because  the matrix $\bS^{-1}_{22}\bS_{21}\bS^{-1}_{11}\bS_{12}$ can be written in form $\bW\bT^{-1}(\frac{\gamma_1}{\gamma_2}\bI +\bW\bT^{-1})^{-1}$ under $H_0$ (see, e.g., \cite{yaobaizheng2015} Section 8.5.1). Thus, in order to visualize difference between them we take $t=10$ and $t=40$ for the test proposed by \cite{yang2015}. Further, the simulations showed that taking larger $t$ will lead only to a slight increase of power in the case where $p_1$  is not equal to $p-p_1$. Nevertheless, it must be mentioned that the regularized test
    is applicable even in the case $p-p_1>n$ while all of other considered tests need both $p_1<n$ and $p-p_1<n$.

   Figure \ref{fig5} justifies our theoretical findings, i.e., none of the tests can detect the alternatives for $\rho \in [0, 0.01]$ (the power function in this region is basically flat and close to the nominal level $0.05$). On the other hand, if the correlation is greater than $0.01$, then all of the tests gain power. For $p_1=30$ (case of equal blocks) all test are powerful enough to reject $H_0$ if the correlation is greater than $0.03$. 
These results are in accordance with the discussion  in Remark \ref{rem3},  because 
in the considered scenario the largest eigenvalue of the matrix $\bR$  is given by 
$$\frac{p_1(p-p_1)\rho^2}{1-p_1(p-p_1)\rho^2}.
$$
 Thus, if the correlation coefficient $\rho$ is close to ${1}/{\sqrt{p_1(p-p_1)}}$  we will get a spike (note that  ${1}/{\sqrt{p_1(p-p_1)}}\approx 0.0333$ if $p_1=30$, $p=60$).
 Moreover, here we have a clear winner - the Lawley-Hotelling's (LH) trace criterion. The test of  \cite{yamadaetal2017} and  Wilk's test with the corrected likelihood ratio (LR) criterion are ranked on the second  and third, respectively. The regularized test of \cite{yang2015} together with the trace criterion of \cite{jiangetal2013} are on the fourth position, while the the Bartlett-Nanda-Pillai's (BNP) trace criterion shows the worst performance. Interestingly, the tests  \cite{jiangetal2013} and \cite{yang2015} can not be visually distinguished neither for $t=10$ nor for $t=40$.

A similar ranking was observed for $p_1=50$ with the difference
 of a decreasing power of all tests and a slight increase of power for \cite{yang2015} with respect to its benchmark for $t=0$, i.e., \cite{jiangetal2013}. Note that  Wilk's test and the LR test have the same power for $p_1=50$. In light of the previous findings obtained under the null hypothesis $H_0$, the case $p_1=10$ is  the 
 most interesting one. Indeed, here we observe that  Wilk's and the LR test are not reliable anymore (they either always (Wilk's) or never (LR) reject $H_0$).
 On the other hand  the  other tests show a similar  behaviour as in the case $p_1=50$. As before, the LH test is the most powerful in all three situations.

\begin{figure}[H]
  \centering
  \includegraphics[scale=0.222]{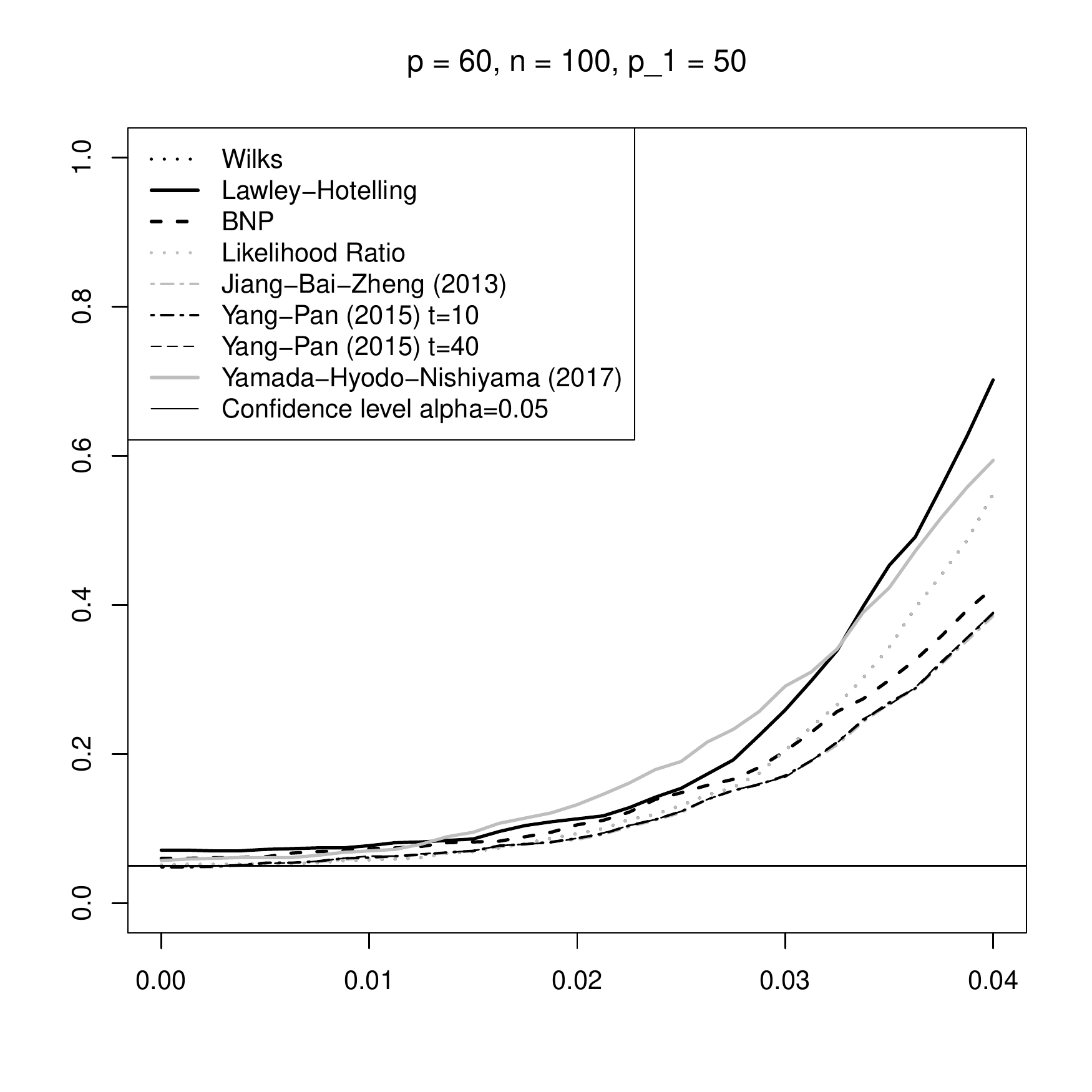} ~~
  \includegraphics[scale=0.222]{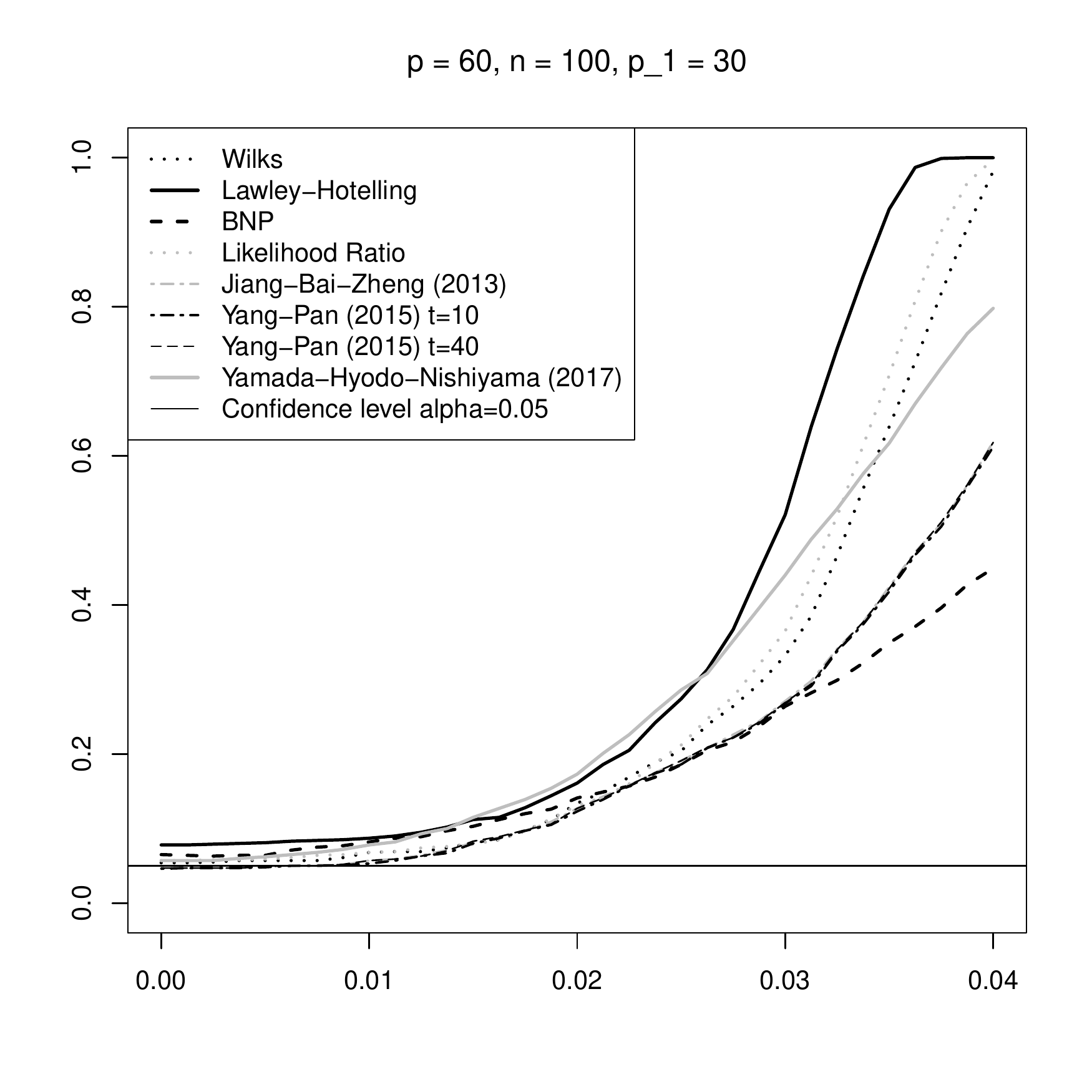} ~~
  \includegraphics[scale=0.222]{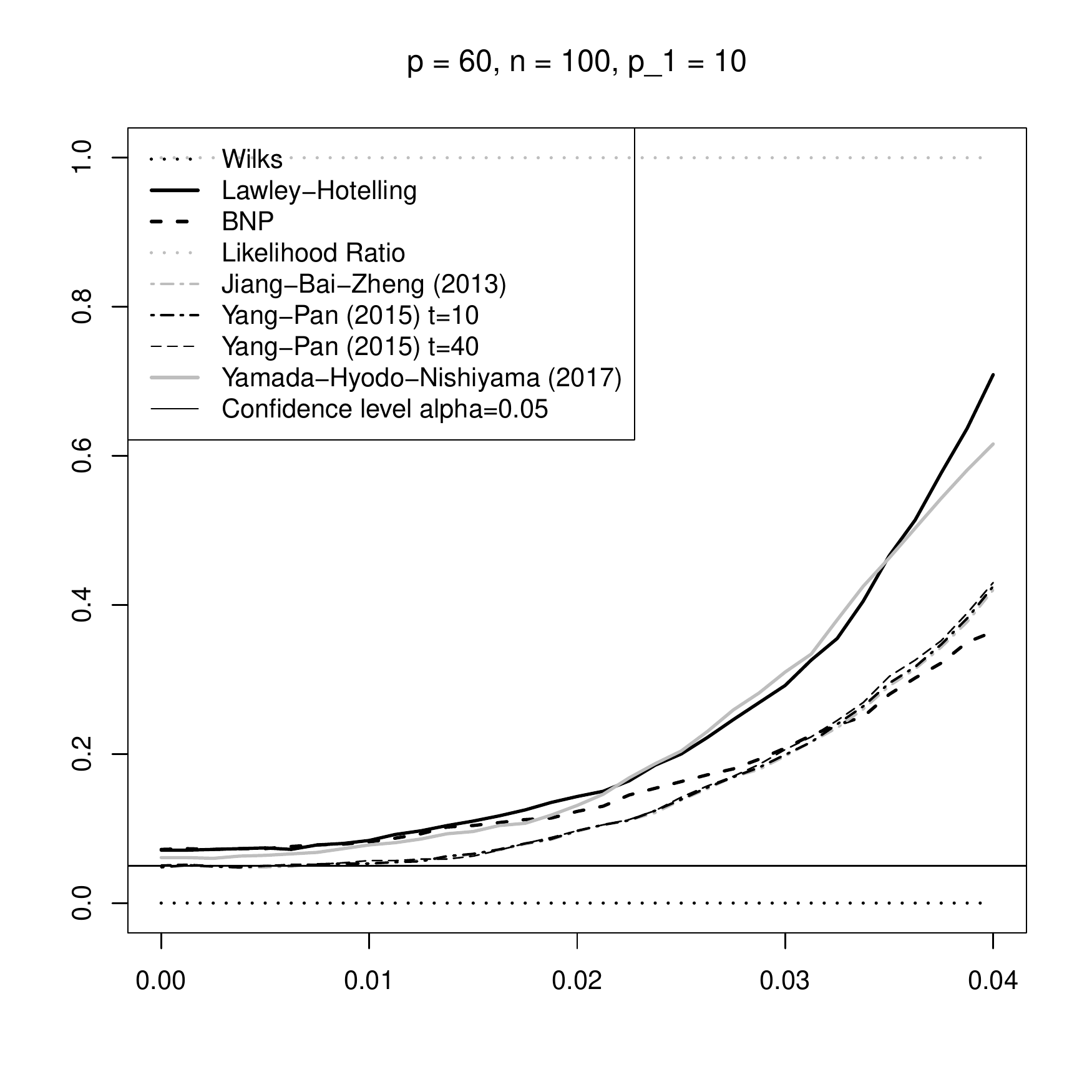}\\
\vspace{-.5cm}
  \caption{\it Empirical power of  different tests for block diagonality for  sample size $n=100$,  dimension $p=60$ and  various values of $p_1=50, 30, 10$ as a function of the correlation coefficient $\rho=\frac{\sigma_{12}}{\sigma}$ in $[0, 0.04]$ and sparsity level of $20\%$.
  \label{fig6}}
\end{figure}
\begin{figure}[H]
  \centering
  \includegraphics[scale=0.222]{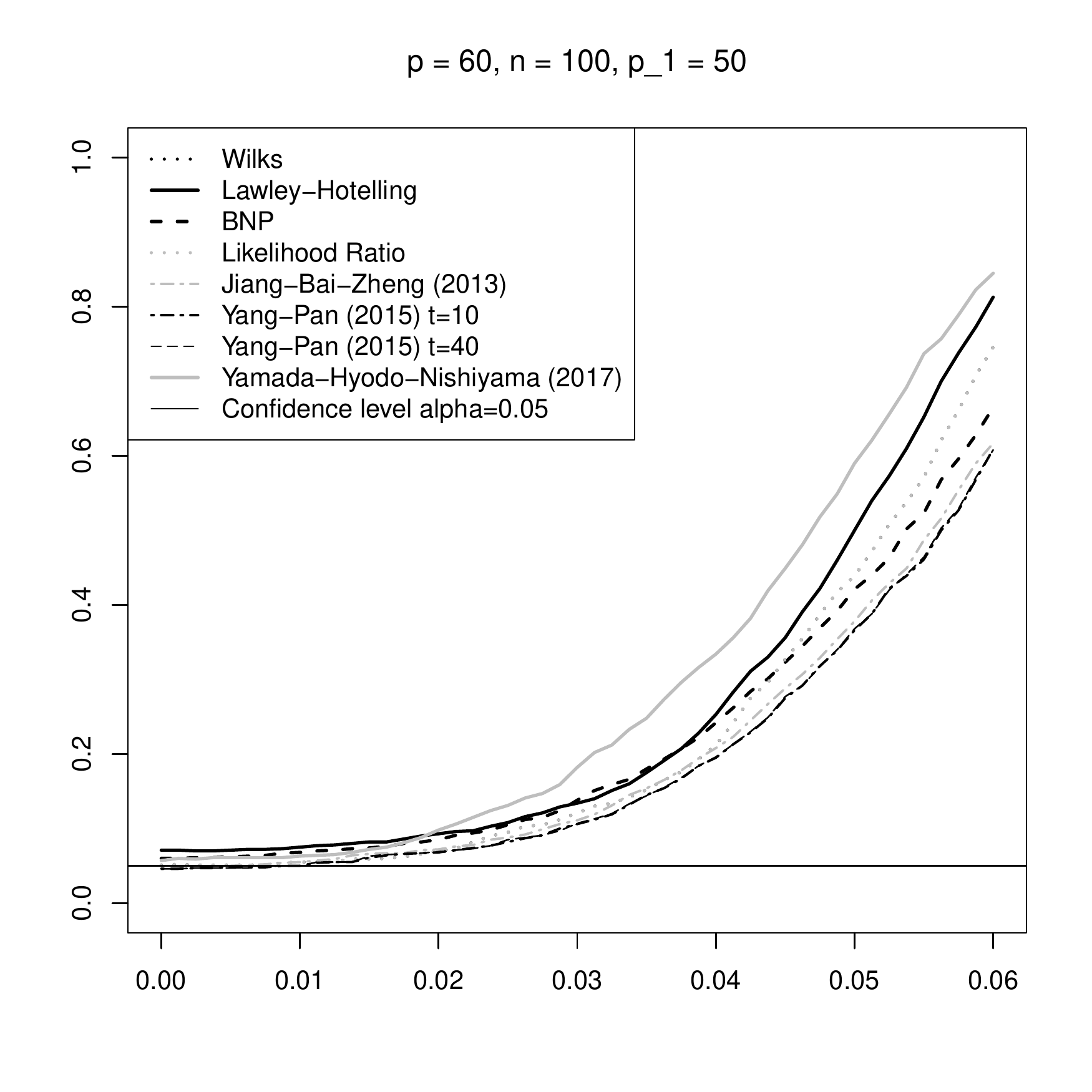} ~~
  \includegraphics[scale=0.222]{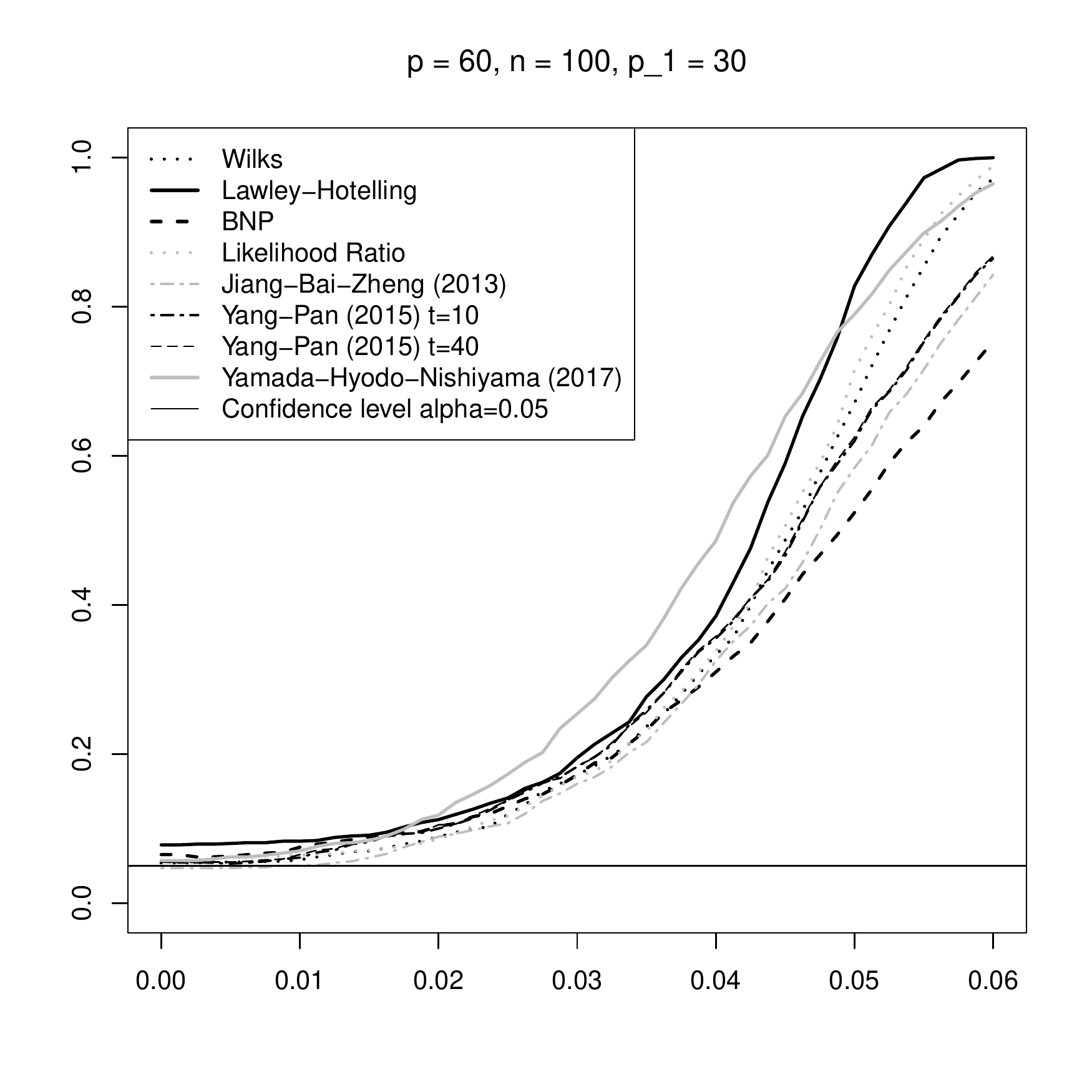} ~~
  \includegraphics[scale=0.222]{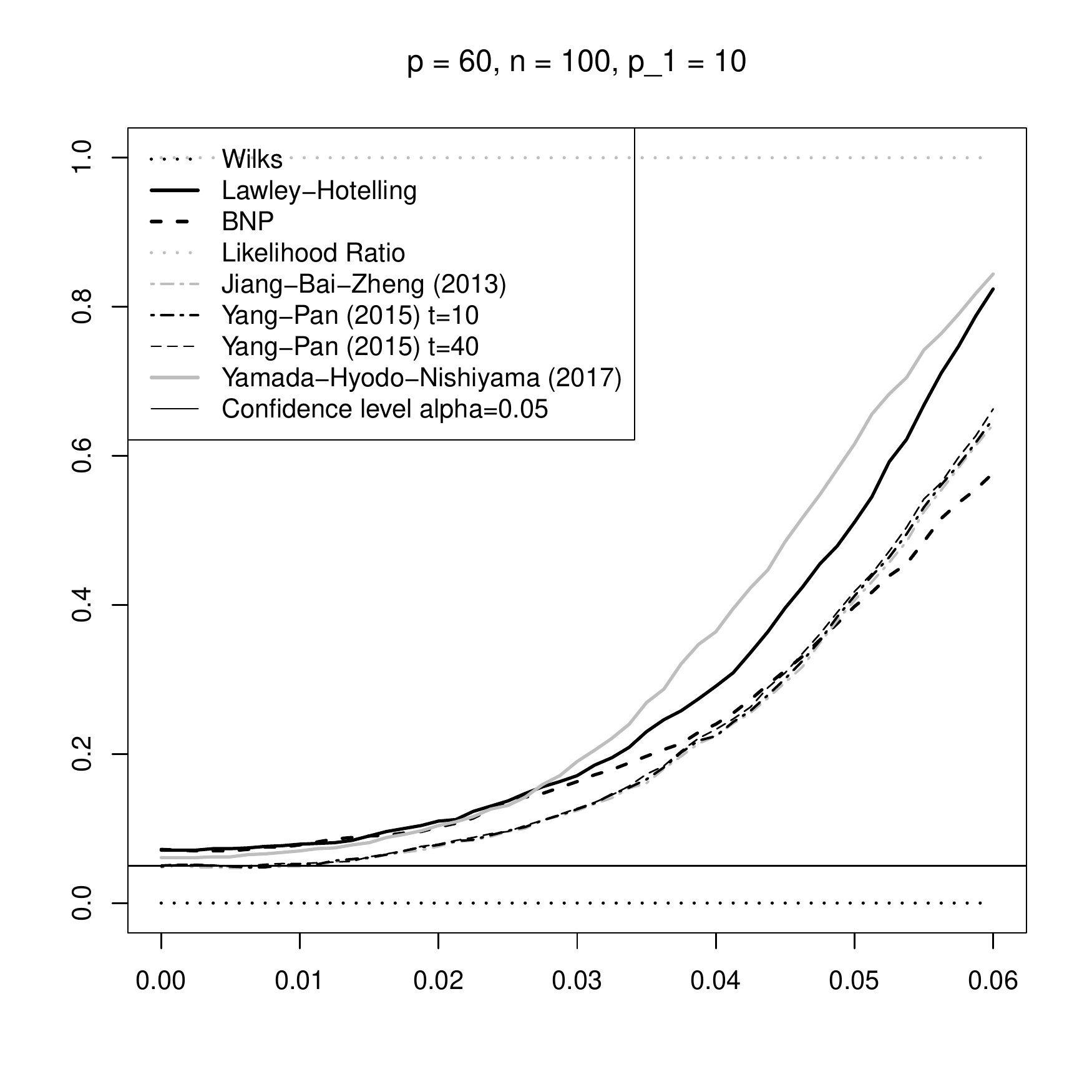}\\
\vspace{-.5cm}
  \caption{\it Empirical power  of  different tests for block diagonality  for  sample size $n=100$,  dimension $p=60$ and  various values of $p_1=50, 30, 10$ as a function of the correlation coefficient $\rho=\frac{\sigma_{12}}{\sigma}$ in $[0, 0.06]$ and sparsity level of $50\%$.
  \label{fig7}}
\end{figure}

In order to investigate the robustness of the tests we increase the sparsity of the matrix $\bSigma_{12}$, where  $20\%$ and $50\%$ of the  elements are set randomly to zero, while all other elements are still equal to $\sigma_{12}$. By this procedure we increase the probability that $\bSigma_{12}$ has full rank. The results are summarized in Figures \ref{fig6} and \ref{fig7}.

We observe a similar behaviour  as in the non-sparse case (see Figure \ref{fig5}).  The LH test and the test proposed in
 \cite{yamadaetal2017}
show the best performance. The latter is slightly  better  than the LH test for the sparsity level $50\%$, while a superiority of the LH test could be observed for
a sparsity level of  $20\%$. 
Of course, by increasing the sparsity level we make the alternative hypothesis harder to detect. For this reason the non-sensitivity interval $[0, 0.01]$ (the interval where the test is not sensitive to the alternative $H_1$) is increased to $[0, 0.02]$ and $[0, 0.03]$ in case of $20\%$ and $50\%$ sparsity levels, respectively.

Moreover, in the supplementary material  (see, \cite{supplement}) we have also investigated the performance of the different tests for $p=10$, $100$, $200$, $250$, $290$ and $n=300$ (see, Figures  \ref{fig10} - \ref{fig14} in Appendix B of the supplementary material). Our findings still remain unchanged except for the  case $p=290$: here all tests have a substantial 
loss in  power  and the  LH test does  not keep the nominal  level. Indeed, this is an expected result,
because in the  case $p=290$ and $n=300$ the ratio $p/n$ is close  to one,
and  the sample covariance matrix is a very unstable estimator.

As a conclusion, although the LH trace criterion is the most simple one among the linear spectral statistics of  the matrix
$\bW\bT^{-1}$ ($f=id$),  it seems to be the most robust and powerful test on the block diagonality of the large-dimensional covariance matrix. On the other hand the corrected LR and Wilk's criteria can not be recommended, if the size of the first block is much smaller than the size of the second one.

  \begin{remark}\rm
 A possible reason for the superior performance of the LH test 
 are the specific  alternatives  considered in our numerical experiments. In particular results  
 coincide with the findings in  \citet{pillai1967}, where the LH test showed the best performance in the case, where  the eigenvalues of the matrix 
 $\bSigma_{22}^{-1}\bSigma_{21}\bSigma_{11}^{-1}\bSigma_{12}$ are far apart. Thus, taking very sparse or other full rank alternatives could have 
 an considerable impact on the dominance property of the LH test.
\end{remark}

\bigskip

\noindent
{\bf Acknowledgements.} The authors would like to thank Professor Tailen Hsing and two anonymous referees for the valuable comments and remarks which have significantly improved the previous version of this paper. Additionally, we want to thank M. Stein who typed parts of this manuscript with considerable technical
expertise. The work of H. Dette was supported by the DFG Research Unit 1735, DE 502/26-2.

%%%%%%%%%%%%%%%%%%%%%%%%%%%%%%%%%%%%%%%%%%%%%%%%%%%%%%%%%%%%%%%%%%%%%%%%%%%%%%%%%%%%%%%%%%%
% {\small
%  \bibliography{blockdiag}
% }

%%%%%%%%%%%%%%%%%%%%%%%%%%%%%%%%%%%%%%%%%%%%%%%%%%%%%%%%%%%%%%%%%%%%%%%%%%%%%%%%%%%%%%%%%
\section{Appendix: Proofs} \label{sec5}
\def\theequation{5.\arabic{equation}}
\setcounter{equation}{0}
 
\noindent 
 \begin{proof}[\bf Proof of Proposition \ref{upper_bound}:]
    Because $\bW$ and $\bT$ are positive semi-definite we have
    {\small
  \begin{eqnarray}\label{ineq_1}
    \mathbbm{P}(\lambda_{max}(\bW\bT^{-1})>l_r)&{\leq}& \mathbbm{P}(\lambda_{max}(\bW)\lambda_{max}(\bT^{-1})>l_r)\nonumber\\
    &\leq&  \mathbbm{P}\left(\lambda_{max}(\bW)> l_r(1-\sqrt{\gamma_2})^2\right)+o\big(n^{-k}\big)\,,
     \end{eqnarray}}
     where the last inequality  follows from inequality (1.9b) in \citet{baisil2004}. Furthermore, using 
     this inequality again and \eqref{tail_ineq} we get{\small
     \begin{eqnarray}\label{ineq_2}
      &&  \mathbbm{P}\left(\lambda_{max}(\bW)> l_r(1-\sqrt{\gamma_2})^2\right)\\
      &{\leq}& \mathbbm{P}\Big (2\lambda_{max}\big(1/p_1\bX\bX^\top\big) + 2\lambda_{max}\big(\bM\bM^\top \big )> l_r(1-\sqrt{\gamma_2})^2\Big )\nonumber\\
      &\leq& \mathbbm{P}\Big (\lambda_{max}\big(1/p_1\bX\bX^\top\big)> \frac{l_r}{2}(1-\sqrt{\gamma_2})^2-\lambda_{max}(\bR)(1+\sqrt{c_1})^2\Big )+o\big(n^{-k}\big).\nonumber
     \end{eqnarray}}
%     where the last inequality is  a consequence of inequality (1.9a)  in \citet{baisil2004}. 
  Finally, combining  \eqref{ineq_1}, \eqref{ineq_2} and using \eqref{rtilde} with $l_r>r$ we arrive at
{\small
     \begin{eqnarray*}
       \mathbbm{P}(\lambda_{max}(\bW\bT^{-1})>l_r)\leq \mathbbm{P}\Big (\lambda_{max}\big(1/p_1\bX\bX^\top\big)> (1+\sqrt{\gamma_1})^2\Big )+o\big(n^{-k}\big)=o\big(n^{-k}\big),
     \end{eqnarray*}
     }
where the last equality follows again from (1.9a) of \citet{baisil2004}.
\end{proof}

\begin{proof}[\bf Proof of Theorem \ref{lsdH1}:]
Since $(n-p_1) \bT \sim W_{p-p_1}(n-p_1,\bSigma_{22\cdot 1}) $, $ p_1 \bW|\bS_{11} \sim W_{p-p_1}(p_1, \bSigma_{22\cdot 1},\bOmega_1) $ with
 $\bOmega_1=\bOmega_1(\bS_{11})=\bSigma_{22\cdot 1}^{-1}\bSigma_{21}\bSigma_{11}^{-1}\bS_{11}\bSigma_{11}^{-1}\bSigma_{12} = \bSigma_{22\cdot 1}^{-1}\bM\bM^\top$, $\bS_{11}\sim W_{p_1}(n, \bSigma_{11})$, and $\bT$ is independent of $\bW$ and $\bS_{11}$ we get the following stochastic representations for $\bT$ and $\bW$ expressed as\footnote{Here we use the definition of the non-central Wishart distribution given by \cite[Definition 10.3.1]{muirhead1982}.}
\begin{eqnarray*}
  \bW &\overset{d}{=}&\frac{1}{p_1}\bSigma^{1/2}_{22\cdot 1}(\bX+\bSigma^{-1/2}_{22\cdot 1}\bM)(\bX+\bSigma^{-1/2}_{22\cdot 1}\bM)^\top\bSigma^{1/2}_{22\cdot 1}\\
  \bT&\overset{d}{=}& \frac{1}{n-p_1}\bSigma^{1/2}_{22\cdot 1}\bY\bY^\top\bSigma^{1/2}_{22\cdot 1}\,,
\end{eqnarray*}
where $\bX\sim \mathcal{N}_{p-p_1, p_1}(\mathbf{O},\bI\otimes\bI)$, $\bY\sim \mathcal{N}_{p-p_1, n-p_1}(\mathbf{O},\bI\otimes \bI)$, and $\bX$, $\bY$, $\bS_{11}$ are mutually independent. Then, the stochastic representation of $\bW\bT^{-1}$ is given by
\begin{eqnarray*}
  \bW\bT^{-1} &\overset{d}{=}& \frac{1}{p_1}\bSigma^{1/2}_{22\cdot 1}(\bX+\bSigma^{-1/2}_{22\cdot 1}\bM)(\bX+\bSigma^{-1/2}_{22\cdot 1}\bM)^\top\bSigma^{1/2}_{22\cdot 1}\\
  &\times&\left(\frac{1}{n-p_1}\bSigma^{1/2}_{22\cdot 1}\bY\bY^\top\bSigma^{1/2}_{22\cdot 1}\right)^{-1}\,.
\end{eqnarray*}
The last equality in distribution implies that the spectral distribution of $\bW\bT^{-1}$ is the same as the spectral distribution of $\widetilde{\bW}\widetilde{\bT}^{-1}$ with
\begin{equation}\label{tbW} \nonumber
  \widetilde{\bW}=\frac{1}{p_1}(\bX+\bSigma^{-1/2}_{22\cdot 1}\bM)(\bX+\bSigma^{-1/2}_{22\cdot 1}\bM)^\top
  ~~\text{and}~~ \widetilde{\bT}=\frac{1}{n-p_1}\bY\bY^\top.
\end{equation}

From Theorem 2.1 of \cite{zhengbaiyao2015} it holds that the Stieltjes transform of $\widetilde{\bW}\widetilde{\bT}^{-1}$ given $\widetilde{\bW}$ $m_{F^{\widetilde{\bW}\widetilde{\bT}^{-1}}|\widetilde{\bW}}(z)$ converges to $s_{\widetilde{\bW}}(z)$ which satisfies the following equation
\begin{eqnarray}\label{stWT}
  zs_{\widetilde{\bW}}(z) = -1 + \int\frac{tdH(t)}{t-z(1+\gamma_2zs_{\widetilde{\bW}}(z))}\,,
\end{eqnarray}
where $H(t)=H_{\widetilde{\bW}}(t)$ is the limiting spectral distribution of the matrix $\widetilde{\bW}$, which is a deterministic function following Theorem 1.1 of \cite{dozsil2007}. Noting that the right hand-side of \eqref{stWT} does not depend on the condition $\widetilde{\bW}$ and rewriting \eqref{stWT}, we get the limiting spectral distribution of $\bW\bT^{-1}$, which is equal to $\widetilde{\bW}\widetilde{\bT}^{-1}$, is given by $s(z)=s_{\widetilde{\bW}}(z)$ expressed as
\begin{eqnarray*}
  zs(z)& =& \int\frac{z\gamma_2(zs(z)+1)dH(t)}{t-z(1+\gamma_2zs(z))} = z(\gamma_2zs(z)+1)m_{H}(z(\gamma_2zs(z)+1)),
\end{eqnarray*}
where (see Theorem 1.1 of \cite{dozsil2007})
\begin{eqnarray*}
    m_H(z)& =& \int \frac{(1+\gamma_1m_H(z))d\tilde{H}(t)}{t-(1+\gamma_1m_H(z))[(1+\gamma_1m_H(z))z-(1-\gamma_1)]}\nonumber\\[2mm]
& =& (1+\gamma_1m_H(z))m_{\tilde{H}}((1+\gamma_1m_H(z))[(1+\gamma_1m_H(z))z-(1-\gamma_1)])
\end{eqnarray*}
with $\tilde{H}$ the limiting spectral distribution of
\begin{eqnarray*}
  \widetilde{\bR}&=&1/p_1\bSigma_{22\cdot 1}^{-1/2}\bSigma_{21}\bSigma_{11}^{-1}\bS_{11}\bSigma_{11}^{-1}\bSigma_{12} \bSigma_{22\cdot 1}^{-1/2}\\
  &=&c^{-1}_{1,n}1/n\bSigma_{22\cdot 1}^{-1/2}\bSigma_{21}\bSigma_{11}^{-1}\bS_{11}\bSigma_{11}^{-1}\bSigma_{12} \bSigma_{22\cdot 1}^{-1/2}.   
\end{eqnarray*}
satisfying the following equation
\begin{eqnarray*}
 m_{\tilde{H}}(z) &=&  \int \frac{(1-(c-c_1)-(c-c_1)zm_{\tilde{H}}(z))^{-1}dG(t)}{c^{-1}_1t-\frac{z}{1-(c-c_1)-(c-c_1)zm_{\tilde{H}}(z)}}\nonumber\\
                  &=& c^{-1}_1(1-(c-c_1)-(c-c_1)zm_{\tilde{H}}(z))^{-1}\nonumber\\
  &\times&m_G\left(\frac{c_1z}{1-(c-c_1)-(c-c_1)zm_{\tilde{H}}(z)} \right)
\end{eqnarray*}
where $G(t)$ is the limiting spectral distribution of the matrix\\ $\bR =\bSigma_{22\cdot 1}^{-1/2}\bSigma_{21}\bSigma_{11}^{-1}\bSigma_{12}\bSigma_{22\cdot 1}^{-1/2}$ which is deterministic as well.
\end{proof}

In the proof of Theorem \ref{cltH1} we make use of the following lemma which simplifies the conditions used in Theorem 2.2 of \cite{zhengbaiyao2015}
and is proved in  Appendix B of the supplementary material (see, \cite{supplement}).

%%%%%%%%%%%%%

\begin{lemma}\label{Taras} Conditionally on $\bS_{11}$ the distribution of the matrix $\bW\bT^{-1}$ solely depends on the eigenvalues of the non-centrality matrix $\bOmega_1(\bS_{11})$ and does not depend on the corresponding eigenvectors. Moreover, the unconditional distribution of the eigenvalues of matrix $\bW\bT^{-1}$ depends only on the eigenvalues of the matrix $\widetilde{\bR}=\bSigma_{11}^{-1}\bSigma_{12}\bSigma_{22\cdot 1}^{-1}\bSigma_{21}$.
\end{lemma}
%%%%%%%%%%
The results of Lemma \ref{Taras} shows that both the unconditional distribution of the eigenvalues of $\bW\bT^{-1}$ and its conditional distribution given $\bS_{11}$ depend only on the eigenvalues of $\bOmega_1(\bS_{11})$ and of $\widetilde{\bR}=\bSigma_{22\cdot 1}^{-1/2}\bSigma_{21}\bSigma_{11}^{-1}\bS_{11}\bSigma_{11}^{-1}\bSigma_{12} \bSigma_{22\cdot 1}^{-1/2}$, respectively, for any fixed dimension $p$ and sample size $n$. Consequently, without loss of generality both matrices $\bOmega_1(\bS_{11})$ and of $\widetilde{\bR}$ can be taken as diagonal. These simplify the validation of the conditions present in Theorem 2.2.1 and Theorem 2.2.2 of \cite{yao2013}.

%The last two terms in \ref{omega} and \ref{zeta} are in general dependent on the eigenvectors of the matrix $\bOmega_1(\bS_{11})$. Fortunately, we show that in case of normally distributed data the dependence on the eigenvectors is actually not present even in the finite dimensional case for any fixed sample size $n$. For that purpose we derive the joint finite-sample distribution of $v_1 \ge v_2 \ge ... \ge v_{p-p_1}$ under alternative $H_1$.}

\begin{proof}[{\bf Proof of Theorem \ref{cltH1}:}]

Throughout the proof of Theorem \ref{cltH1}, we assume that the complex number $z$ belongs to the arbitrary positively oriented contour $\mathcal{C}$, which contains the limiting support $[0, r]$. We consider
\begin{eqnarray}\label{main}
  &&(p-p_1)\left(m_{F^{\bW\bT^{-1}}}(z) - s_n(z)\right)\\
  &=& (p-p_1)\left(m_{F^{\bW\bT^{-1}}}(z) - s^*_n(z)\right)+ (p-p_1)\left(s^*_n(z) - s_n(z)\right)\nonumber\,,
\end{eqnarray}
where $s_n(z)$ and $s_n^*(z)$ are unique roots of the following equations
\begin{eqnarray}
  zs_n(z)  = -1 + \int\frac{tdH_n(t)}{t-z(1+\gamma_{2,n}z s_n(z))}\label{eqn1}\\
  zs^*_n(z) =    -1 + \int\frac{tdF^{\bW}_n(t)}{t-z(1+\gamma_{2,n}z s^*_n(z))}\label{eqn2}\,
\end{eqnarray}
with $\gamma_{2,n}=\frac{p-p_1}{n-p_1}$. The symbol $H_n$ denotes the discretized limiting distribution of $\bW$ with $\gamma_2$ replaced by $\gamma_{2,n}$ and $F^{\bW}_n$ stands for the empirical spectral distribution of $\bW$.
%\begin{eqnarray}\label{step1}
%  (p-p_1)\left(m_{F^{\bW\bT^{-1}}}(z) - s^*_n(z)\right)
%\end{eqnarray}

Following the proof of Theorem 2.2 by \cite{zhengbaiyao2015}, we get that the first summand $ (p-p_1)\left(m_{F^{\bW\bT^{-1}}}(z) - s^*_n(z)\right)$ in \eqref{main} conditionally on the matrix $\bW$ converges to a Gaussian process $M_1(z)$ with the mean function
\begin{eqnarray}\label{clt1_mean}
  \E(M_1(z)) = \frac{\gamma_2b^3(z)}{z^2q^2(z)}\int\frac{tdH(t)}{(t/z-b(z))^3}=\frac{1}{2}(log(q(z))'
\end{eqnarray}
and the covariance function
\begin{eqnarray}\label{clt1_cov}
  Cov(M_1(z_1), M_1(z_2)) &=& 2\frac{(z_1b(z_1))'(z_2b(z_2))'}{(z_1b(z_1)-z_2b(z_2))^2}\\
  &=&2\frac{\partial log((z_1b(z_1)-z_2b(z_2)))}{\partial z_1\partial z_2}\nonumber\,,
\end{eqnarray}
where
\begin{eqnarray}\label{clt1_cov_a} \nonumber
 b(z) &=& 1+\gamma_2zs(z),\\
 q(z)& =&  1- \gamma_2\int \frac{b^2(z)dH(t)}{(t/z-b(z))^2}\,\label{clt1_cov_b}
\end{eqnarray}
for $z_1$ and $z_2$ from $\mathcal{C}$. Since all quantities in \eqref{clt1_mean}-\eqref{clt1_cov_b} do not depend on the condition $\bW$, we get that this is also the unconditional distribution and both summands in \eqref{main} are independent.

Next, we derive the asymptotic distribution of the second summand $(p-p_1)\left(s^*_n(z) - s_n(z)\right)$ in \eqref{main}. Let
\begin{eqnarray}\label{b_n} \nonumber
 b_n^*(z) = 1+\gamma_{2,n}zs^*_n(z) \quad \text{and} \quad b_n(z) = 1+\gamma_{2,n}zs_n(z)\,.
\end{eqnarray}
Then, by using the definition of the Stieltjes transform, \eqref{eqn1}, and \eqref{eqn2} we get
\begin{eqnarray*}
  &&(p-p_1)(s_n^*(z) - s_n(z)) \\
  &=& (p-p_1)\left(b^*_n(z) m_{F_n^{\bW}}(zb_n^*(z)) - b_n(z)m_{H_n}(zb_n(z))\right)\\
 &=& (p-p_1)(b^*_n(z)-b_n(z))m_{F_n^{\bW}}(zb_n^*(z)) \\
  &+& (p-p_1)b_n(z)(m_{F_n^{\bW}}(zb_n^*(z))-m_{F_n^{\bW}}(zb_n(z)) )\\
 &+& (p-p_1)b_n(z)(m_{F_n^{\bW}}(zb_n(z))-m_{H_n}(zb_n(z)))\\
  &=& (p-p_1)\gamma_{2,n}z(s_n^*(z) - s_n(z))m_{F_n^{\bW}}(zb_n^*(z)) \\
  &+& (p-p_1)b_n(z)\gamma_{2,n}z^2(s_n^*(z) - s_n(z)) \int\frac{dF_n^W(t)}{(t-zb_n^*(z))(t-zb_n(z))}\\
  &+& (p-p_1)b_n(z)(m_{F_n^{\bW}}(zb_n(z))-m_{H_n}(zb_n(z))) \,.
\end{eqnarray*}
Hence,
\begin{eqnarray*}
 &&(p-p_1)(s_n^*(z) - s_n(z)) =
 (p-p_1)(m_{F_n^{\bW}}(zb_n(z))-m_{H_n}(zb_n(z)))\\
 &\times&\frac{b_n(z)}{1-\gamma_{2,n}z m_{F_n^{\bW}}(zb_n^*(z))-b_n(z)\gamma_{2,n}z^2\int\frac{dF_n^W(t)}{(t-zb_n^*(z))(t-zb_n(z))}}\,,
\end{eqnarray*}
where

\begin{eqnarray}\label{eq1_as} \nonumber 
&&\frac{b_n(z)}{1-\gamma_{2,n}z m_{F_n^{\bW}}(zb_n^*(z))-b_n(z)\gamma_{2,n}z^2\int\frac{dF_n^W(t)}{(t-zb_n^*(z))(t-zb_n(z))}}\\ \nonumber
&\stackrel{a.s.}{\to}&\theta_{b, H}(z)= \frac{b(z)}{1-\gamma_{2}z m_{H}(zb(z))-b(z)\gamma_{2}z^2\int\frac{dH(t)}{(t-zb(z))^2}}= \frac{b^2(z)}{q(z)}\,,
\end{eqnarray}
where the last equality follows from \eqref{clt1_cov_b} and
  \begin{eqnarray}\label{eq_zb}
    \gamma_{2}zb(z) m_{H}(zb(z))=b(z)-1\,.
  \end{eqnarray}

Next, we derive the asymptotic distribution of $(p-p_1)(m_{F_n^{\bW}}(zb_n(z))-m_{H_n}(zb_n(z)))$. It holds that
\begin{eqnarray}\label{main_W}
&&  (p-p_1)(m_{F_n^{\bW}}(zb_n(z))-m_{H_n}(zb_n(z)))\nonumber\\
&=&(p-p_1)(m_{F_n^{\bW}}(zb_n(z))-m_{H_n^{\bS_{11}}}(zb_n(z)))\\
                                                &+&(p-p_1)(m_{H_n^{\bS_{11}}}(zb_n(z))-m_{H_n}(zb_n(z)))\label{S11}
\end{eqnarray}
where $m_{H_n^{\bS_{11}}}(z)$ and $m_{H_n}(z)$ are the unique solutions of the equations
{\small
\begin{eqnarray}
 \frac{m_{H_n^{\bS_{11}}}(z)}{(1+\gamma_{1,n}m_{H_n^{\bS_{11}}}(z))}&=&\int\frac{d{F_n^{\widetilde{\bR}}(t)}}{t-(1+\gamma_{1,n}m_{H_n^{\bS_{11}}}(z))\left[(1+\gamma_{1,n}m_{H_n^{\bS_{11}}}(z))z-(1-\gamma_{1,n})\right]}\nonumber\\
  &&\label{eqn11}\\
  \frac{m_{H_n}(z)}{(1+\gamma_{1,n}m_{H_n}(z))}&=& \int\frac{d\tilde{H}_n(t)}{t-(1+\gamma_{1,n}m_{H_n}(z))\left[(1+\gamma_{1,n}m_{H_n}(z))z-(1-\gamma_{1,n})\right]},\nonumber\\
  &&\label{eqn22}\,
\end{eqnarray}
}
where $\widetilde{\bR}=1/p_1\bSigma_{22\cdot 1}^{-1/2}\bSigma_{21}\bSigma_{11}^{-1}\bS_{11}\bSigma_{11}^{-1}\bSigma_{12} \bSigma_{22\cdot 1}^{-1/2}$, $\tilde{H}_n(t)$ stands for its discretized limiting spectral distribution, and $F_n^{\widetilde{\bR}}(t)$ is the empirical spectral distribution of $\widetilde{\bR}$.

First, we consider the second summand in \eqref{S11}. Let
\begin{eqnarray}\label{tb_n}
 \nonumber 
  \tilde{b}_n^*(z) = 1+\gamma_{1,n} m_{H_n^{\bS_{11}}}(z) \quad \text{and} \quad \tilde{b}_n(z) = 1+\gamma_{1,n}m_{H_n}(z)\,.
\end{eqnarray}
Similarly, using the definition of Stieltjes transform, \eqref{eqn11} and \eqref{eqn22} one can write
{\footnotesize
\begin{eqnarray*}
  &&(p-p_1)(m_{H_n^{\bS_{11}}}(z)-m_{H_n}(z))\\
  &=&  (p-p_1)\tilde{b}^*_n(z) m_{F_n^{\widetilde{\bR}}}\left(\tilde{b}_n^*(z)(\tilde{b}_n^*(z)z-(1-\gamma_{1,n}))\right) \\
  &-& (p-p_1)\tilde{b}_n(z)m_{\tilde{H}_n}\left(\tilde{b}_n(z)(\tilde{b}_n(z)z-(1-\gamma_{1,n}))\right)\nonumber\\
  &=& (p-p_1)(\tilde{b}^*_n(z)-\tilde{b}_n(z))m_{F_n^{\widetilde{\bR}}}\left(\tilde{b}_n^*(z)(\tilde{b}_n^*(z)z-(1-\gamma_{1,n}))\right)\\ &+& (p-p_1)\tilde{b}_n(z)\left[m_{F_n^{\widetilde{\bR}}}\left(\tilde{b}_n^*(z)(\tilde{b}_n^*(z)z-(1-\gamma_{1,n}))\right)-m_{F_n^{\widetilde{\bR}}}\left(\tilde{b}_n(z)(\tilde{b}_n(z)z-(1-\gamma_{1,n}))\right) \right]\\
 &+& (p-p_1)\tilde{b}_n(z)\left[m_{F_n^{\widetilde{\bR}}}\left(\tilde{b}_n(z)(\tilde{b}_n(z)z-(1-\gamma_{1,n}))\right)-m_{\tilde{H}_n}\left(\tilde{b}_n(z)(\tilde{b}_n(z)z-(1-\gamma_{1,n}))\right)\right]\\
  &=& (p-p_1)\gamma_{1,n}(m_{H_n^{\bS_{11}}}(z)-m_{H_n}(z))m_{F_n^{\widetilde{\bR}}}\left(\tilde{b}_n^*(z)(\tilde{b}_n^*(z)z-(1-\gamma_{1,n}))\right) \\
  &+& (p-p_1)\tilde{b}_n(z)\gamma_{1,n}\left(m_{H_n^{\bS_{11}}}(z)-m_{H_n}(z) \right){(z(\tilde{b}^*_n+\tilde{b}_n)- (1-\gamma_{1,n}))}\\ 
  &\times& \int\frac{dF_n^{\widetilde{\bR}}(t)}{\left[t-\left(\tilde{b}_n^*(z)(\tilde{b}_n^*(z)z-(1-\gamma_{1,n}))\right)\right]\left[t-\left(\tilde{b}_n(z)(\tilde{b}_n(z)z-(1-\gamma_{1,n}))\right)\right]}\\
  &+& (p-p_1)\tilde{b}_n(z)\left[m_{F_n^{\widetilde{\bR}}}\left(\tilde{b}_n(z)(\tilde{b}_n(z)z-(1-\gamma_{1,n}))\right)-m_{\tilde{H}_n}\left(\tilde{b}_n(z)(\tilde{b}_n(z)z-(1-\gamma_{1,n}))\right)\right] \,.
\end{eqnarray*}}
Rearranging terms, we get{\footnotesize
 \begin{eqnarray*}
   &&(p-p_1)(m_{H_n^{\bS_{11}}}(z)-m_{H_n}(z)) \\
   &=&(p-p_1)\left[m_{F_n^{\widetilde{\bR}}}\left(\tilde{b}_n(z)(\tilde{b}_n(z)z-(1-\gamma_{1,n}))\right)-m_{\tilde{H}_n}\left(\tilde{b}_n(z)(\tilde{b}_n(z)z-(1-\gamma_{1,n}))\right)\right]\\
   &\times&\tilde{b}_n(z)\Bigg(1-\gamma_{1,n} m_{F_n^{\widetilde{\bR}}}\left(\tilde{b}_n^*(z)(\tilde{b}_n^*(z)z-(1-\gamma_{1,n}))\right)-\tilde{b}_n(z)\gamma_{1,n}{(z(\tilde{b}^*_n+\tilde{b}_n)- (1-\gamma_{1,n}))}\\
   &\times&\int\frac{dF_n^{\widetilde{\bR}}(t)}{\left[t-\left(\tilde{b}_n^*(z)(\tilde{b}_n^*(z)z-(1-\gamma_{1,n}))\right)\right]
   \left[t-\left(\tilde{b}_n(z)(\tilde{b}_n(z)z-(1-\gamma_{1,n}))\right)\right]}\Bigg)^{-1}\,,
 \end{eqnarray*}
 where
\begin{eqnarray*}\label{eq2_as} 
&&\tilde{b}_n(z)\left(1-\gamma_{1,n} m_{F_n^{\widetilde{\bR}}}\left(\tilde{b}_n^*(z)(\tilde{b}_n^*(z)z-(1-\gamma_{1,n}))\right)-\tilde{b}_n(z)\gamma_{1,n}{(z(\tilde{b}^*_n+\tilde{b}_n)- (1-\gamma_{1,n}))} \right.\nonumber\\
  &\times&\left.\int\frac{dF_n^{\widetilde{\bR}}(t)}{\left[t-\left(\tilde{b}_n^*(z)(\tilde{b}_n^*(z)z-(1-\gamma_{1,n}))\right)\right]\left[t-\left(\tilde{b}_n(z)(\tilde{b}_n(z)z-(1-\gamma_{1,n}))\right)\right]}\right)^{-1}\nonumber\\
&\stackrel{a.s.}{\to}&  \theta_{\tilde{b},\tilde{H}}(z)\nonumber\\
  &=&\frac{\tilde{b}(z)}{1-\gamma_{1} m_{\tilde{H}}\left(\tilde{b}(z)[\tilde{b}(z)z-(1-\gamma_{1})]\right)-\tilde{b}(z)\gamma_{1}{(2z\tilde{b}(z)-(1-\gamma_{1}))}
\int\frac{d\tilde{H}(t)}{\left[t-\left(\tilde{b}(z)(\tilde{b}(z)z-(1-\gamma_{1}))\right)\right]^2}},
\end{eqnarray*}}
where $\tilde{b}(z)$ is given in \eqref{btildez}.

The application of Lemma 1.1 in \cite{baisil2004} proves that $(p-p_1)(m_{H_n^{\bS_{11}}}(zb_n(z))-m_{H_n}(zb_n(z)))$ converges to a Gaussian process $M_3(z)$ with the mean function
\begin{eqnarray*}
  \E(M_3(z)) =   \theta_{\tilde{b},\tilde{H}}(zb(z))\frac{c^2_1\int \underline{m}^3_{\tilde{H}}(zb(z))t^2(c_1+t\underline{m}_{\tilde{H}}(zb(z)))^{-3}dG(t)}{(1-c_1\int\underline{m}^2_{\tilde{H}}(zb(z))t^2(c_1+t\underline{m}_{\tilde{H}}(zb(z)))^{-2}dG(t))^2}\,
\end{eqnarray*}
and the covariance function
\begin{eqnarray*}\label{clt3_cov}
&&  Cov(M_3(z_1), M_3(z_2)) = 2\theta_{\tilde{b},\tilde{H}}(z_1b(z_1))\theta_{\tilde{b},\tilde{H}}(z_2b(z_2)) \\
   \nonumber 
&\times &   \left(\frac{\partial}{\partial (z_1b(z_1))}\frac{\underline{m}_{\tilde{H}}(z_1b(z_1))\frac{\partial}{\partial (z_2b(z_2))}\underline{m}_{\tilde{H}}(z_2b(z_2))}{(\underline{m}_{\tilde{H}}(z_1b(z_1))-\underline{m}_{\tilde{H}}(z_2b(z_2)))^2}-\frac{1}{(z_1b(z_1)-z_2b(z_2))^2} \right) \,,
\end{eqnarray*}
where $\underline{m}_{\tilde{H}}(z)=-\frac{1-c_1}{z}+c_1m_{\tilde{H}}(z)$ and $G(t)$ is the limiting spectral distribution of the matrix $\bR=\bSigma_{22\cdot 1}^{-1/2}\bSigma_{21}\bSigma_{11}^{-1}\bSigma_{12} \bSigma_{22\cdot 1}^{-1/2}$.

In order to derive the asymptotic distribution of the first summand in \eqref{main_W}, we use the results in \cite{yao2013} to the conditional distribution of $(p-p_1)(m_{F_n^{\bW}}(zb_n(z))-m_{F_n^{\bS_{11}}}(zb_n(z)))$ given $\bS_{11}$.

From the proof of Theorem \ref{lsdH1}, we know that the empirical spectral distribution of $\bW$ is the same as of $\widetilde{\bW}$ given by
\begin{eqnarray*}
  \widetilde{\bW}=(\frac{1}{\sqrt{p_1}}\bX+\frac{1}{\sqrt{p_1}}\bSigma^{-1/2}_{22\cdot 1}\bM)(\frac{1}{\sqrt{p_1}}\bX+\frac{1}{\sqrt{p_1}}\bSigma^{-1/2}_{22\cdot 1}\bM)^\top\,.
\end{eqnarray*}
with $\bM=\bSigma_{21}\bSigma_{11}^{-1}\bS_{11}^{1/2}$. Furthermore, following Lemma \ref{Taras} it is enough to consider the case where $\bSigma^{-1/2}_{22\cdot 1}\bM\bM^T\bSigma^{-1/2}_{22\cdot 1}$ is diagonal and, consequently, $\bSigma^{-1/2}_{22\cdot 1}\bM$ is pseudo-diagonal.

Finally, in using that $\bX$ consists of i.i.d. entries which are normally distributed and applying the results of Section 2.2.2 in \cite{yao2013}, we get that $(p-p_1)(m_{F_n^{\bW}}(zb_n(z))-m_{F_n^{\bS_{11}}}(zb_n(z)))$ converges to a Gaussian process $M_2(z)$ with the mean function $\E(M_2(z))$ and $Cov(M_2(z_1), M_2(z_2))$ given in the following lemma which is proved below the proof of the theorem.

\begin{lemma}\label{simplify}
The random process $(p-p_1)(m_{F_n^{\bW}}(zb_n(z))-m_{F_n^{\bS_{11}}}(zb_n(z)))$ converges to a Gaussian process $M_2(z)$ with the mean function $\E(M_2(z))$ and $Cov(M_2(z_1), M_2(z_2))$ given by
\begin{eqnarray}\label{clt2_mean_a}  \nonumber 
  \E(M_2(z)) &=& B(zb(z))
\end{eqnarray}
and the covariance function
\begin{eqnarray}\label{clt2_cov_a}  \nonumber 
  Cov(M_2(z_1), M_2(z_2)) =  2\frac{\partial^2\log(z_1b(z_1)\eta(z_1b(z_1))-z_2b(z_2)\eta(z_2b(z_2)))}{\partial(z_1b(z_1))\partial(z_2b(z_2))}
\end{eqnarray}
which are independent of $\bS_{11}$. The functions $B(z)$, $\delta(z)$, $\Psi(z)$, $\xi(z)$ and $\eta(z)$ are given by \eqref{B}, \eqref{dz}, \eqref{Psiz}, \eqref{xiz} and \eqref{etaz}, respectively.
\end{lemma}
The proof of Lemma \ref{simplify} can be found in  Appendix B of the  supplementary material (see, \cite{supplement}).
Thus, merging the results for the independent asymptotic processes $M_2(z)$ and $M_3(z)$, we get
\begin{eqnarray*}
(p-p_1)(s_n^*(z)-s_n(z))\to \theta_{b,H}(z)\left(M_2(z)+M_3(z)\right)  \,,
\end{eqnarray*}
i.e., converges to a Gaussian process with mean and covariance functions given by
{\footnotesize
\begin{eqnarray}
  &&  \theta_{b,H}(z)\left(B(zb(z))+\theta_{\tilde{b},\tilde{H}}(zb(z))\frac{c^2_1\int \underline{m}^3_{\tilde{H}}(zb(z))t^2(c_1+t\underline{m}_{\tilde{H}}(zb(z)))^{-3}dG(t)}
     {(1-c_1\int\underline{m}^2_{\tilde{H}}(zb(z))t^2(c_1+t\underline{m}_{\tilde{H}}(zb(z)))^{-2}dG(t))^2}\right) \nonumber\\
  &&\label{clt2_mean}\\
  &&~~\text{and} \nonumber  \\
  && 2\theta_{b,H}(z_1)\theta_{b,H}(z_2)\left[\frac{\partial^2 \log(z_1b(z_1)\eta(z_1b(z_1))-z_2b(z_2)\eta(z_2b(z_2)))}{\partial(z_1b(z_1))\partial(z_2b(z_2))}+ \theta_{\tilde{b},\tilde{H}}(z_1b(z_1))\theta_{\tilde{b},\tilde{H}}\right.\nonumber\\
 &\times& \left.(z_2b(z_2))\left(\frac{\frac{\partial}{\partial (z_1b(z_1))}\underline{m}_{\tilde{H}}(z_1b(z_1))\frac{\partial}{\partial (z_1b(z_1))}\underline{m}_{\tilde{H}}(z_2b(z_2))}
 {(\underline{m}_{\tilde{H}}(z_1b(z_1))-\underline{m}_{\tilde{H}}(z_2b(z_2)))^2}-\frac{1}{(z_1b(z_1)-z_2b(z_2))^2} \right)\right],\nonumber\\ && \label{clt2_cov}\,
\end{eqnarray}
}
respectively. Remind that $H$ is the asymptotic spectral distribution of the matrix $\bW$ (and, thus, of $\widetilde{\bW}$).
Furthermore, it holds {\footnotesize
  \begin{eqnarray}\label{simple_cov_eta}
  &&  \theta_{b,H}(z_1)\theta_{b,H}(z_2)\frac{\partial^2 \log(z_1b(z_1)\eta(z_1b(z_1))-z_2b(z_2)\eta(z_2b(z_2)))}{\partial(z_1b(z_1))\partial(z_2b(z_2))}\nonumber\\
  &=&\frac{b^2(z_1)}{q(z_1)(z_1b(z_1))'}\frac{b^2(z_2)}{q(z_2)(z_2b(z_2))'}\frac{\partial^2 \log(z_1b(z_1)\eta(z_1b(z_1))-z_2b(z_2)\eta(z_2b(z_2)))}{\partial z_1\partial z_2}\nonumber\\
    &=& \frac{\partial^2 \log(z_1b(z_1)\eta(z_1b(z_1))-z_2b(z_2)\eta(z_2b(z_2)))}{\partial z_1\partial z_2}\,,
  \end{eqnarray}
  where the last equality in \eqref{simple_cov_eta} follows from \eqref{eq_zb} and
  \begin{eqnarray*}
    q(z)(zb(z))'&=&\left(1-\gamma_2(b(z)z)^2\frac{m'_H(zb(z))}{(zb(z))'} \right)(zb(z))'\\
    &=&(zb(z))'-(zb(z))^2(-\frac{1}{z^2}+\frac{(zb(z))'}{(zb(z))^2})=b^2(z)\,.
  \end{eqnarray*}
  Similarly, we get
  {
  \begin{eqnarray}\label{simple_cov_tildeH}
    && \theta_{b,H}(z_1)\theta_{b,H}(z_2)\nonumber\\
    &\times&\left(\frac{\frac{\partial}{\partial z_1b(z_1)} \underline{m}_{\tilde{H}}(z_1b(z_1))\frac{\partial}{\partial z_2b(z_2)}\underline{m}_{\tilde{H}}(z_2b(z_2))}
    {(\underline{m}_{\tilde{H}}(z_1b(z_1))-\underline{m}_{\tilde{H}}(z_2b(z_2)))^2}-\frac{1}{(z_1b(z_1)-z_2b(z_2))^2} \right)\nonumber\\
   &=& \frac{b^2(z_1)}{q(z_1)(z_1b(z_1))'}\frac{b^2(z_2)}{q(z_2)(z_2b(z_2))'}\frac{\partial^2 \log\left(\frac{\underline{m}_{\tilde{H}}(z_1b(z_1))-\underline{m}_{\tilde{H}}(z_2b(z_2))}{z_1b(z_1)-z_2b(z_2) }\right)}{\partial z_1\partial z_2}\nonumber\\
   &=& \frac{\partial^2 \log\left(\frac{\underline{m}_{\tilde{H}}(z_1b(z_1))-\underline{m}_{\tilde{H}}(z_2b(z_2))}{z_1b(z_1)-z_2b(z_2) }\right)}{\partial z_1\partial z_2}\,.
  \end{eqnarray}
}
}

At last, combining the results \eqref{clt1_mean}, \eqref{clt1_cov}, \eqref{clt2_mean}, \eqref{clt2_cov} together with \eqref{simple_cov_eta} and \eqref{simple_cov_tildeH} we get that the process $(p-p_1)\left(m_{F^{\bW\bT^{-1}}}(z) - s_n(z)\right)$ is asymptotically Gaussian with mean and covariance functions given by
{\footnotesize
\begin{eqnarray*}
  && \frac{1}{2}\diff\log(q(z)) +  \theta_{b,H}(z)\\
  &\times&\left(B(zb(z))+\theta_{\tilde{b},\tilde{H}}(zb(z))\frac{c^2_1\int \underline{m}^3_{\tilde{H}}(zb(z))t^2(c_1+t\underline{m}_{\tilde{H}}(zb(z)))^{-3}dG(t)}{(1-c_1\int\underline{m}^2_{\tilde{H}}(zb(z))t^2(c_1+t\underline{m}_{\tilde{H}}(zb(z)))^{-2}dG(t))^2}\right)\\
  &&~~\text{and}\nonumber\\[0.2cm]
 && 2\left[\frac{\partial^2 \log(z_1b(z_1)\eta(z_1b(z_1))-z_2b(z_2)\eta(z_2b(z_2)))}{\partial z_1\partial z_2}\right.\nonumber\\
  &+&\left.\theta_{\tilde{b},\tilde{H}}(z_1b(z_1))\theta_{\tilde{b},\tilde{H}}(z_2b(z_2))\left(\frac{\partial^2 \log\left(\frac{\underline{m}_{\tilde{H}}(z_1b(z_1))-\underline{m}_{\tilde{H}}(z_2b(z_2))}{z_1b(z_1)-z_2b(z_2) }\right)}{\partial z_1\partial z_2}\right)\right]\,.
\end{eqnarray*}
}
%Note that from the definition of $\tilde{b}(z)$ and proof of Lemma \ref{simplify} follows that $\tilde{b}=1+\delta(z)$.

Since the process of interest $(p-p_1)\left(m_{F^{\bW\bT^{-1}}}(z) - s_n(z)\right)=M_{1,n}+M_{2,n}+M_{3,n}$ forms a tight sequence (see, \cite{baisil2004}, \cite{yao2013} and \cite{zhengbaiyao2015}), the Cauchy integral formula leads to
{
\begin{eqnarray}\label{cauchy}
 && \sum_{i=1}^{p-p_1} f(\lambda_i)-(p-p_1)\int f(x)F_n(dx)  \\ &=&  -\frac{1}{2\pi i}\oint f(z) (p-p_1) (m_{F^{\bW\bT^{-1}}}(z)-s_n(z))dz\nonumber\,,
\end{eqnarray}
}
where $\lambda_i$ is the $i$th eigenvalue of the matrix $\bW\bT^{-1}$ and $f$ is an arbitrary analytic function with 
support containing the interval $[0, r]$, which itself contains the asymptotic spectrum 
of the matrix $\bW\bT^{-1}$. The application of \eqref{cauchy} to our process together with some elementary calculus lead to the result of the theorem.
\end{proof}

\clearpage

\begin{frontmatter}
  \title{Testing for Independence of Large Dimensional Vectors\\
    APPENDIX B: Supplementary material}
\runtitle{Testing for Independence of Large Vectors}
%\subtitle{\vspace{1cm} APPENDIX B: SUPPLEMENTARY MATERIAL}
% \thankstext{T1}{Footnote to the title with the ``thankstext'' command.}

 % \begin{aug}
% \author{\fnms{Taras} \snm{Bodnar}\thanksref{m1}\ead[label=e1]{taras.bodnar@math.su.se}},
% \author{\fnms{Holger} \snm{Dette}\thanksref{t1,m2}\ead[label=e2]{holger.dette@ruhr-uni-bochum.de}}
% \and
% \author{\fnms{Nestor} \snm{Parolya}\thanksref{m3}
% \ead[label=e3]{parolya@statistik.uni-hannover.de}}

%\thankstext{t1}{Some comment}
%\thankstext{t1}{DFG Research Unit 1735, DE 502/26-2}
%\thankstext{t3}{Second supporter of the project}
\runauthor{Bodnar, Dette and Parolya}

\end{frontmatter}

%%%%%%%%%%%%%%%%%%%%%%%%%%%%%%%%%%%%%%%%%%%%%%%%%%%%%%%%%%%%%%%%%%%%%%%%%%%%%%%%%%%%%%%%%%%

%%%%%%%%%%%%%%%%%%%%%%%%%%%%%%%%%%%%%%%%%%%%%%%%%%%%%%%%%%%%%%%%%%%%%%%%%%%%%%%%%%%%%%%%%
\section*{} \label{sec_supp}
\def\theequation{B.\arabic{equation}}
\setcounter{equation}{0}

%%%%%%%%%%%%%

%\section{Proofs of Theorem \ref{th2}, Lemma  \ref{Taras} and Lemma \ref{simplify}}

\begin{proof}[{\bf Proof of Theorem \ref{th2}.}]

The asymptotic properties
of a centred version of \eqref{centralF} have been determined by  \cite{zheng2012}, who showed that
for any functions  $f, g$, which are analytic in an open region of the complex plane containing the interval $[a, b]$,
the random vector
\begin{equation}
 \nonumber
  \Big ((p-p_1)\int_0^\infty f(x) \diff (F_n(x)-F^*_n(x)),  ~ (p-p_1)\int_0^\infty g(x) \diff (F_n(x)-F^*_n(x))   \Big)^\top
\end{equation}
converges weakly to a Gaussian vector $(X_{f}, X_{g})^\top$ with means and covariances given by
\begin{eqnarray*}
  E[X_{f}] &=&\frac{1}{2\pi i}\oint f(z) \diff \log\left( \frac{\frac{1-c}{1-c_1}m_0^2(z)+2m_0(z)+2-c/c_1}{\frac{1-c}{1-c_1}m_0^2(z)+2m_0(z)+1} \right)\nonumber\\
&+& \frac{1}{2\pi i}\oint f(z) \diff \log\left( \frac{1- \frac{c-c_1}{1-c_1}m_0^2(z)}{(1+m_0^2(z))^2} \right)\\
\Cov[X_{f}, X_{g}] &=&-\frac{1}{2 \pi^2} \oint \oint \frac{f(z_1)g(z_2)}{(m_0(z_1)-m_0(z_2))^2}\diff m_0(z_1)\diff m_0(z_2)
\end{eqnarray*}
respectively. Here $m_0(z)=\underline{m}_{\gamma_2}(-\underline{m}(z))$ with $\underline{m}_{\gamma_2}(z)=-\frac{1-\gamma_2}{z}+\gamma_2m_{\gamma_2}(z)$ and $\underline{m}(z)=-\frac{1-\gamma_1}{z}+\gamma_1m(z)$, where $m(z)$ denotes the Stieltjes transform of the function \eqref{centralF} and $m_{\gamma_2}(z)$ is the Stieltjes transformation of the matrix $\bW$ under $H_0$. The integrals are taken over arbitrary positively oriented countur which contains the interval $[a ,b]$.
 Note that this result is only applicable under the null hypothesis  $H_0$, because under $H_1$ the unconditional distribution of  the random matrix
 $\bW$ is no longer a central Wishart distribution. Therefore further investigation are  needed in this situation (see the proof of Theorem  \ref{lsdH1})

 The distributions of the test statistics $T_W$, $T_{LH}$, and $T_{BNP}$ are obtained as special cases using
the functions  $f_W (x)= \log(1+x)$, $f_{LH} (x) = x $ and $ f_{BNP}(x) = \frac{x}{1+x}$, respectively.
 Thus, we need to calculate the asymptotic means, variances and the  terms $\int f(x)dF(x)$ in these cases. The asymptotic means and variances for
 $f_W$ and $f_{LH}$ can be deduced from Examples 4.1 and 4.2 in \cite{zheng2012} and we only need to find the corresponding quantities for $f_{BNP}$.

Let $w$ and $d$ be the positive solutions of the equation
\begin{equation}\label{wd}
|1+hz|^2+(1-\gamma_2)^2=|w+dz|^2
\end{equation}
for any $z\in\mathbbm{C}$ with $|z|=1$ which also satisfy
\begin{equation}  \nonumber
  w^2+d^2=h^2+1+(1-\gamma_2)^2~~\text{and}~~~wd=h\,
\end{equation}
and, consequently, it holds that
\begin{eqnarray*}
(1-\gamma_2)^2&=&(1-d^2)(w^2-1),\\
  (1+h)^2 &=& (w+d)^2-(1-\gamma_2)^2,\\
  (1-h)^2 &=&  (w-d)^2-(1-\gamma_2)^2\,.
\end{eqnarray*}
Further, without loss of generality\footnote{It holds that $|w+dz|^2=|d+wz|^2$ for $|z|=1$.} we assume that $w>d$.

In using that $|1+hz|^2=(1+hz)(h+z)/z$, $|w+dz|^2=(w+dz)(d+wz)$ and due to Corollary 3.2 of \cite{zheng2012}, we get
 {\small
\begin{eqnarray*}
  \E[X_{f_{BNP}}] &=& \lim_{r\downarrow1}\frac{1}{4\pi i} \oint_{|z|=1} \frac{|1+hz|^2/(1-\gamma_2)^2}{|1+hz|^2/(1-\gamma_2)^2+1}\left[ \frac{1}{z-r^{-1}}+\frac{1}{z+r^{-1}}-\frac{2}{z+\gamma_2/h} \right]\diff z\nonumber\\
&=& \lim_{r\downarrow1}\frac{1}{4\pi i} \oint_{|z|=1} \frac{|1+hz|^2}{|w+dz|^2}\left[ \frac{1}{z-r^{-1}}+\frac{1}{z+r^{-1}}-\frac{2}{z+\gamma_2/h} \right]\diff z\nonumber\\
&=& \lim_{r\downarrow1}\frac{1}{4\pi i} \oint_{|z|=1} \frac{(1+hz)(h+z)}{(w+dz)(d+wz)}\left[ \frac{1}{z-r^{-1}}+\frac{1}{z+r^{-1}}-\frac{2}{z+\gamma_2/h} \right]\diff z\nonumber\\
                  &=&  \lim_{r\downarrow1} \frac{1}{2}\left[\left.\frac{(1+hz)(h+z)}{(w+dz)(d+wz)}\right|_{z=r^{-1}}+\left.\frac{(1+hz)(h+z)}{(w+dz)(d+wz)}\right|_{z=-r^{-1}}\right.\nonumber\\
  &-&\left.2\left.\frac{(1+hz)(h+z)}{(w+dz)(d+wz)}\right|_{z=-\frac{\gamma_2}{h}}\right]\nonumber\\
                  &+&  \lim_{r\downarrow1} \frac{1}{2w}\left. \frac{(1+hz)(h+z)}{(w+dz)}\left[ \frac{1}{z-r^{-1}}+\frac{1}{z+r^{-1}}-\frac{2}{z+\gamma_2/h} \right]\right|_{z=-d/w}\nonumber\
\end{eqnarray*}
Thus, after calculating the residuals we obtain
\begin{eqnarray}
 \E[X_{f_{BNP}}]&=&\lim_{r\downarrow1}  \frac{1}{2}\left[ \frac{(1+hr^{-1})(h+r^{-1})}{(w+dr^{-1})(d+wr^{-1})} + \frac{(1-hr^{-1})(h-r^{-1})}{(w-dr^{-1})(d-wr^{-1})}   \right.\nonumber\\
                  &-& 2\left. \frac{(1-\gamma_2)(h-\gamma_2/h)}{(w-\gamma_2d/h)(d-\gamma_2w/h)}+2\frac{(1-dh/w)(h-d/w)}{(w^2-d^2)}\right.\nonumber\\
 &\times& \left.\left( -\frac{d/w}{(d/w)^2-r^{-2}}-\frac{1}{\gamma_2/h-d/w}\right) \right].\nonumber\\
  &&\label{neu}\,
\end{eqnarray}
}

Taking now the limit $r\downarrow1$ in \eqref{neu} we obtain the following result for the mean
{\small
\begin{eqnarray*}
  \E[X_{f_{BNP}}]&=&\frac{1}{2}\left[\frac{(1+h)^2}{(w+d)^2}+\frac{(1-h)^2}{(w-d)^2}\right]-\frac{(1-\gamma_2)(h^2-\gamma_2)}{(w^2-\gamma_2)(d^2-\gamma_2)}\\
  &+& \frac{(1-d^2)(w^2-1)d^2(w^2-\gamma_2)}{(w^2-d^2)^2(d^2-\gamma_2)}\\
                 &=& \left(1-\frac{(1-\gamma_2)^2(w^2+d^2)}{(w^2-d^2)^2}\right) -\left(1+\frac{(1-\gamma_2)^2\gamma_2}{(w^2-\gamma_2)(d^2-\gamma_2)}\right)\\
                 &+& \frac{(1-\gamma_2)^2d^2(w^2-\gamma_2)}{(w^2-d^2)^2(d^2-\gamma_2)}\\
                 &=& -\frac{(1-\gamma_2)^2(w^2+d^2)}{(w^2-d^2)^2}- \frac{(1-\gamma_2)^2\gamma_2}{(w^2-\gamma_2)(d^2-\gamma_2)}\\
                 &+&\left(\frac{(1-\gamma_2)^2\gamma_2(w^2-\gamma_2)}{(w^2-d^2)^2(d^2-\gamma_2)}+\frac{(1-\gamma_2)^2(w^2-\gamma_2)}{(w^2-d^2)^2} \right)\\
&=& \frac{(1-\gamma_2)^2}{(w^2-d^2)^2}\left(-(d^2+\gamma_2)-\gamma_2\left[ \frac{(w^2-d^2)^2 - (w^2-\gamma_2)^2}{(w^2-\gamma_2)(d^2-\gamma_2)}    \right]   \right) \\
% &=& \frac{(1-\gamma_2)^2}{(w^2-d^2)^2}\left(-(d^2+\gamma_2)+\gamma_2\left[ \frac{(d^2-\gamma_2)(2w^2-(d^2 +\gamma_2))}{(w^2-\gamma_2)(d^2-\gamma_2)}    \right]   \right)\\
% &=& \frac{(1-\gamma_2)^2}{(w^2-d^2)^2(w^2-\gamma_2)}(w^2\gamma_2-d^2w^2)\\
&=& -\frac{(1-\gamma_2)^2w^2(d^2-\gamma_2)}{(w^2-d^2)^2}\,.
\end{eqnarray*}
}

Similarly, we have for the variance
{\small
\begin{eqnarray*}
 && \Var[X_{f_{BNP}}]\\
 &=&-\lim_{r\downarrow1} \frac{1}{2\pi^2}\oint_{|z_2|=1}\frac{(1+hz_2)(h+z_2)}{(w+dz_2)(d+wz_2)}\oint_{|z_1|=1} \frac{(1+hz_1)(h+z_1)}{(w+dz_1)(d+wz_1)(z_1-rz_2)^2}\diff z_1\diff z_2 \\
&=& -\lim_{r\downarrow1} \frac{i}{\pi}\oint_{|z_2|=1}\frac{(1+hz_2)(h+z_2)}{(w+dz_2)(d+wz_2)}\left(\left.\frac{(1+hz_1)(h+z_1)}{w(w+dz_1)(z_1-rz_2)^2}\right|_{z_1=-\frac{d}{w}} \right)\diff z_2 \\
 &=& -\frac{i}{\pi}\oint_{|z_2|=1}\frac{(1+hz_2)(h+z_2)}{(w+dz_2)(d+wz_2)}\left(-\frac{dw(1-d^2)(w^2-1)}{(w^2-d^2)(d+wz_2)^2}\right)\diff z_2\,,
\end{eqnarray*}
which simplifies to
\begin{eqnarray*}  
\Var[X_{f_{BNP}}]&=& -\frac{h(1-\gamma_2)^2}{w^3(w^2-d^2)}\left[\left.\frac{\partial^2}{\partial z^2_2}\frac{(1+hz_2)(h+z_2)}{(w+dz_2)}\right|_{z_2=-d/w}\right]\\
%&=&  -\frac{h(1-\gamma_2)^2}{w^3(w^2-d^2)}\left[\left.\frac{\partial}{\partial z_2}\left(\frac{(h+z_2)^2}{(w+dz_2)}-\frac{(1+hz_2)(h+z_2)d}{(w+dz_2)^2}\right)\right|_{z_2=-d/w}\right]\\
% &=&  -2\frac{h(1-\gamma_2)^2}{w^3(w^2-d^2)}\left[\left.\frac{(h+z_2)}{(w+dz_2)}-\frac{(h+z_2)^2d}{(w+dz_2)^2}+\frac{(1+hz_2)(h+z_2)d^2}{(w+dz_2)^3}\right|_{z_2=-d/w}\right]\\
&=& -2\frac{hd(w^2-1)(1-\gamma_2)^2}{w^3(w^2-d^2)^2} \left[ 1-\frac{(w^2-1)d}{w^2-d^2} + \frac{d^2w^2(1-d^2)}{(w^2-d^2)^2}    \right]\\
% &=& -2\frac{hd(w^2-1)(1-\gamma_2)^2}{w^3(w^2-d^2)^4}\Big( (w^2-d^2)(w^2+d)(1-d)+d^2w^2(1-d^2)\Big)  \\
&=&  2\frac{d^2(1-\gamma_2)^4}{w^2(1+d)(w^2-d^2)^4}(w^2(w^2+d)+d^3(w^2-1))\,.
\end{eqnarray*}
}

Due to  Theorem 2.23 in \cite{yaobaizheng2015}, the terms $\int_b^b f(x)\diff F(x)$ can be calculated in the following way
\begin{eqnarray*}
  &&\int_b^b f(x)\diff F(x) \\
  &=& -\frac{h^2(1-\gamma_2)}{4\pi i}\oint_{|z|=1}f\left(\frac{|1+hz|^2}{(1-\gamma_2)^2}\right)\frac{(1-z^2)^2}{z(1+hz)(z+h)(\gamma_2z+h)(\gamma_2+hz)}\diff z\,,
\end{eqnarray*}
where the interval $[a, b]$ is the support of limiting spectral distribution  $F$ of the Fisher matrix $\bW\bT^{-1}$ defined in \eqref{WTlsd}.
The function $f_{LH}$  has  already been considered  in \cite{yaobaizheng2015}, Example 2.25, that is
$
  s_{LH} = \int_b^b x\diff F(x) = \frac{1}{1-\gamma_2}
$.
Next we determine the corresponding terms for $f_W$ and $f_{BNP}$ noting that
\begin{eqnarray*}
  s_{W}
  %&=& -\frac{h^2(1-\gamma_2)}{4\pi i}\oint_{|z|=1}\frac{\log\left(\frac{|1+hz|^2}{(1-\gamma_2)^2}+1\right)(1-z^2)^2}{z(1+hz)(z+h)(\gamma_2z+h)(\gamma_2+hz)}\diff z\\
   &=& -\frac{h^2(1-\gamma_2)}{4\pi i}\oint_{|z|=1}\frac{\log\left((1-\gamma_2)^{-2}|w+dz|^2\right)(1-z^2)^2}{z(1+hz)(z+h)(\gamma_2z+h)(\gamma_2+hz)}\diff z \\
  &=& -\log\left((1-\gamma_2)^2 \right)+ I_1 + I_2\,,
\end{eqnarray*}
where we used again \eqref{wd} and the terms $I_1$ and $I_2$ are defined by
\begin{eqnarray*}
  I_1 &=& -\frac{h^2(1-\gamma_2)}{4\pi i}\oint_{|z|=1}\frac{\log(w+dz)(1-z^2)^2}{z(1+hz)(z+h)(\gamma_2z+h)(\gamma_2+hz)}\diff z\\
I_2 &=&  -\frac{h^2(1-\gamma_2)}{4\pi i}\oint_{|z|=1}\frac{\log(w+d\bar{z})(1-z^2)^2}{z(1+hz)(z+h)(\gamma_2z+h)(\gamma_2+hz)}\diff z\,.
\end{eqnarray*}
A  change of variables yields  $I_1 = I_2$ and we obtain  (see \cite{yaobaizheng2015} for detailed calculation)
{\small
\begin{eqnarray}
 2I_1& =&   -\frac{h^2(1-\gamma_2)}{2\pi i}\oint_{|z|=1}\frac{\log(w+dz)(1-z^2)^2}{z(1+hz)(z+h)(\gamma_2z+h)(\gamma_2+hz)}\diff z
 % \\
%&=& -\frac{1-\gamma_2}{\gamma_2}\log(w) + \frac{\gamma_1+\gamma_2}{\gamma_1\gamma_2}\log(w-d\gamma_2/h)-
  \nonumber\\
  &=&\left\{
\begin{array}{ll}
\frac{1-\gamma_1}{\gamma_1}\log(w-dh), & \gamma_1 \in (0, 1)\\
0, & \gamma_1=1\\
-\frac{1-\gamma_1}{\gamma_1}\log(w-d/h), & \gamma_1 >1\,,
\end{array}
\right.\nonumber
\end{eqnarray}
}
which yields the desired representation of $s_W$.
Similarly, we obtain
{\small
\begin{eqnarray*}
  s_{BNP} &=&  -\frac{h^2(1-\gamma_2)}{4\pi i}\oint_{|z|=1}\frac{(1+hz)(h+z)}{(w+dz)(d+wz)}\frac{(1-z^2)^2}{z(1+hz)(z+h)(\gamma_2z+h)(\gamma_2+hz)}\diff z\\
&=&   -\frac{h^2(1-\gamma_2)}{4\pi i}\oint_{|z|=1}\frac{(1-z^2)^2}{z(w+dz)(d+wz)(\gamma_2z+h)(\gamma_2+hz)}\diff z\\
          &=& -\frac{h^2(1-\gamma_2)}{2}\left(\left. \frac{(1-z^2)^2}{hz(w+dz)(d+wz)(\gamma_2z+h)}\right|_{z=-\frac{\gamma_2}{h}}\right.\\
  &+&\left.\left. \frac{(1-z^2)^2}{zw(w+dz)(\gamma_2z+h)(\gamma_2+hz)}\right|_{z=-\frac{d}{w}} \right.\\
&+& \left. \left.\frac{(1-z^2)^2}{(w+dz)(d+wz)(\gamma_2z+h)(\gamma_2+hz)}\right|_{z=0}  \right)\\
%&=&\frac{(1-\gamma_2)}{2}\left( \frac{(h^2-\gamma_2^2)^2}{\gamma_2(w^2-\gamma_2)(d^2-\gamma_2)(h^2-\gamma_2^2)} + \frac{(w^2-d^2)^2}{(w^2-d^2)(w^2-\gamma_2)(\gamma_2-d^2)}-\frac{h^2}{\gamma_2hwd}\right)\\
%&=&  \frac{(1-\gamma_2)}{2}\left( \frac{h^2-\gamma_2^2}{\gamma_2(w^2-\gamma_2)(d^2-\gamma_2)} + \frac{w^2-d^2}{(w^2-\gamma_2)(\gamma_2-d^2)}-\frac{1}{\gamma_2}\right)\\
%&=& \frac{(1-\gamma_2)}{2\gamma_2(w^2-\gamma_2)(d^2-\gamma_2)}\left( -\gamma_2(w^2-d^2)+h^2-\gamma_2^2-(w^2-\gamma_2)(d^2-\gamma_2)    \right)\\
%&=& \frac{(1-\gamma_2)}{2\gamma_2(w^2-\gamma_2)(d^2-\gamma_2)} (2\gamma_2d^2-2\gamma_2^2)\\
&=& \frac{1-\gamma_2}{w^2-\gamma_2}\,,
\end{eqnarray*}
} which completes the proof of Theorem \ref{th2}.
\end{proof}

%%%%%%%%%% 
\begin{proof}[{\bf Proof of Lemma \ref{Taras}.}]
First, we note that this distribution is independent of $\bSigma_{22\cdot 1}$ since the eigenvalues of $\bW\bT^{-1}$ coincide with the eigenvalues of 
the matrix  $\widetilde{\bW}\widetilde{\bT}^{-1}$, where 
\[\widetilde{\bW}=\bSigma_{22\cdot 1}^{-1/2}\bW\bSigma_{22\cdot 1}^{-1/2} ~~\text{and} ~~\widetilde{\bT}=\bSigma_{22\cdot 1}^{-1/2}\bT\bSigma_{22\cdot 1}^{-1/2}\,,\]
$\bSigma_{22\cdot 1}^{1/2}$ denotes the symmetric   square root of the matrix  $\bSigma_{22\cdot 1}$; $\widetilde{\bT}\sim W_{p-p_1}(n-p_1,\mathbf{I})$; $\widetilde{\bW}|\bS_{11}\sim W_{p-p_1}(p_1,\mathbf{I},\bOmega_1(\bS_{11}))$;  and the random variables $\widetilde{\bT}$ and $(\widetilde{\bW},\bS_{11})$ are independent. Moreover, from the proof of Theorem \ref{lsdH1}, we obtain the  stochastic representations  
\begin{equation*}\label{tbW} \nonumber
  \widetilde{\bW}=\frac{1}{p_1}(\bX+\bSigma^{-1/2}_{22\cdot 1}\bM)(\bX+\bSigma^{-1/2}_{22\cdot 1}\bM)^\top
  ~~\text{and}~~ \widetilde{\bT}=\frac{1}{n-p_1}\bY\bY^\top.
\end{equation*}
where $\bX\sim \mathcal{N}_{p-p_1, p_1}(\mathbf{O},\bI\otimes\bI)$, $\bY\sim \mathcal{N}_{p-p_1, n-p_1}(\mathbf{O},\bI\otimes \bI)$, and $\bX$, $\bY$, $\bS_{11}$ are mutually independent; $\bM=\bSigma_{21}(\bSigma_{11}^{-1/2})^\top\widetilde{\bS}_{11}^{1/2}$ with $\bS_{11}=\bSigma_{11}^{1/2}\widetilde{\bS}_{11}(\bSigma_{11}^{1/2})^\top$ where $\widetilde{\bS}_{11}\sim W_{p_1}(n,\mathbf{I}) $, $\bSigma_{11}^{1/2}$ is the Cholesky root of the positive definite symmetric matrix $\bSigma_{11}$, and $\widetilde{\bS}_{11}^{1/2}$ is the symmetric square root of $\widetilde{\bS}_{11}$. Moreover, it holds that $\widetilde{\bS}_{11}= \bZ\bZ^\top$ with $\bZ\sim \mathcal{N}_{p_1, n}(\mathbf{O},\bI\otimes \bI)$.\\

\noindent \textbf{Case (a): conditional distribution: }
Let 
\[\bSigma^{-1/2}_{22\cdot 1}\bM=\bU_{11}\bD_{11}\bP_{11}^\top\]
be the singular-value decomposition of the matrix $\bSigma^{-1/2}_{22\cdot 1}\bM$ where
\begin{itemize}
\item $\bD_{11}=\bD(\bS_{11})$ is a rectangular diagonal matrix with non-negative real numbers on the diagonal equal to the non-zero eigenvalues of $\bSigma^{-1/2}_{22\cdot 1}\bM \bM^\top \bSigma^{-1/2}_{22\cdot 1}$ or of $\bM^\top\bSigma^{-1}_{22\cdot 1}\bM$;
\item $\bU_{11}=\bU(\bS_{11})$ is the orthogonal matrix consisting of eigenvectors of $\bSigma^{-1/2}_{22\cdot 1}\bM \bM^\top \bSigma^{-1/2}_{22\cdot 1}$;
\item $\bP_{11}=\bP(\bS_{11})$ is the orthogonal matrix consisting of eigenvectors of $\bM^\top\bSigma^{-1}_{22\cdot 1}\bM$.
\end{itemize}

Since $\bU_{11}$ and $\bP_{11}$ are orthogonal
%, that is $\bU_{11}\bU_{11}^\top=\bU_{11}^\top\bU_{11}=\bI$ and $\bP_{11}\bP_{11}^\top=\bP_{11}^\top\bP_{11}=\bI$, 
it follows 
\begin{eqnarray*}
\widetilde{\bW}\widetilde{\bT}^{-1}&=&
\frac{1}{p_1}(\bX+\bSigma^{-1/2}_{22\cdot 1}\bM)(\bX+\bSigma^{-1/2}_{22\cdot 1}\bM)^\top
\Big (\frac{1}{n-p_1}\bY\bY^\top\Big )^{-1}\\
&=&\frac{1}{p_1}(\bU_{11}\bU_{11}^\top\bX\bP_{11}\bP_{11}^\top+\bU_{11}\bD_{11}\bP_{11}^\top)\\
&\times&(\bU_{11}\bU_{11}^\top\bX\bP_{11}\bP_{11}^\top+\bU_{11}\bD_{11}\bP_{11}^\top)^\top
\Big(\frac{1}{n-p_1}\bU_{11}\bU_{11}^\top\bY\bY^\top\bU_{11}\bU_{11}^\top\Big)^{-1}\\
&=&\frac{1}{p_1}\bU_{11}(\bU_{11}^\top\bX\bP_{11}+\bD_{11})\bP_{11}^\top\bP_{11}
    (\bU_{11}^\top\bX\bP_{11}+\bD_{11})^\top\bU_{11}^\top \bU_{11}\\
  &\times&
\Big (\frac{1}{n-p_1}\bU_{11}^\top\bY\bY^\top\bU_{11}\Big )^{-1}\bU_{11}^\top.
\end{eqnarray*}
The  eigenvalues of this matrix  coincide with  the eigenvalues of the matrix
\[\frac{1}{p_1}(\bU_{11}^\top\bX\bP_{11}+\bD_{11})(\bU_{11}^\top\bX\bP_{11}+\bD_{11})^\top
\Big(\frac{1}{n-p_1}\bU_{11}^\top\bY(\bU_{11}^\top\bY)^\top\Big)^{-1},
\]
where $\bU_{11}^\top\bX\bP_{11} \sim \mathcal{N}_{p-p_1, p_1}(\mathbf{O},\bI\otimes\bI)$ and $\bU_{11}^\top\bY\sim \mathcal{N}_{p-p_1, n-p_1}(\mathbf{O},\bI\otimes \bI)$ since the matrix-variate standard normal distribution is invariant to orthogonal transformations. Hence, conditionally on $\bS_{11}$ the distribution of the matrix $\bW\bT^{-1}$ depends only on the eigenvalues of the non-centrality matrix $\bOmega_1(\bS_{11})$.\\

\noindent \textbf{Case (b): unconditional distribution:}
Similarly to  case (a), let
\[\bSigma^{-1/2}_{22\cdot 1}\bSigma_{21}(\bSigma_{11}^{-1/2})^\top=\bU\bD\bP^\top\]
be the singular-value decomposition of $\bSigma^{-1/2}_{22\cdot 1}\bSigma_{21}(\bSigma_{11}^{-1/2})^\top$ where
\begin{itemize}
\item $\bD$ is a rectangular diagonal matrix with non-negative real numbers on the diagonal equal to the non-zero eigenvalues of $\bSigma^{-1/2}_{22\cdot 1}\bSigma_{21}\bSigma_{11}^{-1}\bSigma_{21}^\top \bSigma^{-1/2}_{22\cdot 1}$ or of $\bSigma_{11}^{-1/2}\bSigma_{21}^\top\bSigma^{-1}_{22\cdot 1}\bSigma_{21}(\bSigma_{11}^{-1/2})^\top$;
\item $\bU$ is the orthogonal matrix consisting of eigenvectors of $\bSigma^{-1/2}_{22\cdot 1}\bSigma_{21}\bSigma_{11}^{-1}\bSigma_{21}^\top \bSigma^{-1/2}_{22\cdot 1}$;
\item $\bP$ is the orthogonal matrix consisting of eigenvectors of $\bSigma_{11}^{-1/2}\bSigma_{21}^\top\bSigma^{-1}_{22\cdot 1}\bSigma_{21}(\bSigma_{11}^{-1/2})^\top$.
\end{itemize}

Since $\bU$ and $\bP$ are orthogonal
%, that is $\bU\bU^\top=\bU^\top\bU=\bI$ and $\bP\bP^\top=\bP^\top\bP=\bI$,
 we obtain 
\begin{eqnarray*}
\widetilde{\bW}\widetilde{\bT}^{-1}&=&
\frac{1}{p_1}(\bX+\bSigma^{-1/2}_{22\cdot 1}\bSigma_{21}(\bSigma_{11}^{-1/2})^\top(\bZ \bZ^\top))(\bX+\bSigma^{-1/2}_{22\cdot 1}\bSigma_{21}(\bSigma_{11}^{-1/2})^\top(\bZ \bZ^\top))^\top
  \\
  &\times&\Big(\frac{1}{n-p_1}\bY\bY^\top\Big)^{-1}\\
&=&\frac{1}{p_1}(\bU\bU^\top\bX\bP\bP^\top+\bU\bD\bP^\top (\bZ \bZ^\top)\bP\bP^\top)\\
&\times&(\bU\bU^\top\bX\bP\bP^\top+\bU\bD\bP^\top (\bZ \bZ^\top)\bP\bP^\top)^\top
\Big(\frac{1}{n-p_1}\bU\bU^\top\bY\bY^\top\bU\bU^\top\Big)^{-1}\\
&=&\frac{1}{p_1}\bU(\bU^\top\bX\bP+\bD(\bP^\top\bZ\bZ^\top\bP))\bP^\top\bP
(\bU^\top\bX\bP+\bD(\bP^\top\bZ\bZ^\top\bP))^\top\bU^\top \bU\\
&\times&\Big(\frac{1}{n-p_1}\bU^\top\bY\bY^\top\bU\Big)^{-1}\bU^\top . 
\end{eqnarray*}
The  eigenvalues of this matrix  coincide with the eigenvalues of 
\[\frac{1}{p_1}(\bU^\top\bX\bP+\bD(\bP^\top\bZ(\bP^\top\bZ)^\top))
(\bU^\top\bX\bP+\bD(\bP^\top\bZ(\bP^\top\bZ)^\top))^\top
\Big(\frac{1}{n-p_1}\bU^\top\bY(\bU^\top\bY)^\top\Big)^{-1},
\]
where $\bU^\top\bX\bP \sim \mathcal{N}_{p-p_1, p_1}(\mathbf{O},\bI\otimes\bI)$, $\bU^\top\bY\sim \mathcal{N}_{p-p_1, n-p_1}(\mathbf{O},\bI\otimes \bI)$, and 
$\bP^\top\bZ\sim \mathcal{N}_{p_1, n}(\mathbf{O},\bI\otimes \bI)$,
since the matrix-variate standard normal distribution is invariant to orthogonal transformations. Hence, the unconditional distribution of the matrix $\bW\bT^{-1}$ depends only on the eigenvalues of the matrix $\widetilde{\bR}=\bSigma_{11}^{-1}\bSigma_{12}\bSigma_{22\cdot 1}^{-1}\bSigma_{21}$.

\end{proof}

\begin{proof}[{\bf Proof of Lemma \ref{simplify}.}]
Let
\begin{eqnarray*}
T_n(z) &=& \left(\frac{1}{1+\delta_n(z)}\bSigma^{-1/2}_{22\cdot1}\bM\bM^\top\bSigma^{-1/2}_{22\cdot1} -z(1+\tilde{\delta}_n(z))\bI_{p-p_1}\right)^{-1}   \\
\tilde{T}_n(z) &=& \left( \frac{1}{1+\tilde{\delta}_n(z)}\bM^\top\bSigma^{-1}_{22\cdot1}\bM -z(1+\delta_n(z))\bI_{p_1}\right)^{-1}\,,
\end{eqnarray*}
where $\delta_n(z)$ and $\tilde{\delta}_n(z)$ are the unique solutions of the following system of equations
\begin{eqnarray}  \nonumber 
  \delta_n(z) &=&\frac{1}{p_1}\text{tr}\Bigl(T_n(z) \Bigr),~~~\tilde{\delta}_n(z) = \frac{1}{p_1}\text{tr}\left(\tilde{T}_n(z) \right)\,
\end{eqnarray}
in the class of Stieltjes transforms of non-negative measures\footnote{In fact, $\delta_n$ is the Stieltjes transform of a measure with total mass equal to $\frac{p-p_1}{p_1}$ while $\tilde{\delta}_n$ is the Stieltjes transform of a measure with total mass equal to $1$ (see,  \cite{hachem2012})} with support in $\mathbbm{R}^+$.

The functions $T_n(z) $ and $\tilde{T}_n(z) $ are the deterministic approximations of the resolvents 
$$Q_n(z) = \left((\bX+\bSigma^{-1/2}_{22\cdot1}\bM)(\bX+\bSigma^{-1/2}_{22\cdot1}\bM)^\top-z\bI_{p-p_1} \right)^{-1},$$
$$\tilde{Q}_n(z) =\left((\bX+\bSigma^{-1/2}_{22\cdot1}\bM)^\top(\bX+\bSigma^{-1/2}_{22\cdot1}\bM)-z\bI_{p_1} \right)^{-1},$$
respectively, in the sense that
\begin{eqnarray}
  \frac{1}{p-p_1}\text{tr}(Q_n(z)-T_n(z)) \overset{a.s.}{\longrightarrow} 0~~\text{and}~~\frac{1}{p_1}\text{tr}(\tilde{Q}_n(z)-\tilde{T}_n(z)) \overset{a.s.}{\longrightarrow} 0~~\text{as}~n\to\infty\,. \nonumber 
\end{eqnarray}
%The matrix $\bX$ is a corresponding $(p-p_1)\times p_1$-dimensional  random matrix independent of $\bT$ with standard matrix-variate normal distribution.

First, we find the connection between $\delta_n(z)$ and the Stieltjes transform $m_{H}(z)$, where $H$ is the limiting spectral distribution of $\bW$. For that reason, we consider the asymptotic values of $\delta_n(z)$ and $\tilde{\delta}_n(z)$ given by
\begin{eqnarray*}
  \delta_n(z) &=& \frac{p-p_1}{p_1}(1+\delta_n(z))\frac{1}{p-p_1}\text{tr}\left(\bSigma^{-1/2}_{22\cdot1}\bM\bM^\top\bSigma^{-1/2}_{22\cdot1} -z(1+\tilde{\delta}_n(z))(1+\delta_n(z))\bI_{p-p_1}\right)^{-1}\nonumber\\
&=& (1+\delta_n(z))\frac{p-p_1}{p_1}\int\frac{d\tilde{H}_n(t)}{t-z{\eta_n(z)}} \longrightarrow \delta(z)= (1+\delta(z))\gamma_1m_{\tilde{H}}(z\eta(z))\label{delta_inf}\\
\tilde{\delta}_n(z) &=& \frac{1}{p_1}(1+\tilde{\delta}_n(z))\text{tr}\left(\bM^\top\bSigma^{-1}_{22\cdot1}\bM -z(1+\tilde{\delta}_n(z))(1+\delta_n(z))\bI_{p_1}\right)^{-1}\nonumber\\
&=& (1+\tilde{\delta}_n(z))\frac{1}{p_1}\int\frac{d\tilde{\underline{H}}_n(t)}{t-z\eta(z)} \longrightarrow \tilde{\delta}(z)= (1+\tilde{\delta}(z))m_{\tilde{\underline{H}}}(z{\eta_n(z)})\nonumber\\
&=&(1+\tilde{\delta}(z))\left(-\frac{1-\gamma_1}{z\eta(z)}+\gamma_1m_{\tilde{H}}(z\eta(z)) \right)\label{tilde_delta_inf}
\end{eqnarray*}
with 
$$\eta_n(z) = (1+\delta_n(z))(1+\tilde{\delta}_n(z)) ~~ \text{and} ~~\eta(z) = (1+\delta(z))(1+\tilde{\delta}(z)).$$
Equivalently, we have
{
  \begin{eqnarray}
    && \frac{\delta(z)}{1+\delta(z)} = \gamma_1m_{\tilde{H}}\big(z\eta(z) \big)~~\text{and}~~\frac{\tilde{\delta}(z)}{1+\tilde{\delta}(z)}=m_{\tilde{\underline{H}}}(z\eta(z)) =-\frac{1-\gamma_1}{z\eta(z)} + \gamma_1m_{\tilde{H}}\big(z\eta(z) \big) \nonumber\\ && \label{inf_delta}
\end{eqnarray}
}
which leads to
\begin{eqnarray}\label{link_delta}
  \tilde{\delta}(z) = -\frac{1-\gamma_1}{z}+\delta(z)\,.
\end{eqnarray}
We claim that in fact we have
\begin{eqnarray}\label{claim}
\delta(z)=\gamma_1m_{H}(z)\,
\end{eqnarray}
and, consequently, $\tilde{\delta}(z)=-\frac{1-\gamma_1}{z}+\gamma_1m_{H}(z)$.
 In order to prove \eqref{claim}, we plug $\delta(z)=\gamma_1m_{H}(z)$ into \eqref{inf_delta} and use \eqref{link_delta}. It leads to
 \begin{eqnarray*}
   \frac{m_{H}(z)}{1+\gamma_1m_{H}(z)}&=& m_{\tilde{H}}\left(z[1+\gamma_1m_{H}(z)]\left[1+\gamma_1m_{H}(z)-\frac{1-\gamma_1}{z}\right]   \right)\\
&=& m_{\tilde{H}}\Big([1+\gamma_1m_{H}(z)][z(1+\gamma_1m_{H}(z))-(1-\gamma_1)] \Big)\,,
 \end{eqnarray*}
which is exactly \eqref{H}. From the uniqueness of the solution the claim \eqref{claim} follows.
Thus, in light of Theorem {\color{blue} 2} we get as $n\to\infty$
\begin{eqnarray*}
&&  \delta_n(z) \longrightarrow \gamma_1m_{H}(z),\label{d_1}\\
&&  \tilde{\delta}_n(z) \longrightarrow -\frac{1-\gamma_1}{z}+\gamma_1m_{H}(z)=m_{\underline{H}}(z)\label{d_2}\,.
\end{eqnarray*}

For all $z\in \mathbbm{C}^+$, we define
\begin{eqnarray*}
  \Delta_n(z) &=& \Big 
  (1 - \frac{\frac{1}{p_1}\text{tr}[T^2_n(z)\bSigma_{22\cdot1}^{-1/2}\bM\bM^\top\bSigma_{22\cdot1}^{-1/2}]}{(1+\delta_n(z))^2}    \Big )^2 - z^2\xi_n(z)\tilde{\xi}_n(z)\\
  \Psi_n(z) &=& \Big (1-z\xi_n(z) - \frac{\frac{1}{p_1}\text{tr}\left(T_n^2(z)\bSigma_{22\cdot1}^{-1/2}\bM\bM^\top\bSigma_{22\cdot1}^{-1/2}\right)}{(1+\delta_n(z))^2} \Big )^{-1}\\
\omega_n(z) &=& \frac{1}{p_1}\sum_{j=1}^{p_1}z^2\tilde{t}^2_{jj}\label{omega}\\
\zeta_n(z)&=&\frac{1}{p-p_1}\underset{k\neq l}{\sum\limits_{k=1}^{p-p_1}\sum\limits_{l=1}^{p-p_1}}(\mathbf{m}_k^\top T_n(z) \mathbf{m}_l)^2\label{zeta}\,
\end{eqnarray*}
with $\tilde{t}_{jj}$ being the diagonal elements of the matrix $\tilde{T}_n$ and $\mathbf{m}_k$ - the $k$th column of matrix $\bSigma_{22\cdot1}^{-1/2}\bM$, while
{\small
\begin{eqnarray}
  \xi_n(z_1,z_2) &=&  \frac{1}{p_1}\text{tr}(T_n(z_1)T_n(z_2)),~~\tilde{\xi}_n(z_1,z_2) =   \frac{1}{p_1}\text{tr}(\tilde{T}_n(z_1)\tilde{T}_n(z_2)) \,
\end{eqnarray}}
and, obviously, $\xi_n(z)\equiv\xi_n(z,z)$ and $\tilde{\xi}_n(z)\equiv\tilde{\xi}_n(z,z)$. 

Next, we simplify the above expressions. In using \eqref{inf_delta}, we get
{\small
\begin{eqnarray*}
  \xi_n(z_1,z_2) &=&  \frac{p-p_1}{p_1}\frac{(1+\delta_n(z_1))(1+\delta_n(z_2))}{p-p_1}\\
  &\times&\text{tr}\left(\left[\bSigma^{-1/2}_{22\cdot1}\bM\bM^\top\bSigma^{-1/2}_{22\cdot1} - z_1\eta_n(z_1)\bI_{p-p_1}\right]^{-1}\left[\bSigma^{-1/2}_{22\cdot1}\bM\bM^\top\bSigma^{-1/2}_{22\cdot1} - z_2\eta_n(z_2)\bI_{p-p_1}\right]^{-1}     \right)\nonumber\\
&\longrightarrow&\xi(z_1,z_2) = \gamma_1(1+\delta(z_1))(1+\delta(z_2))\int\frac{d\tilde{H}(t)}{(t-z_1\eta(z_1))(t-z_2\eta(z_2))}\nonumber\\
&=&\gamma_1(1+\delta(z_1))(1+\delta(z_2)) \frac{m_{\tilde{H}}(z_1\eta(z_1))-m_{\tilde{H}}(z_2\eta(z_2))}{z_1\eta(z_1)-z_2\eta(z_2)}\nonumber\\
&=& \frac{\delta(z_1)-\delta(z_2)}{z_1\eta(z_1)-z_2\eta(z_2)}\,.
\end{eqnarray*}
}
In the case of $z_1=z_2=z$, we obtain
\begin{eqnarray*}
\xi_n(z) &\longrightarrow& \xi(z) =  \gamma_1(1+\delta(z))^2\int\frac{d\tilde{H}(t)}{(t-z\eta(z))^2}\nonumber\\
&=& \gamma_1\frac{m'_{\tilde{H}}(z\eta(z))}{z\eta'(z)+\eta(z)}(1+\delta(z))^2 = \frac{\delta'(z)}{(z\eta(z))'}\,.
\end{eqnarray*}
Similarly, using \eqref{inf_delta} we get for $\tilde{\xi}_n(z_1,z_2)$, i.e.,
{\small
\begin{eqnarray*}
 \tilde{\xi}_n(z_1,z_2) &\longrightarrow&\tilde{\xi}(z_1,z_2) =\frac{\tilde{\delta}(z_1)-\tilde{\delta}(z_2)}{z_2\eta(z_2)-z_1\eta(z_1)} = (1-\gamma_1)\frac{z_1-z_2}{(z_1\eta(z_1)-z_2\eta(z_2))z_1z_2}+\frac{\delta(z_1)-\delta(z_2)}{z_1\eta(z_1)-z_2\eta(z_2)}\nonumber\\
&=& \frac{(1-\gamma_1)}{z_1z_2}\frac{z_1-z_2}{z_1\eta(z_1)-z_2\eta(z_2)}+\xi(z_1,z_2) \\
\tilde{\xi}_n(z)& \longrightarrow& \tilde{\xi}(z) =  \frac{\left(\frac{1-\gamma_1}{z^2}+ \delta'(z)\right)}{(z\eta(z))'}
%= (z\eta(z))^\top \left(\frac{1-\gamma_1}{z}+\delta(z) \right)^\top
= \frac{(1-\gamma_1)}{z^2(z\eta(z))'} + \xi(z)\,.
\end{eqnarray*}
} 
In using these results as well as
{\small
\begin{eqnarray}
 m_{\tilde{H}}'(z\eta(z)) &=& \frac{\partial}{\partial z}\int\frac{d\tilde{H}(t)}{(t-z\eta(z))} = \int\frac{d\tilde{H}(t)}{(t-z\eta(z))^2}(z\eta(z))', \label{first_derST}\\
\gamma_1m_{\tilde{H}}'(z\eta(z)) &=& \frac{\delta'(z)}{(1+\delta)^2}\label{first_derST2}\,
\end{eqnarray}
}
and applying
{\footnotesize
\begin{eqnarray*}
  \gamma_1\int\frac{d\tilde{H}(t)}{(t-z\eta(z))^2}&=& \frac{\delta'(z)}{(1+\delta(z))^2}\frac{1}{(z\eta(z))'}=\frac{\xi(z)}{(1+\delta(z))^2},\\
  \gamma_1\int\frac{td\tilde{H}(t)}{(t-z\eta(z))^2} &=&\gamma_1\int\frac{\tilde{H}(t)}{(t-z\eta(z))} + \gamma_1z\eta(z)\int\frac{\tilde{H}(t)}{(t-z\eta(z))^2}= \frac{\delta(z)}{1+\delta(z)} + \frac{\xi(z)}{(1+\delta(z))^2}z\eta(z)\,,
\end{eqnarray*}
}
 we get
{
\begin{eqnarray}\label{Delta_inf}
 \Delta_n(z)& \longrightarrow& \Delta(z)= \left(1- \gamma_1\int\frac{t\tilde{H}(t)}{(t-z\eta(z))^2} \right)^2-z^2\xi(z)\tilde{\xi}(z)\nonumber\\
%&=& \left(1- \gamma_1\int\frac{\tilde{H}(t)}{(t-z\eta(z))} - \gamma_1z\eta(z)\int\frac{\tilde{H}(t)}{(t-z\eta(z))^2} \right)^2--z^2\xi(z)\tilde{\xi}(z)\nonumber\\
%&=& \left(\frac{1}{1+\delta(z)}- \frac{z\eta(z)}{(z\eta(z))^\top}\frac{\delta^\top(z)}{(1+\delta(z))^2} \right)^2-z^2\left(\frac{\delta^\top(z)}{(z\eta(z))^\top}\right)^2+\frac{1-\gamma_1}{\left[(z\eta(z))^\top\right]^2}\delta^\top(z)
&=& \left(\frac{1}{1+\delta(z)} - \frac{z\eta(z)}{(1+\delta(z))^2}\xi(z)\right)^2-z^2\xi(z)\tilde{\xi}(z)\nonumber\\
&=& \left(\frac{1}{1+\delta(z)} - z\xi(z)+\frac{1-\gamma_1}{1+\delta(z)}\xi(z)\right)^2-z^2\xi^2(z)-\frac{1-\gamma_1}{(z\eta(z)')}\xi(z)\,.
\end{eqnarray}
}
Moreover, the term $(z\eta(z))'$ can be rewritten further as follows
\begin{eqnarray*}
(z\eta(z))' &=& (z(1+\delta(z))(1+\tilde{\delta}(z)))' =\left(z(1+\delta(z))\left(1+\delta(z)-\frac{1-\gamma_1}{z}\right)\right)'\nonumber\\
&=& (z(1+\delta(z))^2-(1-\gamma_1)(1+\delta(z)))' = (1+\delta(z))^2 + 2(1+\delta(z))\delta'(z)z-(1-\gamma_1)\delta'(z)\nonumber\\
&=& (1+\delta(z))^2 + 2(z\eta(z))'(1+\delta(z))\xi(z)z-(1-\gamma_1)(z\eta(z))'\xi(z)\,,
\end{eqnarray*}
which yields to
\begin{eqnarray}\label{zeta}  \nonumber 
  \frac{1}{(z\eta(z))'} =  \frac{1}{1+\delta(z)}\left(\frac{1}{1+\delta(z)}-2\xi(z)z+\frac{1-\gamma_1}{1+\delta(z)}\xi(z)\right)
\end{eqnarray}
Similarly, 
\begin{eqnarray}
  \Psi^{-1}_n(z)\longrightarrow \Psi^{-1}(z)& =& \frac{1}{1+\delta(z)} - \frac{z\eta(z)}{(1+\delta(z))^2}\xi(z) - z\xi(z)\nonumber\\
&=&\frac{1}{1+\delta(z)}-2\xi(z)z+\frac{1-\gamma_1}{1+\delta(z)}\xi(z)\,,
\end{eqnarray}
which is exactly equal to $(1+\delta(z))/(z\eta(z))'$.
Now, \eqref{zeta} and \eqref{Delta_inf} lead to
{\small
\begin{eqnarray}
  \Delta(z) &=& \left(\frac{1}{1+\delta(z)}+\frac{1-\gamma_1}{1+\delta(z)}\xi(z)\right)\left(\frac{1}{1+\delta(z)}-2\xi(z)z+\frac{1-\gamma_1}{1+\delta(z)}\xi(z) \right) -\frac{1-\gamma_1}{(z\eta(z))'}\xi(z)\nonumber\\
&=& \left(\frac{1}{1+\delta(z)}+\frac{1-\gamma_1}{1+\delta(z)}\xi(z)\right)\Psi^{-1}(z)-\frac{1-\gamma_1}{1+\delta(z)}\xi(z)\Psi^{-1}(z)\nonumber\\
            &=& \frac{1}{1+\delta(z)}\Psi^{-1}(z).\nonumber\\
  && \label{Delta_simple}\,
\end{eqnarray}
}
%%%%%%%%%%%%%%%%%%%%%%%%%%%%%%%%%%%%%%%%%%%%%%%%%%%%%%%%%%%%%%%%%%%%%%%%%%%%%%%%%%%%%%%%%
From Lemma {\color{blue} 1} we get that the matrices $T_n(z)$ and $\tilde{T}_n(z)$ could be chosen without loss of generality as diagonal matrices, which implies
\begin{eqnarray*}
  \omega_n(z)&=&z^2\frac{1}{p_1}\sum\limits_{j=1}^{p_1}\tilde{t}^2_{jj}{=} z^2\text{tr}(\tilde{T}^2_n(z))\rightarrow z^2\tilde{\delta}^2(z),\\
  \zeta_n(z) &=& \frac{1}{p-p_1}\underset{k\neq l}{\sum\limits_{k=1}^{p-p_1}\sum\limits_{l=1}^{p-p_1}}(\mathbf{m}_k^\top T_n(z) \mathbf{m}_l)^2=0\,.
\end{eqnarray*}

Now, Theorems 2.2.1 and 2.2.2 by  \cite{yao2013} reveal that $M_{2,n}(z)=(p-p_1)(m_{F^{\bW}}(z)-m_{H_n}(z))$ converges to a Gaussian process $M_2(z)$ with mean function and covariance function given by
{\scriptsize
 \begin{eqnarray*}
   {\rm E}(M_2(z))& =& \frac{\Psi_n(z)}{\Delta_n(z)}\left(z^2\tilde{\xi}_n(z)\frac{1}{p_1}\text{tr}\big(T^3_n(z)\big)+\zeta_n(z)\frac{1}{p_1}\text{tr}\big(T^3_n(z)\big) \right.\nonumber\\
   &+&\left. 2\frac{\frac{1}{p_1}\text{tr}\big(\bSigma_{22\cdot1}^{-1/2}\bM\bM^\top\bSigma_{22\cdot1}^{-1/2} T^3_n(z)\big)}{(1+\delta_n(z))^2}\left(1   -  \frac{\frac{1}{p_1}\text{tr}\left(\bSigma_{22\cdot1}^{-1/2}\bM\bM^\top\bSigma_{22\cdot1}^{-1/2} T_n^2(z) \right)}{(1+\delta_n(z))^2} \right) \right.\nonumber\\
&+& \left. \frac{\omega_n(z)}{(1+\delta_n(z))^2}\left( \frac{1}{p_1}\text{tr}\big(\bSigma_{22\cdot1}^{-1/2}\bM\bM^\top\bSigma_{22\cdot1}^{-1/2} T^3_n(z)\big) \xi_n(z) - \frac{1}{p_1^2}\text{tr}\big(\bSigma_{22\cdot1}^{-1/2}\bM\bM^\top\bSigma_{22\cdot1}^{-1/2} T^2_n(z)\big)\text{tr}\big(T^3_n(z)\big) \right)   \right)  \nonumber\\
 {\rm Cov}(M_2(z_1), M_2(z_2))& =& 2\frac{(z_1\eta_n(z_1))'(z_2\eta_n(z_2))'}{(z_1\eta_n(z_1)-z_2\eta_n(z_2))^2} \,.
 \end{eqnarray*}
}
Since $\eta_n(z)\to \eta(z)$, we get that
\begin{eqnarray*}
  {\rm Cov}(M_2(z_1), M_2(z_2))& \longrightarrow& 2\frac{\eta'(z_1)\eta'(z_2)}{(\eta(z_1)-\eta(z_2))^2}=2\frac{\partial \log(\eta(z_1)-\eta(z_2))}{\partial z_1\partial z_2}\,.
\end{eqnarray*}
Furthermore, 
{\scriptsize
\begin{eqnarray*}
    {\rm E}(M_2(z))  &\longrightarrow& B(z) = \frac{\Psi(z)}{\Delta(z)}\left(z^2\tilde{\xi}(z)\gamma_1\int\frac{d\tilde{H}(t)}{(t-z\eta(z))^3}(1+\delta(z))^3 + 2\gamma_1\frac{\int\frac{td\tilde{H}(t)}{(t-z\eta(z))^3}(1+\delta(z))^3}{(1+\delta(z))^2}\left[\Psi^{-1}(z)+z\xi(z)\right] \right.\nonumber\\
&+& \left. \frac{z^2\tilde{\delta}^2(z)}{(1+\delta(z))^2}\left( \gamma_1\int\frac{td\tilde{H}(t)}{(t-z\eta(z))^3}(1+\delta(z))^3 \xi_n(z) - \gamma_1^2\int\frac{td\tilde{H}(t)}{(t-z\eta(z))^2}\int\frac{d\tilde{H}(t)}{(t-z\eta(z))^3}(1+\delta(z))^5 \right)   \right)\nonumber\,,
\end{eqnarray*}
} 
where from \eqref{Delta_simple} it follows that
\begin{eqnarray*}
  \frac{\Psi(z)}{\Delta(z)} = (1+\delta(z))\Psi^2(z)\,.
\end{eqnarray*}
Because of \eqref{first_derST}, \eqref{first_derST2} and
\begin{eqnarray*}
  \frac{\partial^2}{\partial z^2}\int\frac{td\tilde{H}(t)}{(t-z\eta(z))} &=& 2    \int\frac{d\tilde{H}(t)}{(t-z\eta(z))^3} \left[(z\eta(z))^\top\right]^2 +  \int\frac{d\tilde{H}(t)}{(t-z\eta(z))^2} (z\eta(z))^{''}\nonumber\\
&=& 2  \int\frac{d\tilde{H}(t)}{(t-z\eta(z))^3} \left[(z\eta(z))'\right]^2 + m'_{\tilde{H}}(z\eta(z))\frac{(z\eta(z))^{''}}{(z\eta(z))^{'}}\nonumber\\
&=&2  \int\frac{d\tilde{H}(t)}{(t-z\eta(z))^3} \left[(z\eta(z))'\right]^2 + \gamma^{-1}\frac{\delta'(z)}{(1+\delta(z))^2}\frac{(z\eta(z))^{''}}{(z\eta(z))^{'}},\\
\gamma_1m^{''}_{\tilde{H}}(z\eta(z))&=& \frac{\delta^{''}(z)}{(1+\delta(z)^2} - 2 \frac{\delta^{'\; 2}(z)}{(1+\delta(z))^3},
\end{eqnarray*}
we obtain
\begin{eqnarray*}
  \gamma_1 \int\frac{d\tilde{H}(t)}{(t-z\eta(z))^3} &=& \frac{ \frac{\delta^{''}(z)}{(1+\delta(z))^2} - 2 \frac{\delta^{'\; 2}(z)}{(1+\delta(z))^3}-\frac{\delta'(z)}{(1+\delta(z))^2}\frac{(z\eta(z))^{''}}{(z\eta(z))^{'}}}{2(z\eta(z))^{'\; 2}} \,.
\end{eqnarray*}
On the other hand, it holds that
\begin{eqnarray*}
  \xi'(z) = \frac{\delta^{''}(z)}{(z\eta(z))'} - \delta'(z)\frac{(z\eta(z))^{''}}{(z\eta(z))^{'\; 2}} =  \frac{\delta^{''}(z)}{(z\eta(z))'} -\xi(z)\frac{(z\eta(z))^{''}}{(z\eta(z))^{'}}\,.
\end{eqnarray*}
Thus, we have
\begin{eqnarray*}
  \gamma_1 \int\frac{d\tilde{H}(t)}{(t-z\eta(z))^3}& =& \frac{\xi'(z)}{2(1+\delta(z))^2(z\eta(z))'}-\frac{\xi^2(z)}{(1+\delta(z))^3}\nonumber\\
&=& \frac{1}{(1+\delta(z))^3}\left(\frac{\xi'(z)\Psi^{-1}(z)}{2} - \xi^2(z) \right)\,.
\end{eqnarray*}
and
{\small
\begin{eqnarray*}
  \gamma_1 \int\frac{td\tilde{H}(t)}{(t-z\eta(z))^3} &=&\gamma_1 \int\frac{d\tilde{H}(t)}{(t-z\eta(z))^2}+\gamma_1z\eta(z) \int\frac{d\tilde{H}(t)}{(t-z\eta(z))^3}\nonumber\\
%&=& \frac{\delta^\top(z)}{(1+\delta(z))^2}\frac{1}{(z\eta(z)^\top}+(z\eta(z))\frac{ \frac{\delta^{''}(z)}{(1+\delta(z))^2} - 2 \frac{\delta^{\top\; 2}(z)}{(1+\delta(z))^3}-\frac{\delta^\top(z)}{(1+\delta(z))^2}\frac{(z\eta(z))^{''}}{(z\eta(z))^{'}}}{2(z\eta(z))^{\top\; 2}}
%&=& \frac{1}{(1+\delta(z))^2}\left(\xi(z) + z\eta(z)\left(\frac{\xi^\top(z)}{2(z\eta(z))^\top}-\frac{\xi^2(z)}{(1+\delta(z))} \right)\right)\nonumber\\
&=&\frac{1}{(1+\delta(z))^2}\left(\xi(z) + \left(z-\frac{(1-\gamma_1)}{1+\delta(z)}\right)\left(\frac{\xi'(z)\Psi^{-1}(z)}{2}-\xi^2(z)\right)\right)\,.
\end{eqnarray*}
}
As a result, it holds that
{\footnotesize
\begin{eqnarray*}
  B(z)&=&(1+\delta(z))\Psi^2(z)\left( \left(z^2\xi(z) + \frac{1-\gamma_1}{1+\delta(z)}\Psi^{-1}(z) \right)\left(\frac{\xi^\top(z)\Psi^{-1}(z) }{2} - \xi^2(z)\right)\right.\nonumber\\
&+&\left. \frac{2}{1+\delta(z)}\left( \xi(z)+\left(z-\frac{1-\gamma_1}{1+\delta(z)} \right)\left(\frac{\xi'(z)\Psi^{-1}(z) }{2} - \xi^2(z)\right)    \right)\left(\Psi^{-1}(z)+z\xi(z) \right)\right.\nonumber\\
&+& \left.\frac{\omega(z)}{(1+\delta(z))}\left(\xi(z)\left( \xi(z)+\left(z-\frac{1-\gamma_1}{1+\delta(z)} \right)\left(\frac{\xi'(z)\Psi^{-1}(z) }{2} - \xi^2(z)\right)\right)\right.\right.\nonumber\\
&-&\left.\left. \left(\delta(z)+\xi(z)(z(1+\delta(z))-(1-\gamma_1))\right) \left( \frac{\xi'(z)\Psi^{-1}(z) }{2} - \xi^2(z)  \right)\right)\right)\nonumber\\
&=& (1+\delta(z))\Psi^2(z) \left[\widetilde{\omega}(z)N(z) + \frac{2}{1+\delta(z)}\left(\xi(z)+\left(z-\frac{1-\gamma_1}{1+\delta(z)}\right)N(z)\right)\left[\Psi^{-1}(z) + z\xi(z)\right]\right.\nonumber\\
&+& \left.\frac{\omega(z)}{1+\delta(z)} \left( \xi^2(z) - \delta(z)N(z)\left(z-\frac{1-\gamma_1}{1+\delta(z)}+1\right)\right) \right]\nonumber\\
&=& (1+\delta(z))\Psi^2(z) \left[ \widetilde{\omega}(z)N(z) + \frac{2}{1+\delta(z)}\Bigl(-\widetilde{\omega}(z)N(z)+\xi(z)(\Psi^{-1}(z)+z\xi(z))\right.\nonumber\\
&+&\left.z\underbrace{\left(2z\xi(z)-\frac{1-\gamma_1}{1+\delta(z)}\xi(z)-\frac{1}{1+\delta(z)}\right)}_{-\Psi^{-1}(z)}N(z) + \frac{1}{1+\delta(z)}N(z) + z\Psi^{-1}(z)N(z)   \Bigr)\right.\nonumber\\
 &+&   \left.   \frac{\omega(z)}{1+\delta(z)} \left( \xi^2(z) - \delta(z)N(z)\left(z-\frac{1-\gamma_1}{1+\delta(z)}+1\right)\right)      \right]\nonumber\\
&=& (1+\delta(z))\Psi^2(z) \left[ \widetilde{\omega}(z)N(z)\frac{\delta(z)-1}{1+\delta(z)} \right.\nonumber\\
&+&\left.\frac{1}{1+\delta(z)}\left(\frac{1}{(1+\delta(z))}N(z)+\xi(z)\left(\Psi^{-1}(z)+z\xi(z)\right)\right)\right.\nonumber\\
&+&  \left.   \frac{\omega(z)}{1+\delta(z)} \left( \xi^2(z) - \delta(z)N(z)\left(z-\frac{1-\gamma_1}{1+\delta(z)}+1\right)\right)      \right]
\end{eqnarray*}
with
\begin{eqnarray*}
N(z) = \frac{\xi'(z)\Psi^{-1}(z) }{2} - \xi^2(z)~~~\text{and}~~~\widetilde{\omega}(z)=z^2\xi(z) + \frac{1-\gamma_1}{1+\delta(z)}\Psi^{-1}(z)\,.
\end{eqnarray*}
}
\end{proof}

%\section{Additional simulations}

\begin{figure}[p]
  \centering
  \includegraphics[scale=0.28]{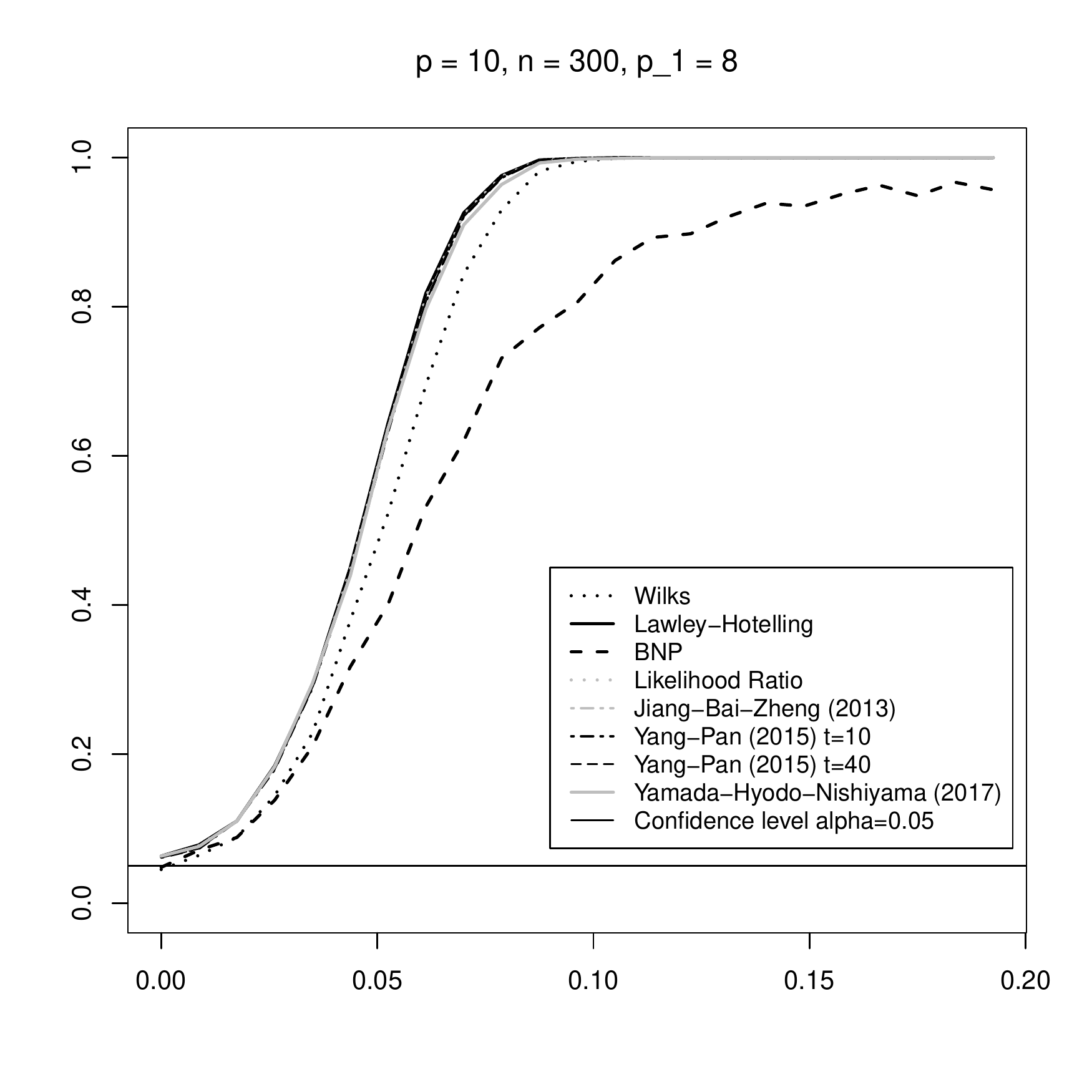} ~~
  \includegraphics[scale=0.28]{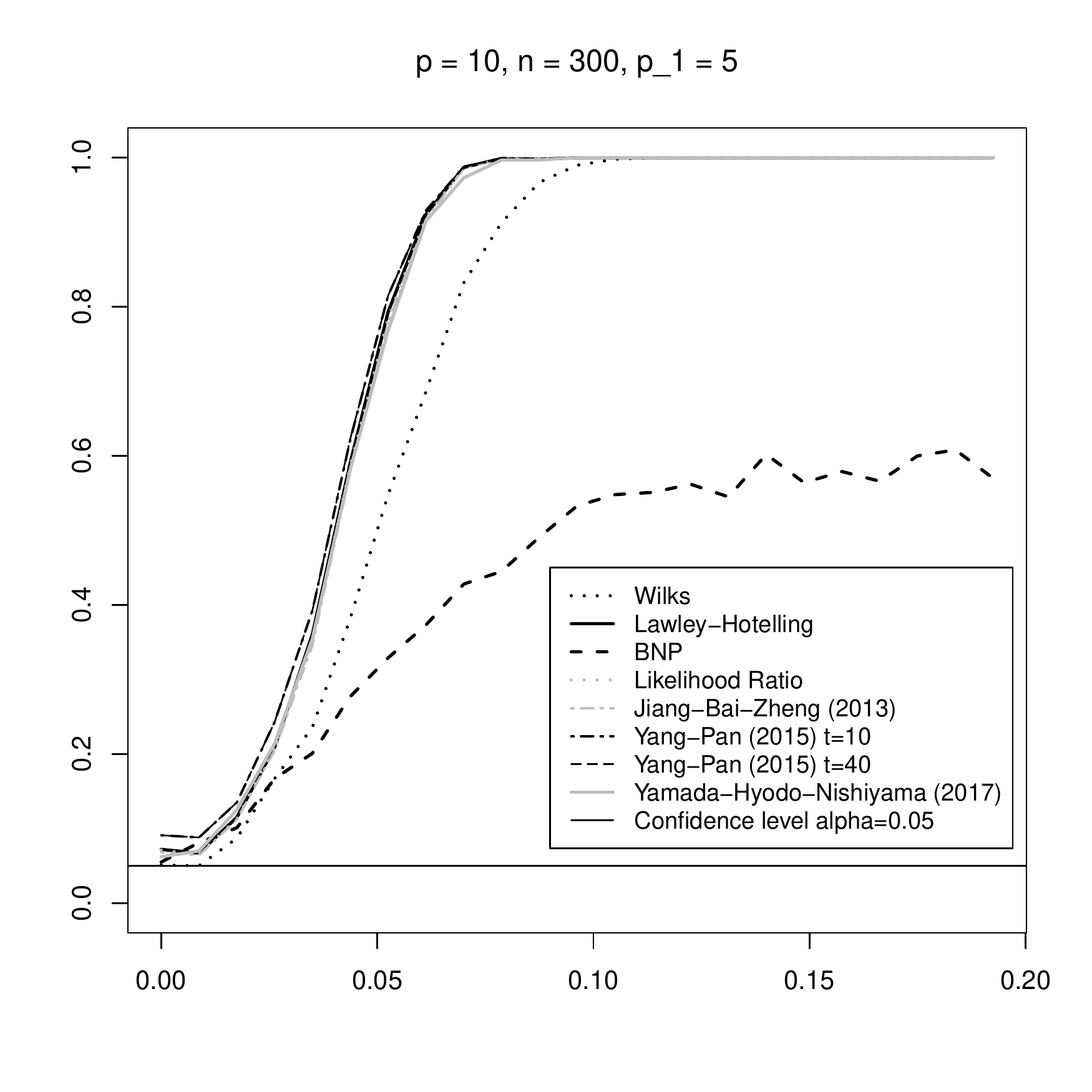} ~~
  \includegraphics[scale=0.28]{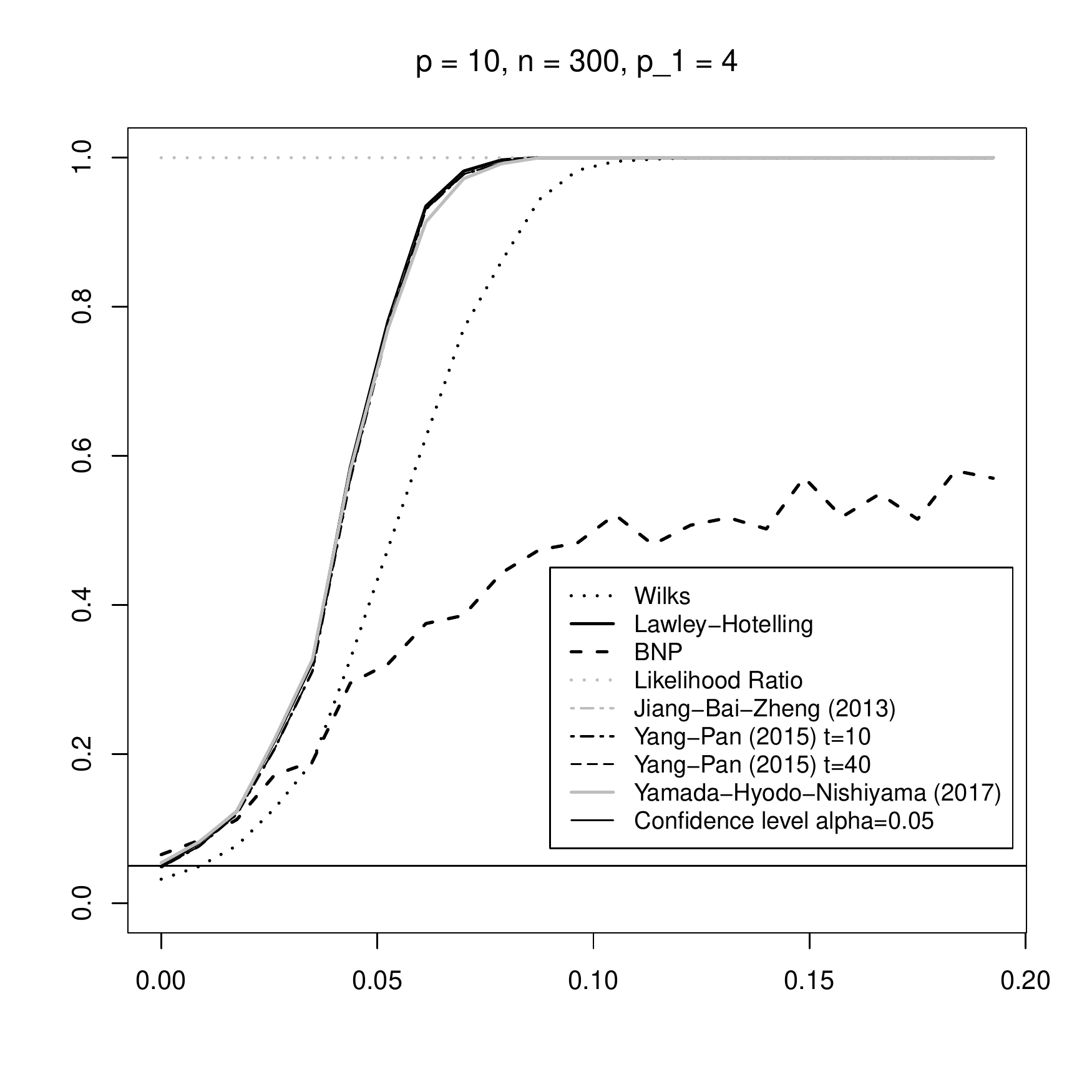}\\

\vspace{-.5cm}
  \caption{\it Empirical power  of  different tests for block diagonality   for  sample size $n=300$,  dimension $p=10$ and  various values of $p_1=8, 5, 4$ as a function of the correlation coefficient $\rho=\frac{\sigma_{12}}{\sigma}$ in $[0, 0.2]$.
  \label{fig10}}
\end{figure}

\begin{figure}[p]
  \centering
  \includegraphics[scale=0.28]{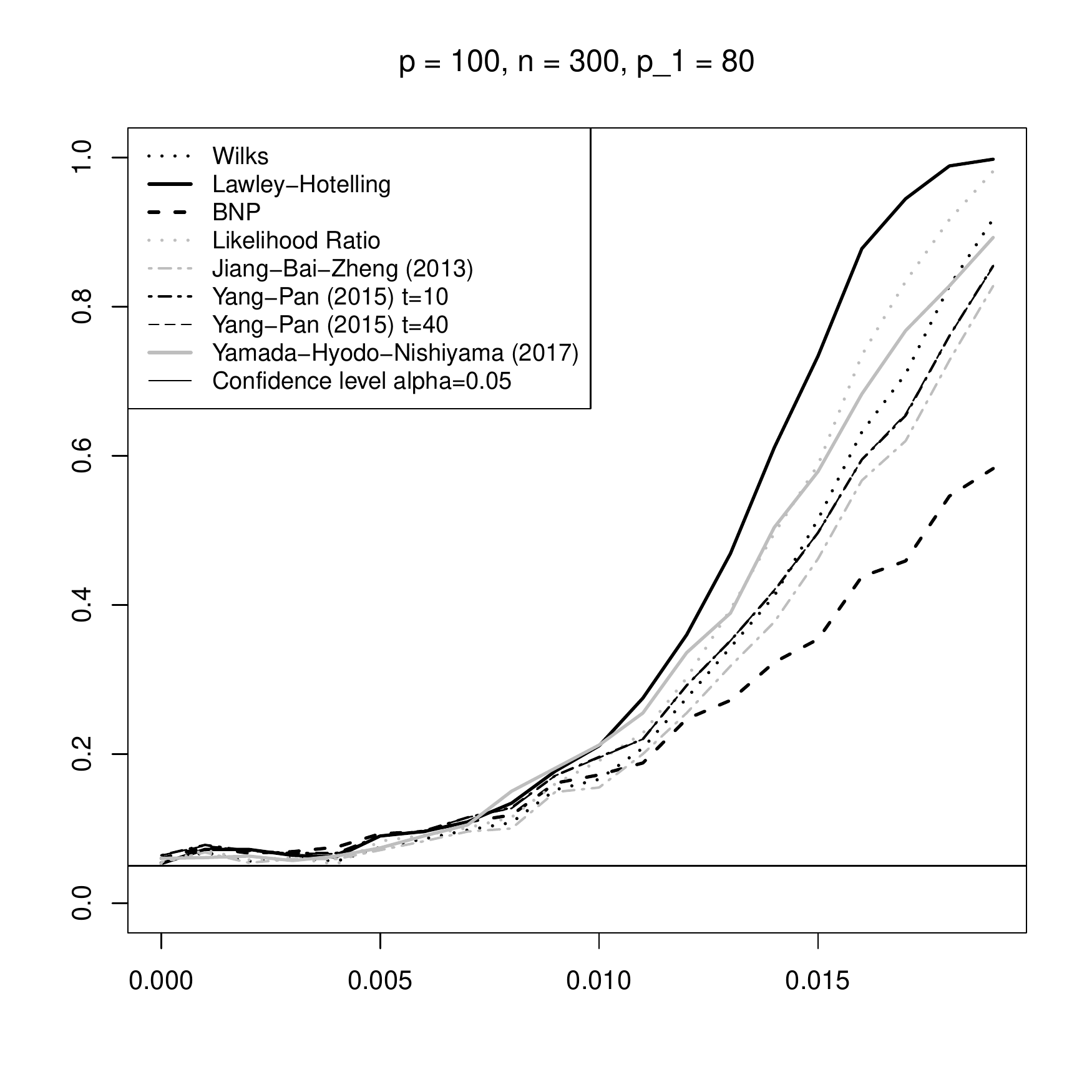} ~~
  \includegraphics[scale=0.28]{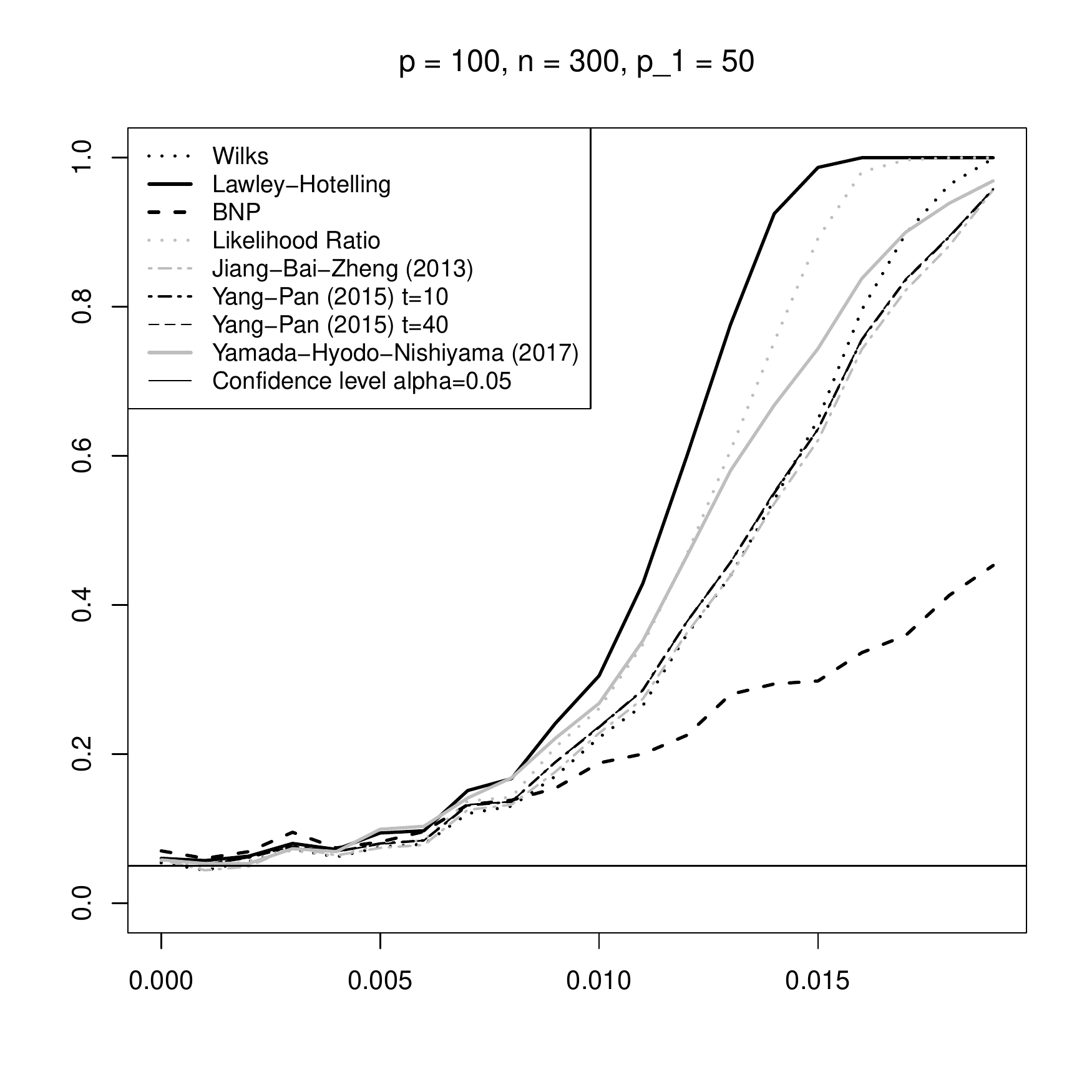} ~~
  \includegraphics[scale=0.28]{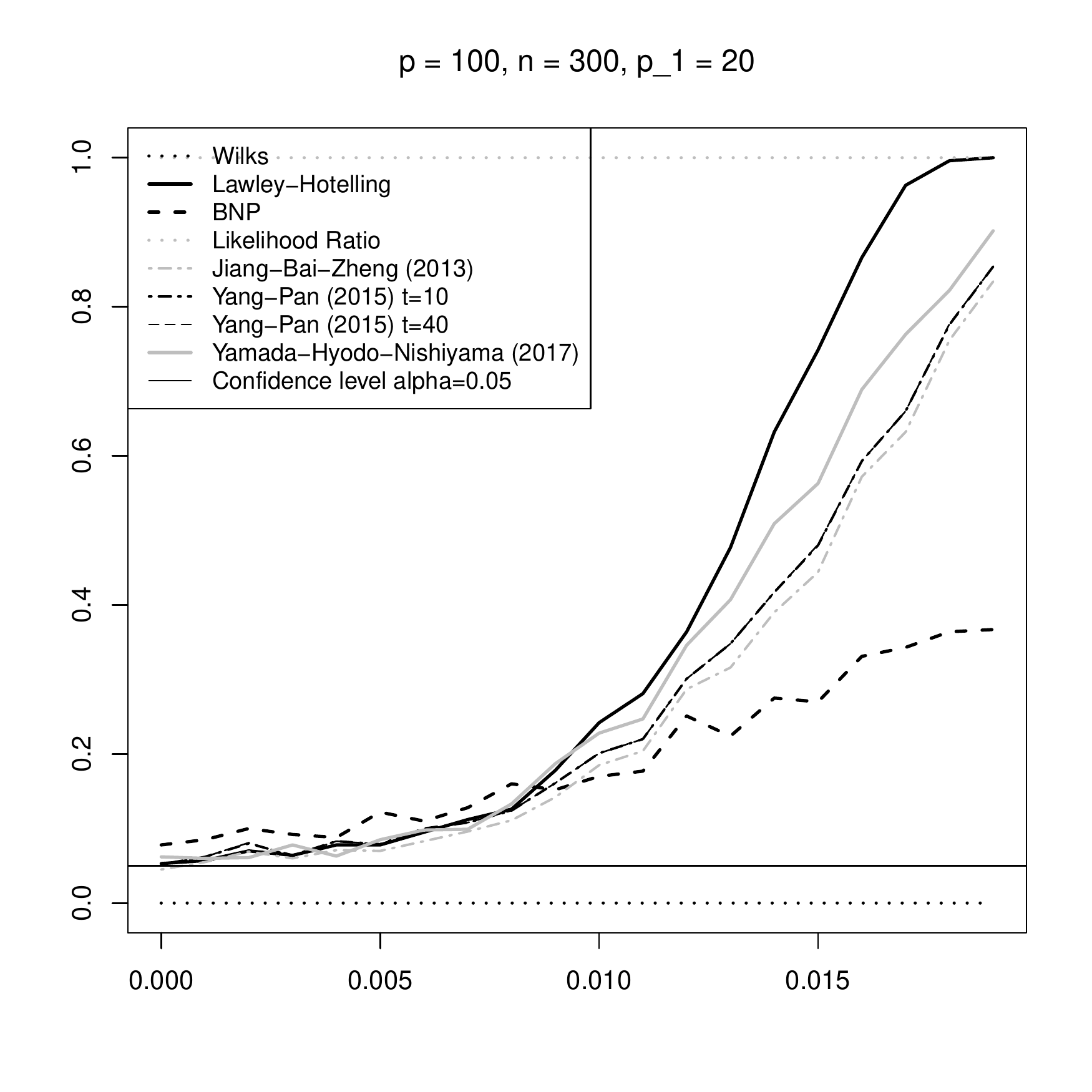}\\

\vspace{-.5cm}
  \caption{\it Empirical power  of  different tests for block diagonality   for  sample size $n=300$,  dimension $p=100$ and  various values of $p_1=80, 50, 20$ as a function of the correlation coefficient $\rho=\frac{\sigma_{12}}{\sigma}$ in $[0, 0.02]$.
  \label{fig11}}
\end{figure}

\begin{figure}[p]
  \centering
  \includegraphics[scale=0.28]{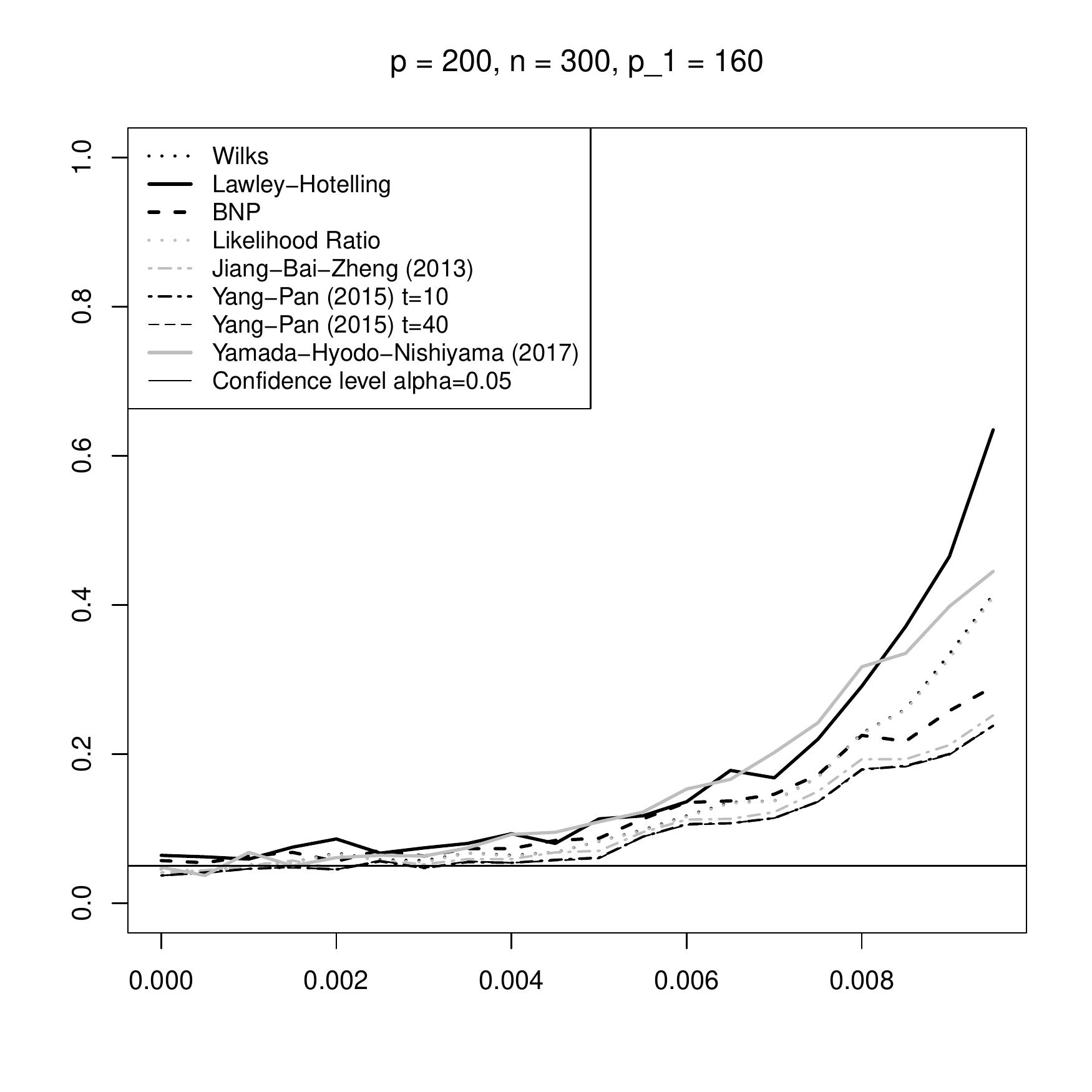} ~~
  \includegraphics[scale=0.28]{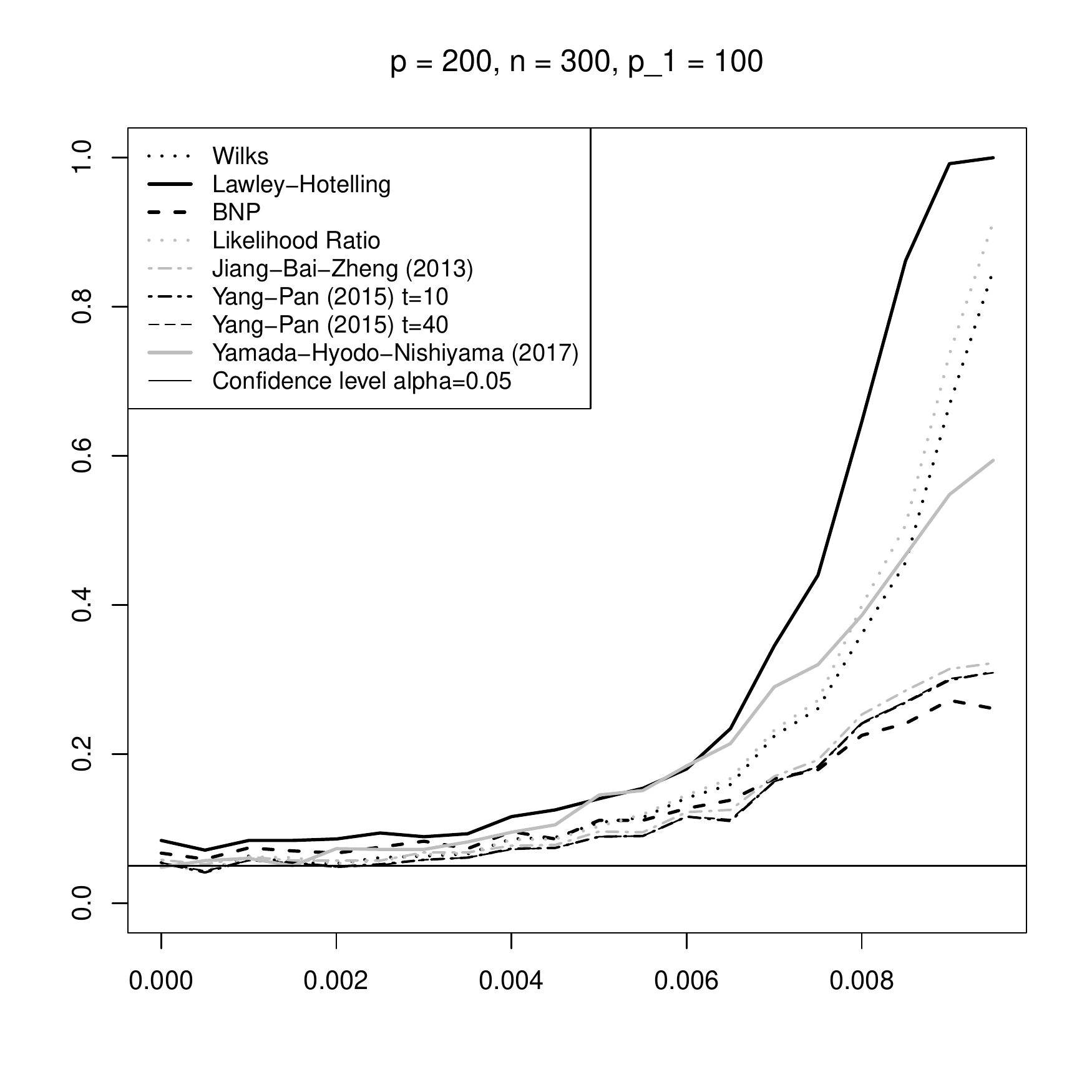} ~~
  \includegraphics[scale=0.28]{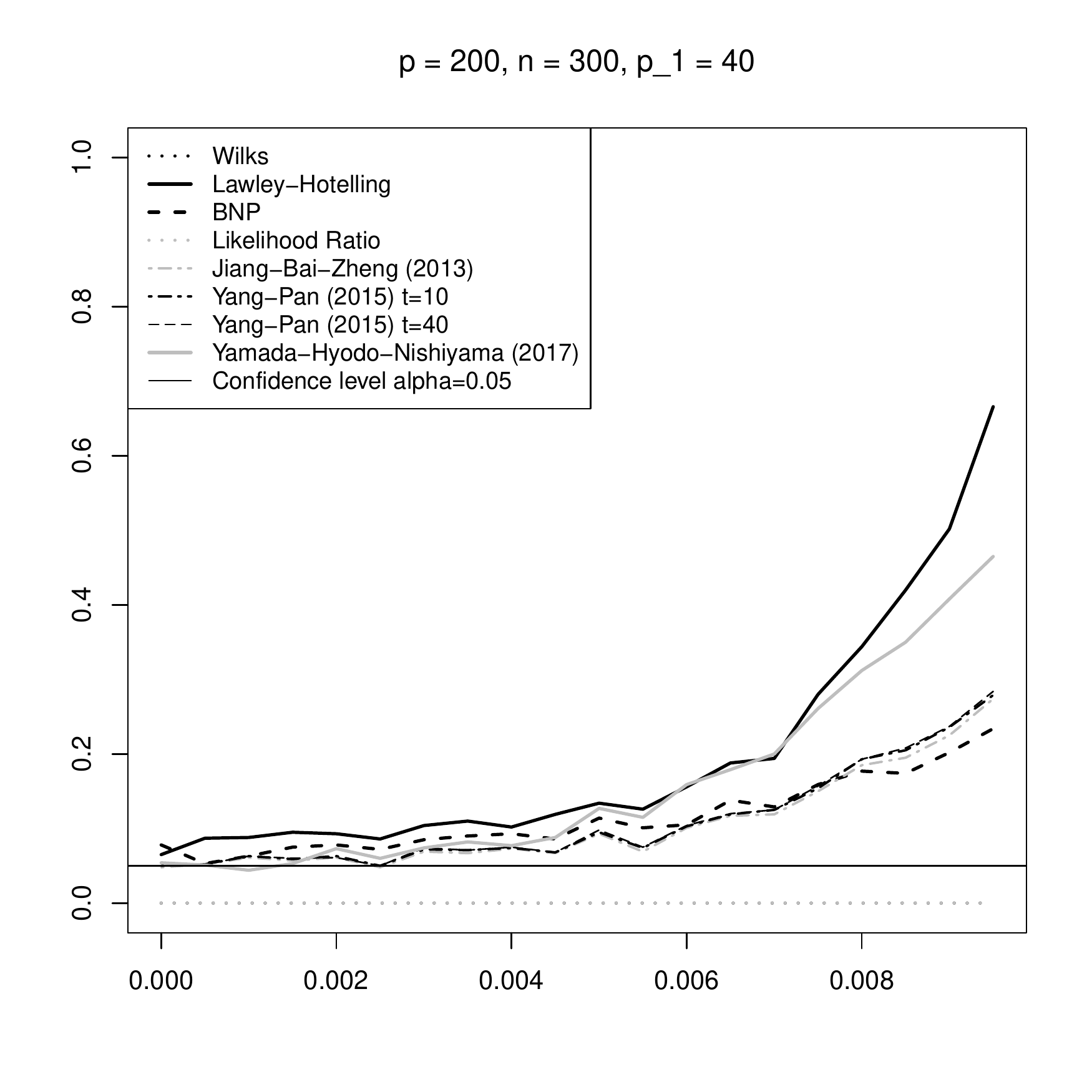}\\
\vspace{-.5cm}
  \caption{\it Empirical power  of  different tests for block diagonality   for  sample size $n=300$,  dimension $p=200$ and  various values of $p_1=160, 100, 40$ as a function of the correlation coefficient $\rho=\frac{\sigma_{12}}{\sigma}$ in $[0, 0.019]$.
  \label{fig12}}
\end{figure}

\begin{figure}[p]
  \centering
  \includegraphics[scale=0.28]{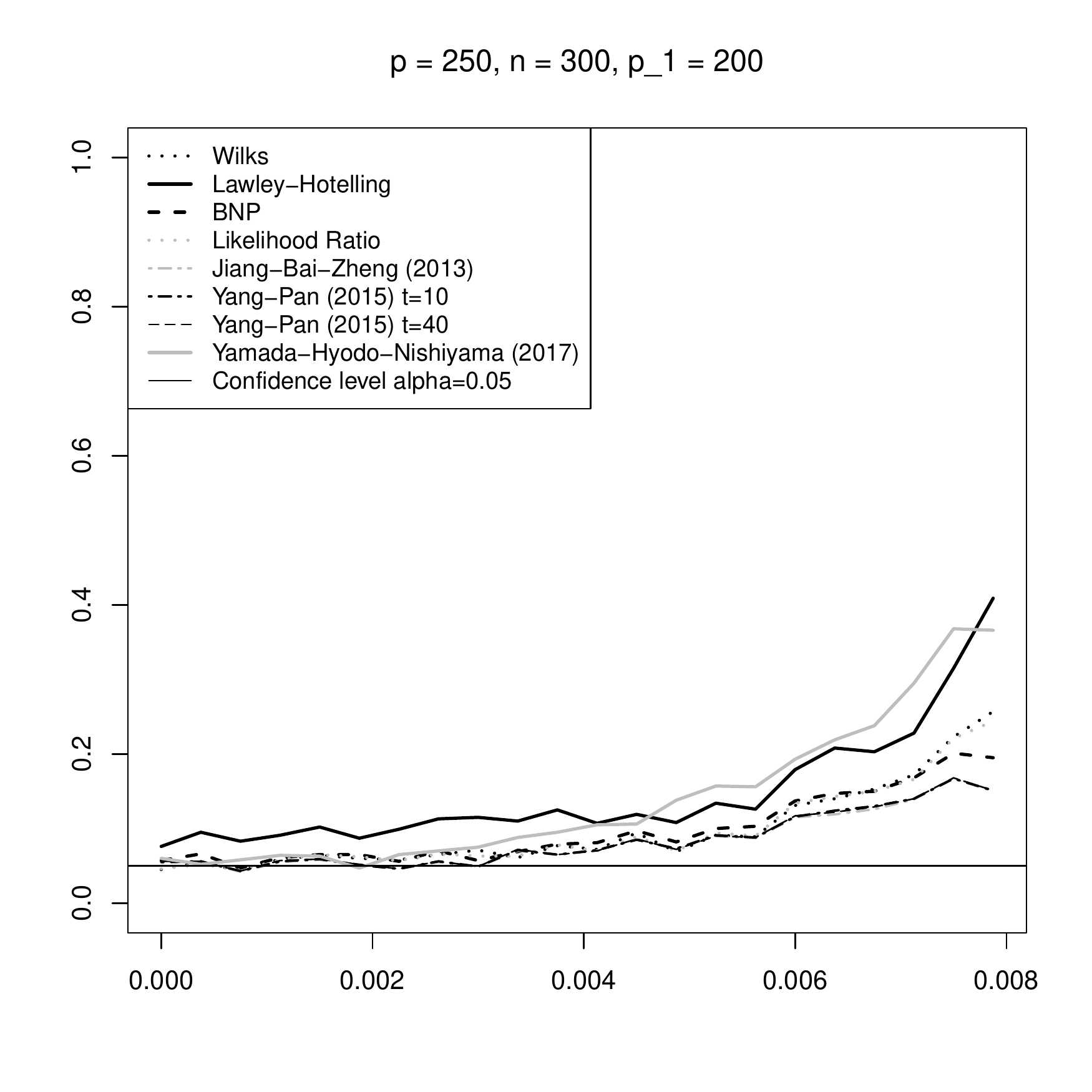} ~~
  \includegraphics[scale=0.28]{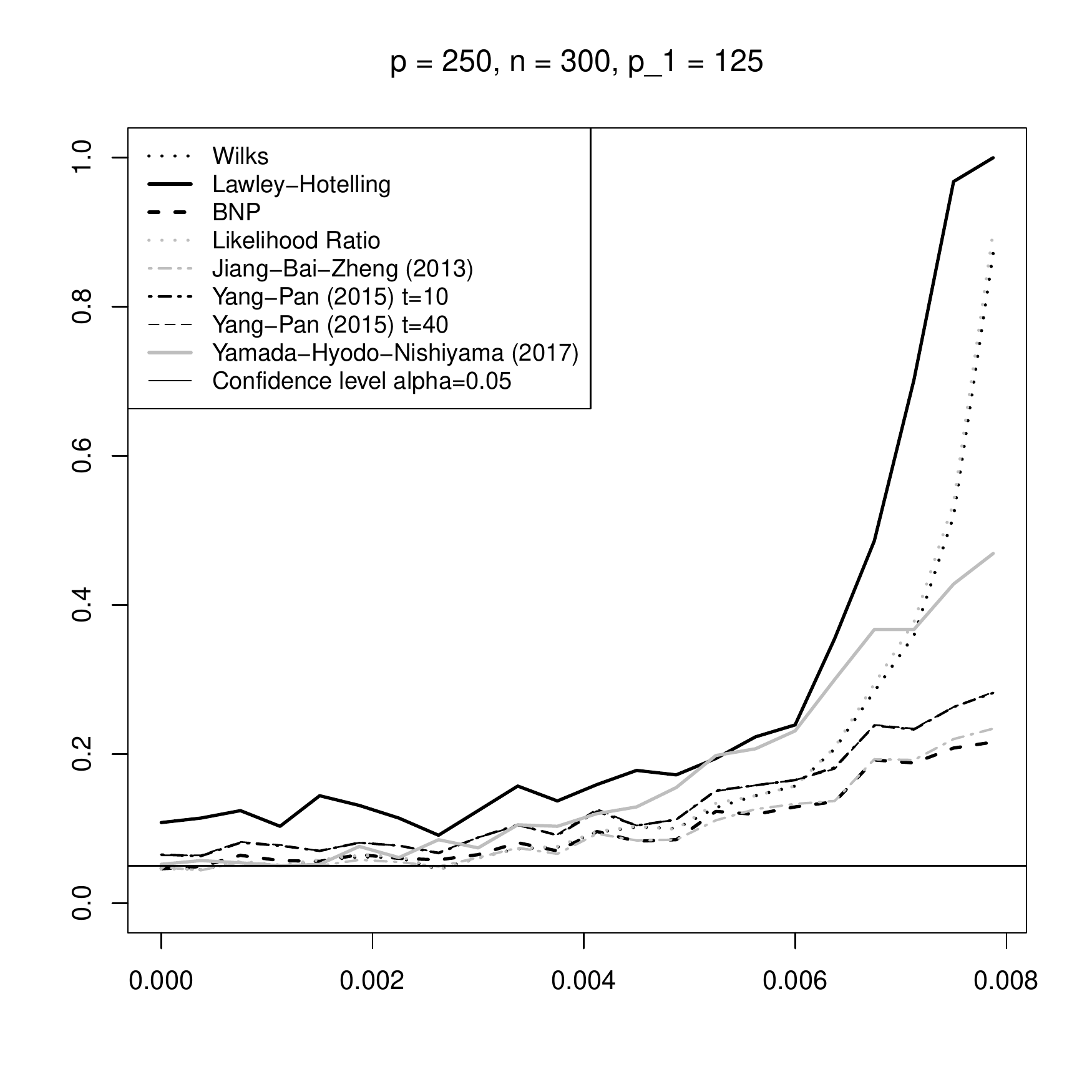} ~~
  \includegraphics[scale=0.28]{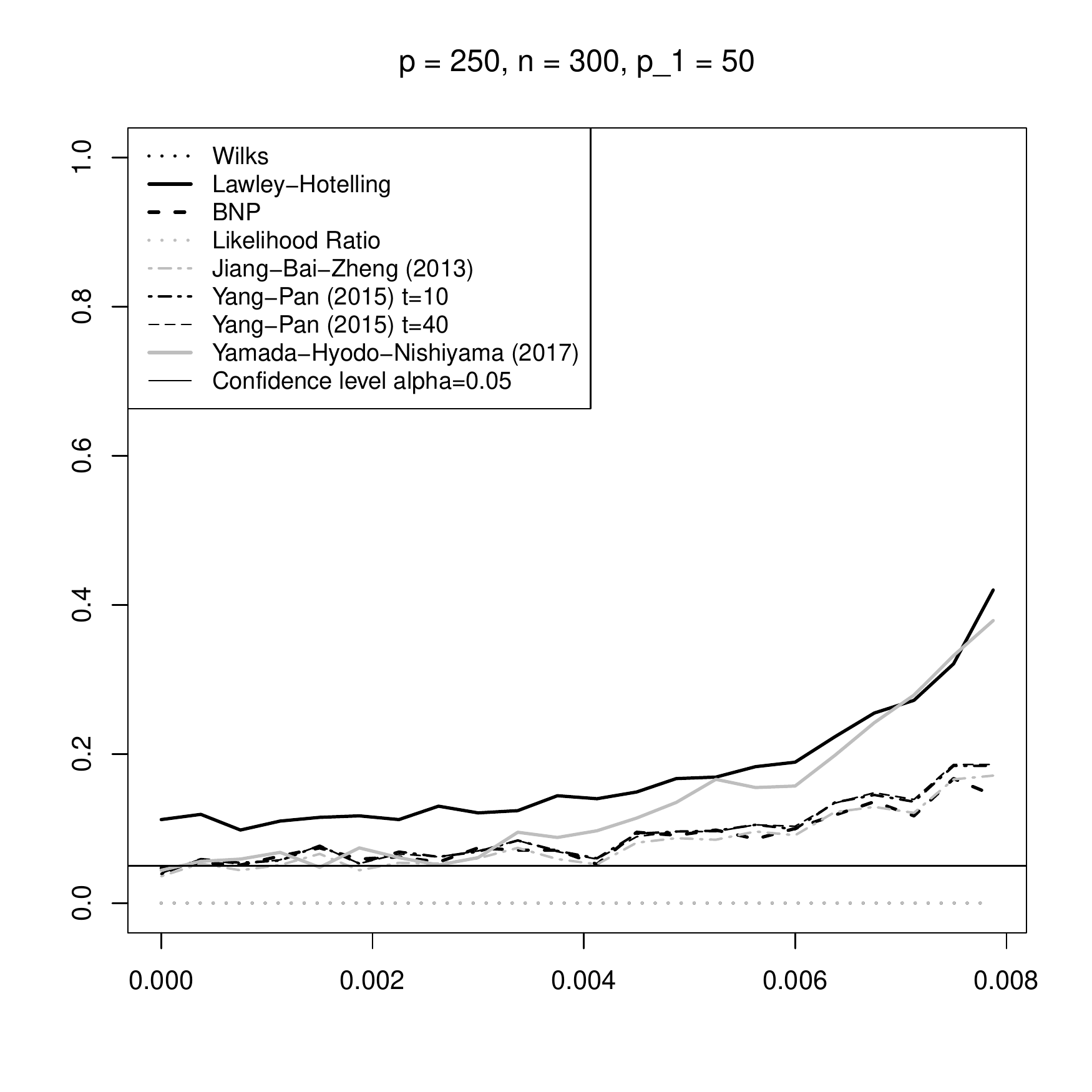}\\

\vspace{-.5cm}
  \caption{\it Empirical power  of  different tests for block diagonality   for  sample size $n=300$,  dimension $p=250$ and  various values of $p_1=200, 125, 50$ as a function of the correlation coefficient $\rho=\frac{\sigma_{12}}{\sigma}$ in $[0, 0.0095]$.
  \label{fig13}}
\end{figure}

\begin{figure}[p]
  \centering
  \includegraphics[scale=0.28]{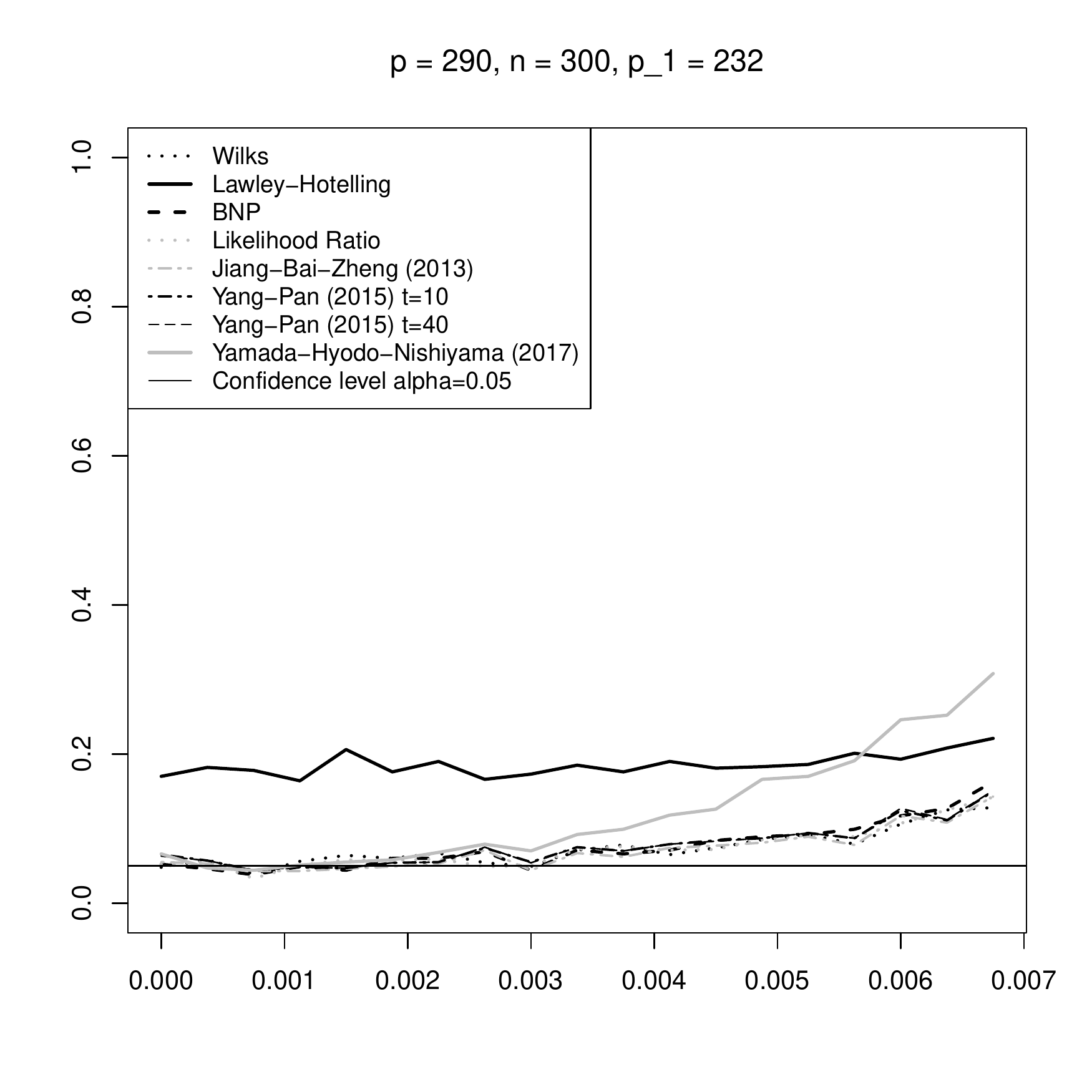} ~~
  \includegraphics[scale=0.28]{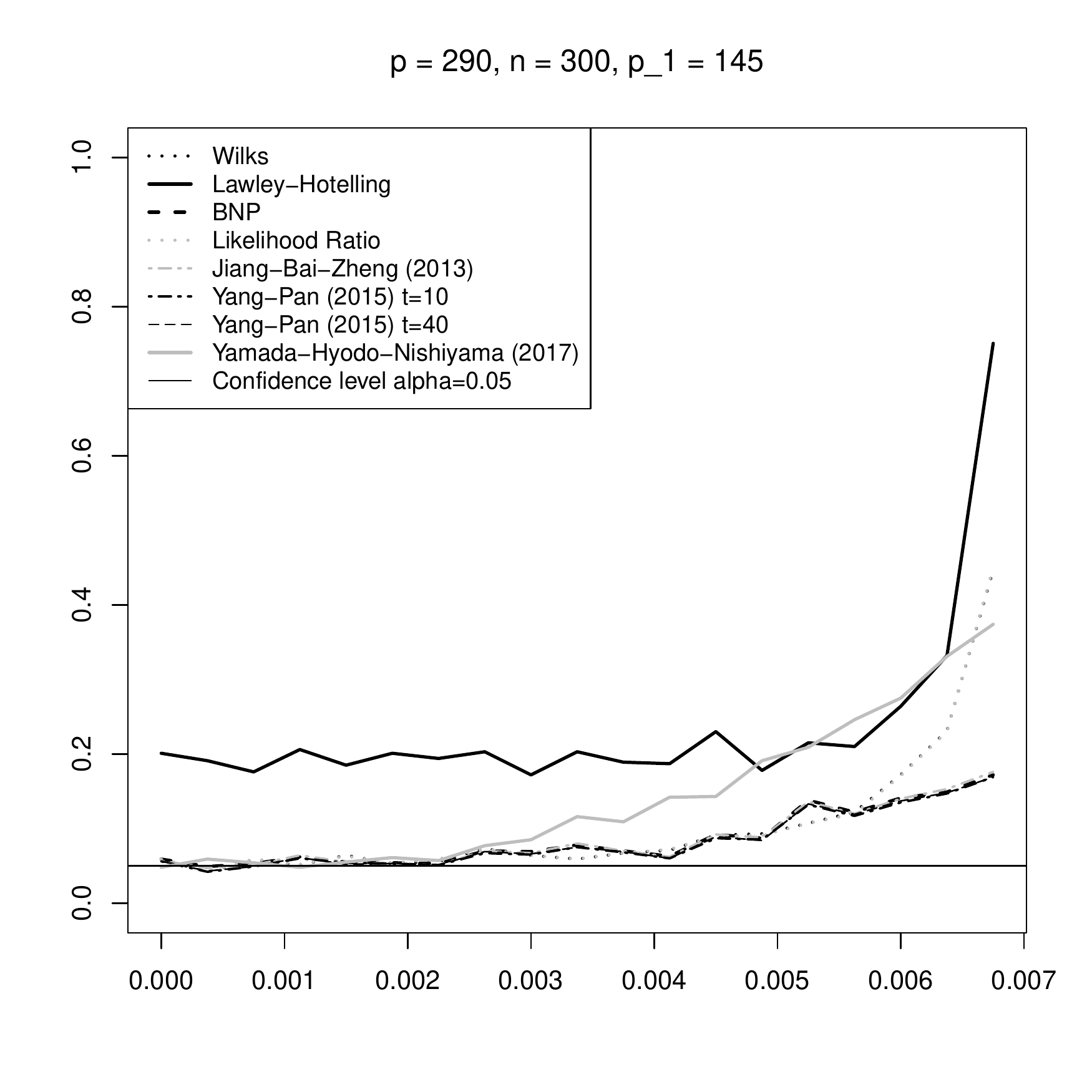} ~~
  \includegraphics[scale=0.28]{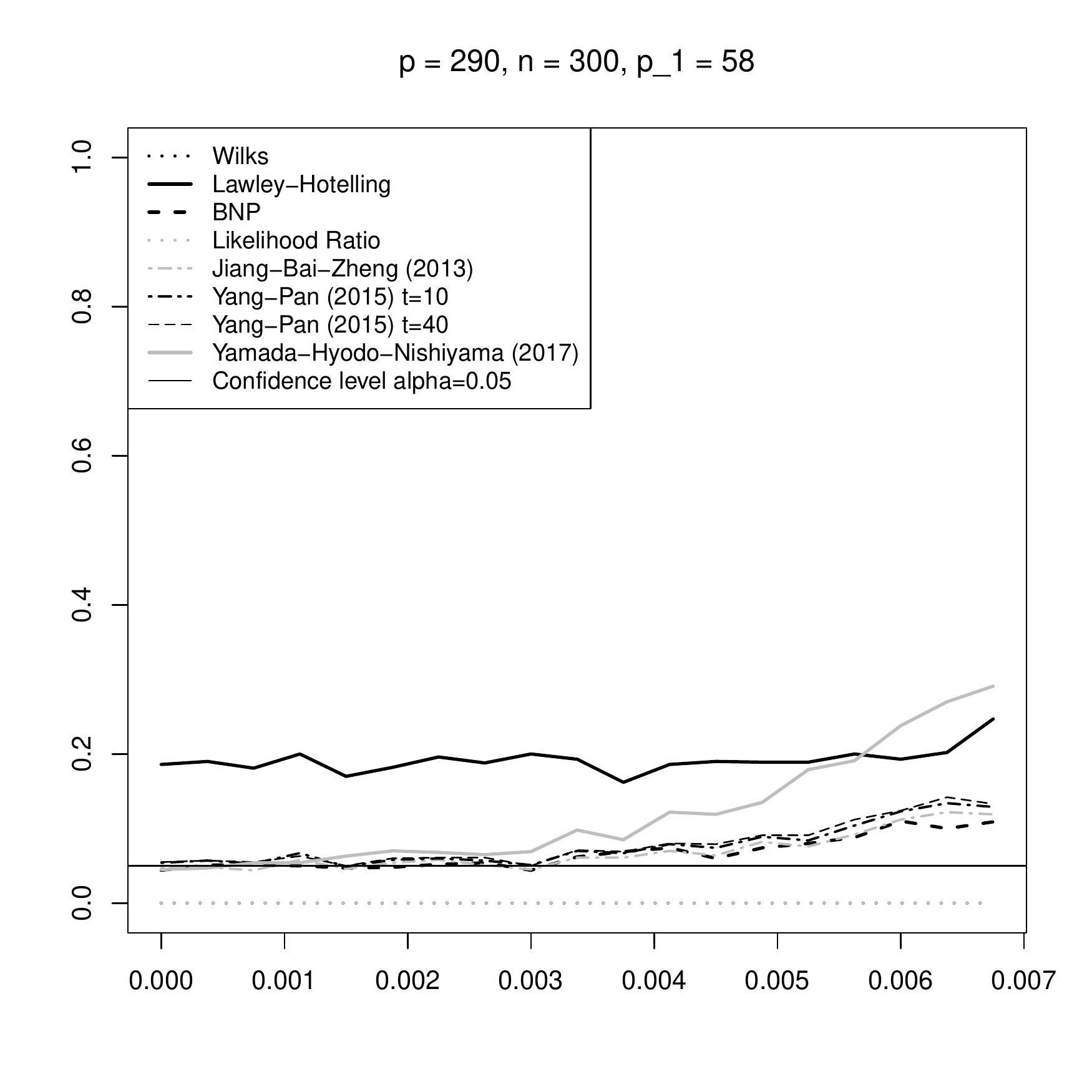}\\

\vspace{-.5cm}
  \caption{\it Empirical power  of  different tests for block diagonality   for  sample size $n=300$,  dimension $p=290$ and  various values of $p_1=232, 145, 58$ as a function of the correlation coefficient $\rho=\frac{\sigma_{12}}{\sigma}$ in $[0, 0.00675]$.
  \label{fig14}}
\end{figure}

{\small
 \bibliography{blockdiag}

\begin{thebibliography}{}

\bibitem[Anderson, 2003]{anderson2003}
Anderson, T.~W. (2003).
\newblock {\em Multivariate Statistical Analysis}.
\newblock John Wiley \& Sons, New York.

\bibitem[Bai et~al., 2009]{baietal2009}
Bai, Z., Jiang, D., Yao, J.-F., and Zheng, S. (2009).
\newblock Corrections to {LRT} on large-dimensional covariance matrix by {RMT}.
\newblock {\em Annals of Statistics}, 37:3822--3840.

\bibitem[Bai and Silverstein, 2004]{baisil2004}
Bai, Z.~D. and Silverstein, J.~W. (2004).
\newblock {CLT} for linear spectral statistics of large dimensional sample
  covariance matrices.
\newblock {\em Annals of Probability}, 32:553--605.

\bibitem[Bai and Silverstein, 2010]{baisilverstein2010}
Bai, Z.~D. and Silverstein, J.~W. (2010).
\newblock {\em Spectral Analysis of Large Dimensional Random Matrices}.
\newblock Springer, New York.

\bibitem[Bickel and Levina, 2008]{bicklevi2008}
Bickel, P.~J. and Levina, E. (2008).
\newblock Regularized estimation of large covariance matrices.
\newblock {\em Annals of Statistics}, 36:199--227.

\bibitem[Birke and Dette, 2005]{birkedette2005}
Birke, M. and Dette, H. (2005).
\newblock A note on testing the covariance matrix for large dimension.
\newblock {\em Statistics {\&} Probability Letters}, 74:281--289.

\bibitem[Bodnar et~al., 2018]{supplement}
Bodnar, T., Dette, H., and Parolya, N. (2018).
\newblock Supplement to ''{T}esting for independence of large dimensional
  vectors''.
\newblock {\em Annals of Statistics}.

\bibitem[Bodnar et~al., 2014]{bodnaretal2014}
Bodnar, T., Gupta, A.~K., and Parolya, N. (2014).
\newblock On the strong convergence of the optimal linear shrinkage estimator
  for large dimensional covariance matrix.
\newblock {\em Journal of Multivariate Analysis}, 132:215--228.

\bibitem[Bodnar et~al., 2016]{bodnaretal2016}
Bodnar, T., Gupta, A.~K., and Parolya, N. (2016).
\newblock Direct shrinkage estimation of large dimensional precision matrix.
\newblock {\em Journal of Multivariate Analysis}, 146:223--236.

\bibitem[Cai et~al., 2011]{cailuiluo2011}
Cai, T., Lui, W., and Luo, X. (2011).
\newblock A constrained {$l_1$} minimization approach to sparse precision
  matrix estimation.
\newblock {\em Journal of the American Statistical Association}, 106:594--607.

\bibitem[Cai et~al., 2013]{cairenzhou2013}
Cai, T., Ren, Z., and Zhou, H.~H. (2013).
\newblock Optimal rates of convergence for estimating {T}oeplitz covariance
  matrices.
\newblock {\em Probability Theory and Related Fields}, 156(1):101--143.

\bibitem[Cai and Shen, 2011]{caishen2011}
Cai, T. and Shen, X. (2011).
\newblock {\em High-Dimensional Data Analysis, Volume 2, Frontiers of
  Statistics}.
\newblock World Scientific, Singapore.

\bibitem[Cai and Zhou, 2012]{caizhou2012}
Cai, T. and Zhou, H. (2012).
\newblock Minimax estimation of large covariance matrices under {$l_1$} norm.
\newblock {\em Statistica Sinica}, 22:1319--1378.

\bibitem[Chen et~al., 2010]{chenzhangzhong2010}
Chen, S.~X., Zhang, L.-X., and Zhong, P.-S. (2010).
\newblock Testing high dimensional covariance matrices.
\newblock {\em Journal of the American Statistical Association}, 105:810--819.

\bibitem[Devijver and Gallopin, 2016]{devigall2016}
Devijver, E. and Gallopin, M. (2016).
\newblock Block-diagonal covariance selection for high-dimensional {G}aussian
  graphical models.
\newblock {\em Journal of the American Statistical Association, to appear}.

\bibitem[Dozier and Silverstein, 2007]{dozsil2007}
Dozier, R.~B. and Silverstein, J.~W. (2007).
\newblock On the empirical distribution of eigenvalues of large dimensional
  information-plus-noise-type matrices.
\newblock {\em Journal of Multivariate Analysis}, 98:678--694.

\bibitem[Fan and Li, 2006]{Fan2006}
Fan, J. and Li, R. (2006).
\newblock Statistical challenges with high dimensionality: feature selection in
  knowledge discovery.

\bibitem[Fisher, 1939]{Fisher1939}
Fisher, R.~A. (1939).
\newblock The sampling distribution of some statistics obtained from non-linear
  equations.
\newblock {\em Annals of Eugenics}, 9(3):238--249.

\bibitem[Fisher et~al., 2010]{fishsungall2010}
Fisher, T.~J., Sun, X., and Gallagher, C.~M. (2010).
\newblock A new test for sphericity of the covariance matrix for high
  dimensional data.
\newblock {\em Journal of Multivariate Analysis}, 101:2554--2570.

\bibitem[Fujikoshi et~al., 2004]{fujhimwak2004}
Fujikoshi, Y., Himeno, T., and Wakaki, H. (2004).
\newblock Asymptotic results of a high dimensional {MANOVA} test and power
  comparison when the dimension is large compared to the sample size.
\newblock {\em Journal of the Japan Statistical Society}, 34(1):19--26.

\bibitem[Gupta and Xu, 2006]{guptaxu2006}
Gupta, A.~K. and Xu, J. (2006).
\newblock On some tests of the covariance matrix under general conditions.
\newblock {\em Annals of the Institute of Statistical Mathematics},
  58:101--114.

\bibitem[Hachem et~al., 2012]{hachem2012}
Hachem, W., Kharouf, J., Najim, J., and Silverstein, J.~W. (2012).
\newblock A {CLT} for information-theoretic statistics of non-centered {G}ram
  random matrices.
\newblock {\em Random Matrices: Theory and Applications}, 1(2):1150010.

\bibitem[Hyodo et~al., 2015]{hyodoetal2015}
Hyodo, M., Shutoh, N., Nishiyama, T., and Pavlenko, T. (2015).
\newblock Testing block-diagonal covariance structure for high-dimensional
  data.
\newblock {\em Statistica Neerlandica}, 69(4):460--482.

\bibitem[Jiang et~al., 2013]{jiangetal2013}
Jiang, D., Bai, Z., and Zheng, S. (2013).
\newblock Testing the independence of sets of large-dimensional variables.
\newblock {\em Sci. China Math.}, 56(1):135--147.

\bibitem[Jiang and Yang, 2013]{jiangyang2013}
Jiang, T. and Yang, F. (2013).
\newblock Central limit theorems for classical likelihood ratio tests for
  high-dimensional normal distributions.
\newblock {\em Annals of Statistics}, 41:2029--2074.

\bibitem[John, 1971]{john1971}
John, S. (1971).
\newblock Some optimal multivariate tests.
\newblock {\em Biometrika}, 58:123--127.

\bibitem[Johnstone, 2006]{Johnstone}
Johnstone, I. (2006).
\newblock High dimensional statistical inference and random matrices.

\bibitem[Johnstone, 2001]{johnstone2001}
Johnstone, I.~M. (2001).
\newblock On the distribution of the largest eigenvalue in principal components
  analysis.
\newblock {\em Annals of Statistics}, 29:295--327.

\bibitem[Johnstone, 2008]{johnstone2008}
Johnstone, I.~M. (2008).
\newblock Multivariate analysis and {J}acobi ensembles: {L}argest eigenvalue,
  {T}racy-{W}idom limits and rates of convergence.
\newblock {\em Annals of Statistics}, 36:2638--2716.

\bibitem[Kakizawa and Iwashita, 2008]{kakiiwas2008}
Kakizawa, Y. and Iwashita, T. (2008).
\newblock A comparison of higher-order local powers of a class of one-way
  {MANOVA} tests under general distributions.
\newblock {\em Journal of Multivariate Analysis}, 99(6):1128--1153.

\bibitem[Ledoit and Wolf, 2002]{ledoitwolf2002}
Ledoit, O. and Wolf, M. (2002).
\newblock Some hypothesis tests for the covariance matrix when the dimension is
  large compared to the sample size.
\newblock {\em Annals of Statistics}, 30:1081--1102.

\bibitem[Ledoit and Wolf, 2012]{ledoitwolf2012}
Ledoit, O. and Wolf, M. (2012).
\newblock Nonlinear shrinkage estimation of large-dimensional covariance
  matrices.
\newblock {\em Annals of Statistics}, 40:1024--1060.

\bibitem[Markowitz, 1952]{markowitz1952}
Markowitz, H. (1952).
\newblock Portfolio selection.
\newblock {\em Journal of Finance}, 7(1):77--91.

\bibitem[Mauchly, 1940]{mauchly1940}
Mauchly, J.~W. (1940).
\newblock Significance test for sphericity of a normal {N}-variate
  distribution.
\newblock {\em Annals of Mathematical Statistics}, 11:204--209.

\bibitem[Muirhead, 1982]{muirhead1982}
Muirhead, R.~J. (1982).
\newblock {\em Aspects of Multivariate Statistical Theory}.
\newblock John Wiley \& Sons, New York.

\bibitem[Nagao, 1973]{nagao1973}
Nagao, H. (1973).
\newblock On some test criteria for covariance matrix.
\newblock {\em Annals of Statistics}, 1:700--709.

\bibitem[Pillai and Jayachandran, 1967]{pillai1967}
Pillai, K. C.~S. and Jayachandran, K. (1967).
\newblock Power comparisons of tests of two multivariate hypotheses based on
  four criteria.
\newblock {\em Biometrika}, 54(1/2):195--210.

\bibitem[Schott, 2007]{schott2007}
Schott, J.~R. (2007).
\newblock Some high-dimensional tests for a one-way {MANOVA}.
\newblock {\em Journal of Multivariate Analysis}, 98(9):1825--1839.

\bibitem[Wang et~al., 2015]{wangpantongzhu2015}
Wang, C., Pan, G.~M., Tong, T.~J., and Zhu, L. (2015).
\newblock Shrinkage estimator of large dimensional precision matrix using
  random matrix theory.
\newblock {\em Statistica Sinica}, 25:993--1008.

\bibitem[Yamada et~al., 2017]{yamadaetal2017}
Yamada, Y., Hyodo, M., Shutoh, N., and Nishiyama, T. (2017).
\newblock Testing block-diagonal covariance structure for high-dimensional data
  under non-normality.
\newblock {\em Journal of Multivariate Analysis}, 155:305--316.

\bibitem[Yang and Pan, 2015]{yang2015}
Yang, Y. and Pan, G. (2015).
\newblock Independence test for high dimensional data based on regularized
  canonical correlation coefficients.
\newblock {\em Ann. Statist.}, 43(2):467--500.

\bibitem[Yao, 2013]{yao2013}
Yao, J. (2013).
\newblock {Estimation and fluctuations of functionals of large random
  matrices.}
\newblock {Telecom ParisTech, tel-00909521v1}.

\bibitem[Yao et~al., 2015]{yaobaizheng2015}
Yao, J., Bai, Z., and Zheng, S. (2015).
\newblock {\em Large Sample Covariance Matrices and High-Dimensional Data
  Analysis (No. 39)}.
\newblock Cambridge University Press, New York.

\bibitem[Zheng, 2012]{zheng2012}
Zheng, S. (2012).
\newblock Central limit theorems for linear spectral statistics of large
  dimensional {F}-matrices.
\newblock {\em Annales de l'Institut Henri Poincar{\'{e}}, Probabilit{\'{e}s et
  Statistiques}}, 48(2):444--476.

\bibitem[Zheng et~al., 2017]{zheng2017}
Zheng, S., Bai, Z., and Yao, J. (2017).
\newblock Clt for eigenvalue statistics of large-dimensional general fisher
  matrices with applications.
\newblock {\em Bernoulli}, 23(2):1130--1178.

\bibitem[Zheng et~al., 2015a]{zhengbaiyao2015}
Zheng, Z., Bai, Z.~D., and Yao, J. (2015a).
\newblock {CLT} for linear spectral statistics of a rescaled sample precision
  matrix.
\newblock {\em Random Matrices: Theory and Applications}, 04(04):1550014.

\bibitem[Zheng et~al., 2015b]{zhengbaiyao2015b}
Zheng, Z., Bai, Z.~D., and Yao, J. (2015b).
\newblock Substitution principle for {CLT} of linear spectral statistics of
  high-dimensional sample covariance matrices with applications to hypothesis
  testing.
\newblock {\em Annals of Statistics}, 43:546--591.

\end{thebibliography}
}

\end{document}